\documentclass{siamltex}

\usepackage{amsmath,amsfonts,graphicx,fullpage,pdflscape,color,multirow}
\usepackage{subfig,booktabs,array}
\usepackage[pdfborder={0 0 0 0}]{hyperref}

\newcolumntype{C}[1]{>{\centering\let\newline\\\arraybackslash\hspace{0pt}}m{#1}}

\newcommand{\vect}[1]{\mathbf{#1}}

\newcommand{\clawpack}{{\sc clawpack}}
\newcommand{\amrclaw}{{\sc amrclaw}}

\newcommand{\refcartpaper}{\cite{gil-ou-rjl:poro-2d-cartesian}}

\newcommand{\beq}{\begin{eqnarray}}
\newcommand{\eeq}{\end{eqnarray}}
\newcommand{\beqs}{\begin{eqnarray*}}
\newcommand{\eeqs}{\end{eqnarray*}}
\newcommand{\bC}{\mathbf{C}}

\newcommand{\btau}{\mbox{\boldmath $\tau$}}
\newcommand{\bA}{\mbox{$\mathbf{A}$}}
\newcommand{\bB}{\mbox{$\mathbf{B}$}}
\newcommand{\bD}{\mbox{$\mathbf{D}$}}
\newcommand{\bE}{\mbox{$\mathbf{E}$}}

\newcommand{\bbeta}{\mbox{\boldmath $\beta$}}

\newcommand{\bQ}{\mbox{$\mathbf{Q}$}}

\newcommand{\B}{\mathbf}
\newcommand{\warn}[1]{}

\newcommand{\archivelink}{\url{http://dx.doi.org/10.6084/m9.figshare.701483}}

\author{Grady\ I.\ Lemoine\footnotemark[1] \and M.\
Yvonne\ Ou\footnotemark[2]}

\begin{document}

\title{Finite Volume Modeling of Poroelastic-Fluid Wave Propagation with Mapped Grids}
\maketitle

\renewcommand{\thefootnote}{\fnsymbol{footnote}}

\footnotetext[1]{Corresponding author.  Department of Applied
  Mathematics, University of Washington, Seattle, WA 98195
  (gl@uw.edu).  This author's work was supported in part by NIH grant
  5R01AR53652-2 and NSF grants DMS-0914942 and DMS-1216732.}

\footnotetext[2]{Department of Mathematical Sciences, University of
  Delaware, Newark, DE 19716.  This author's work was supported in
  part by NSF-DMS Mathematical Biology Grant 0920852 and NIH CBER
  PILOT Grant 322 159.}

\begin{abstract}
In this work we develop a high-resolution mapped-grid finite volume
method code to model wave propagation in two dimensions in systems of
multiple orthotropic poroelastic media and/or fluids, with curved
interfaces between different media.  We use a unified formulation to
simplify modeling of the various interface conditions --- open pores,
imperfect hydraulic contact, or sealed pores --- that may exist
between such media.  Our numerical code is based on the \clawpack{}
framework, but in order to obtain correct results at a material
interface we use a modified transverse Riemann solution scheme, and at
such interfaces are forced to drop the second-order correction term
typical of high-resolution finite volume methods.  We verify our code
against analytical solutions for reflection and transmission of waves
at a material interface, and for scattering of an acoustic wave train
around an isotropic poroelastic cylinder.  For reflection and
transmission at a flat interface, we achieve second-order convergence
in the 1-norm, and first-order in the max-norm; for the cylindrical
scatterer, the highly distorted grid mapping degrades performance but
we still achieve convergence at a reduced rate.  We also simulate an
acoustic pulse striking a simplified model of a human femur bone, as
an example of the capabilities of the code.  To aid in
reproducibility, at the web site
\archivelink{} we provide all of
the code used to generate the results here.
\end{abstract}

\begin{keywords}
poroelastic, wave propagation, finite-volume, high-resolution,
operator splitting, mapped grid, transverse solve, interface
condition, cylindrical scatterer
\end{keywords}

\begin{AMS}
  65M08, 74S10, 74F10, 74J10, 74L05, 74L15, 86-08
\end{AMS}

\section{Introduction}
\label{sec:intro}

Poroelasticity theory was developed by Maurice A. Biot to model the
mechanics of a fluid-saturated porous medium.  It models the medium in
a homogenized fashion, with solid portion treated with linear
elasticity, and the fluid with linearized compressible fluid dynamics
combined with Darcy's law to relate its pressure gradient to its flow
rate.  Biot's work is summarized in his 1956 and 1962
papers~\cite{biot:56-1, biot:56-2, biot:62}, and Carcione also
provides an excellent discussion of poroelasticity in chapter 7 of his
book~\cite{carcione:wave-book}.  While it was originally developed to
model fluid-saturated rock and soil, Biot theory has also found
applications in modeling of {\em in vivo} bone~\cite{cowin:bone-poro,
  cowin-cardoso:fabric-2010, gilbert-guyenne-ou:bone-2012} and
underwater acoustics with a porous sea floor~\cite{bgwx:2004,
  gilbert-lin:1997, gilbert-ou:2003}.

Biot theory predicts three different families of propagating waves
within a poroelastic medium.  In order of decreasing speed, these are:
fast P waves, where the fluid and solid parts of the medium move
roughly parallel to the propagation direction --- exactly parallel for
an isotropic medium --- and are typically in phase with each other; S
waves, where the motion of the medium is transverse to the propagation
direction; and slow P waves, where the motion is again roughly
parallel to the wavevector but the fluid and solid typically move 180
degrees out of phase, so that the fluid is leaving a region as the
solid contracts, and vice versa.  Because Biot theory includes viscous
drag between the walls of the pore structure and the fluid, all three
families of waves dissipate as they propagate through the medium; the
slow P wave typically involves much more relative motion between the
solid and fluid than the other two waves, so it is much more strongly
damped.

Researchers have used a variety of numerical methods to model
poroelasticity.  Carcione, Morency, and Santos provide a detailed
chronicle~\cite{carcione-morency-santos:comp-poro-rev}, and we provide
a brief review in the introduction to our previous
paper~\refcartpaper{}.  In this paper we wish to particularly draw
attention to the work of Chiavassa and
Lombard~\cite{chiavassa-lombard:acou-poro,
  chiavassa-lombard:poro-fdm-2011}, who model problems very similar to
the ones we analyze.  They use a fourth-order finite difference method
to solve a first-order velocity-stress form of Biot's equations,
coupled with operator splitting to handle the viscous dissipation
term, and employ an immersed interface approach at the boundary
between the fluid and solid parts of the domain.

This paper extends our previous work in finite volume modeling of a
first-order velocity-stress formulation of Biot poroelasticity
theory~\refcartpaper{} to include non-rectilinear mapped grids and
bounded interfaces between poroelastic media and fluids.  Both of these
capabilities are important for applications such as ocean bottom
acoustics, where the bottom bathymetry may be uneven and waves
propagate from the water into the ocean floor and back again, and
applications to wave propagation in bone, since bones are rarely
rectangular in shape and can benefit from a fitted grid, and wave
propagation between bone and the surrounding fluid or soft tissue is
again important.  We also implement explicit interface conditions
between distinct poroelastic media, such as those proposed by Deresiewicz
and Skalak~\cite{deresiewicz-skalak:poro-uniqueness-1963}.  We
continue to employ the \clawpack{} finite volume method
package~\cite{claw.org}, which substantially reduces the programming
time and effort required for this work; \clawpack{} also offers
built-in operator splitting, as well as
adaptive mesh refinement if desired~\cite{mjb-rjl:amrclaw}.  We
develop a technique to solve the Riemann problem efficiently for grid
interfaces oriented at arbitrary angles in anisotropic poroelastic
media, as well as modifications to the classical high-resolution
finite volume method that are necessary to obtain qualitatively
correct solutions to problems with fluid-poroelastic interfaces.  We
then proceed to verify our code against known analytical solutions
involving interfaces in both rectilinear and curved geometries, and
analyze a sample biological problem of an acoustic wave striking a
simplfied model of a human femur bone.

\section{Governing PDEs and interface conditions}
\label{sec:poro-review}

A poroelastic medium consists of a porous solid {\em skeleton} or {\em
  matrix}, saturated with a fluid.  We use
essentially the same first-order formulation of Biot poroelasticity
theory in two dimensions as in our previous work~\refcartpaper{}, with
the order of the variables in the state vector permuted in order to
emphasize the block structure of the system.  The state variables at a
point are the total stress tensor $\btau$ (the combined stress from
the solid skeleton and the fluid pressure), the fluid pressure $p$,
the solid velocity $\vect{v}$, and the fluid volumetric flow velocity
$\vect{q}$ relative to the solid, which is equal to the fluid velocity
relative to the solid divided by the porosity $\phi$.  This section gives only
a quick summary of the system of PDEs and ancillary functionals;
readers desiring more detail are encouraged to refer
to~\refcartpaper{}, Biot's papers of 1956~\cite{biot:56-1, biot:56-2}
and 1962~\cite{biot:62}, or Carcione's book~\cite{carcione:wave-book}.

In addition to purely poroelastic systems, we also model systems
composed of a combination of poroelastic and fluid media.  We model
the fluid parts of these systems using linear acoustics; this subject
is well-understood, so we discuss it only briefly here.  Of more
interest are the conditions that must be satisifed at the interfaces
between different media, both between a fluid and a poroelastic medium
and between two distinct poroelastic media.  This section covers these
interface conditions in their root form at the PDE level, though we
defer discussion of their implementation in the numerical code to
Section \ref{sec:mapped-fvm}.

\subsection{Poroelasticity system}

First, define the state vector $\bQ$ at a point in a poroelastic
medium as
\begin{equation} \label{eq:state-vector-new}
  \bQ = \begin{pmatrix}
    p & \tau_{xx} & \tau_{zz} & \tau_{xz} & v_x & v_z & q_x & q_z
  \end{pmatrix}^T.
\end{equation}
Here $p$ is the pressure of the fluid in the pores, $\btau$ is the
total stress tensor, $\mathbf{v}$ is the velocity of the solid
component of the medium, and $\mathbf{q}$ is the flow rate of the
fluid relative to the solid.  The $x$ and $z$ axes are the global axes
in which the problem is defined.  Note that this state vector contains
the same quantities as the one defined in~\refcartpaper{}, but
arranged in a different order.

For this work, we deal with orthotropic media: media that possess a
set of three perpendicular planes of symmetry, so that in the axes
defined by these planes, extensional and shear deformation are
decoupled.  We call these the {\em principal axes} of the medium, and
denote them by the numbers $1$, $2$, and $3$ to distinguish them from
the $x$ and $z$ axes.  We also make the assumption that the our
poroelastic media are transversely isotropic --- specifically,
isotropic in the $1$-$2$ plane --- and that axes $1$ and $3$ lie in
the $x$-$z$ plane.  This type of anisotropy is common in
engineering composites~\cite{gibson:composites} and in some biological
materials~\cite{cowin-mehrabadi:elastic-symmetry}, as well as being
present in certain types of stone.  However, this assumption is not a
fundamental requirement of our work --- it simplifies the system
matrices, but it would be straightforward to extend our formulation to
general anisotropic materials.

Following~\refcartpaper{}, in a homogeneous medium with no sources the
equations of Biot poroelasticity theory at low frequency can be
expressed as an $8 \times 8$ first-order linear system of PDEs,
\begin{equation} \label{eq:first-order-system}
  \partial_t \bQ + \bA \partial_x \bQ + \bB \partial_z \bQ = \bD \bQ.
\end{equation}
We also frequently work with the homogeneous form of this system,
\begin{equation} \label{eq:first-order-homogeneous}
  \partial_t \bQ + \bA \partial_x \bQ + \bB \partial_z \bQ = 0.
\end{equation}
In particular, the homogeneous form is the form of the system we
use when solving the Riemann problem for our high-resolution finite
volume scheme.

With the state vector in \eqref{eq:state-vector-new}, the coefficient
matrices can be written in block form,
\begin{equation} \label{eq:system-matrices-block}
  \bA = \begin{pmatrix} 0_{4\times 4} & \bA_{sv} \\ \bA_{vs} &
    0_{4\times 4} \end{pmatrix}, \quad
  \bB = \begin{pmatrix} 0_{4\times 4} & \bB_{sv} \\ \bB_{vs} &
    0_{4\times 4} \end{pmatrix}, \quad
  \bD = \begin{pmatrix} 0_{4\times 4} & 0_{4\times 4} \\ 0_{4\times 4}
    & \bD_{v} \end{pmatrix}.
\end{equation}
If the computational $x$-$z$ axes coincide with the principal $1$-$3$
axes of the material, the blocks in \eqref{eq:system-matrices-block}
are
\begin{align}
  \label{eq:A-blocks} \bA_{sv} &= \begin{pmatrix}
    \alpha_1 M & 0 & M & 0\\     
    -c_{11}^u & 0 & -\alpha_1 M & 0\\
    -c_{13}^u & 0 & -\alpha_3 M & 0\\
    0 & -c_{55}^u & 0 & 0
  \end{pmatrix} &
  \bA_{vs} &= \begin{pmatrix}
    -\frac{\rho_f}{\Delta_1} & -\frac{m_1}{\Delta_1} & 0 & 0 \\
    0 & 0 & 0 & -\frac{m_3}{\Delta_3} \\
    \frac{\rho}{\Delta_1} & \frac{\rho_f}{\Delta_1} & 0 & 0\\
    0 & 0 & 0 & \frac{\rho_f}{\Delta_3}
  \end{pmatrix}\\
  \label{eq:B-blocks} \bB_{sv} &= \begin{pmatrix}
    0 & \alpha_3 M & 0 & M \\
    0 & -c_{13}^u & 0 & -\alpha_1 M\\
    0 & -c_{33}^u & 0 & -\alpha_3 M\\
    -c_{55}^u & 0 & 0 & 0
  \end{pmatrix} &
  \bB_{vs} &= \begin{pmatrix}
    0 & 0 & 0 & -\frac{m_1}{\Delta_1}\\
    -\frac{\rho_f}{\Delta_3} & 0 & -\frac{m_3}{\Delta_3} & 0\\
    0 & 0 & 0 & \frac{\rho_f}{\Delta_1}\\
    \frac{\rho}{\Delta_3} & 0 & \frac{\rho_f}{\Delta_3} & 0
  \end{pmatrix}\\
  \label{eq:D-blocks} \bD_v &= \begin{pmatrix}
    0 & 0 & \frac{\rho_f\eta}{\Delta_1\kappa_1} & 0 \\
    0 & 0 & 0 & \frac{\rho_f\eta}{\Delta_3\kappa_3}\\
    0 & 0 & -\frac{\rho\eta}{\Delta_1\kappa_1} & 0\\
    0 & 0 & 0 & -\frac{\rho\eta}{\Delta_3\kappa_3}
  \end{pmatrix}.
\end{align}
Here the subscripts $s$ and $v$ denote the stress variables ($p$ and
$\btau$) and velocity variables ($\B{v}$ and $\B{q}$), respectively.
The entries in these matrices are determined from the physical
properties in Table \ref{tab:matprops}; the parameters $c_{ij}^u$ are
the undrained elastic stiffness constants, which are determined from the
drained stiffness constants $c_{ij}$, the effective stress
coefficients $\alpha_i$, and the bulk compressibility parameter $M$.
Ordering the variables to highlight this block structure emphasizes
the underlying physics --- gradients of stress produce changes
in velocity, and gradients of velocity produce changes in stress ---
but it also proves useful mathematically.  If the $x$-$z$ axes do
not coincide with the material principal axes, we can obtain
appropriate $\bA$, $\bB$, and $\bD$ matrices by
transforming the state variables $\bQ$ into the new axes and applying the
chain rule of partial differentiation.  We refer the reader to our
previous paper~\refcartpaper{}, or to Carcione's
book~\cite{carcione:wave-book}, for full explanation of the parameters
in these matrices.

\subsection{Energy density for poroelasticity}

In~\refcartpaper{} we derived the energy density associated with a
state vector $\bQ$ as
\begin{equation}
  \mathcal{E} = \frac{1}{2} \bQ^T \bE \bQ.
\end{equation}
The Hessian of the energy density $\bE$ remains symmetric after the
elements of $\bQ$ are permuted to the order used here; with the new
ordering, $\bE$ takes the block diagonal form
\begin{equation} \label{eq:energy-matrix-block}
  \bE = \begin{pmatrix} \bE_s & 0_{4\times 4} \\ 0_{4\times 4} &
    \bE_v \end{pmatrix},
\end{equation}
where in the principal material axes the diagonal blocks of $\bE$ are
\begin{equation} \label{eq:E-blocks}
  \begin{aligned}
  \bE_s &= \begin{pmatrix}
    \frac{1}{M} + \frac{\alpha_1^2 c_{33} + \alpha_3^2 c_{11} -
    2\alpha_1\alpha_3 c_{13}}{c_{11}c_{33} -
      \left(c_{13}\right)^2} & \frac{\alpha_1 c_{33} -
      \alpha_3 c_{13}}{c_{11}c_{33} -
      \left(c_{13}\right)^2} &
    \frac{\alpha_3 c_{11} - \alpha_1 c_{13}}{c_{11}c_{33} - \left(c_{13}\right)^2} &
    0 \\
    \frac{\alpha_1 c_{33} - \alpha_3 c_{13}}{c_{11}c_{33} -
      \left(c_{13}\right)^2} & \frac{c_{33}}{c_{11}c_{33} - \left(c_{13}\right)^2} &
    -\frac{c_{13}}{c_{11}c_{33} - \left(c_{13}\right)^2} & 0 \\
    \frac{\alpha_3 c_{11} - \alpha_1 c_{13}}{c_{11}c_{33} -
      \left(c_{13}\right)^2} & -\frac{c_{13}}{c_{11}c_{33} - \left(c_{13}\right)^2} &
    \frac{c_{11}}{c_{11}c_{33} - \left(c_{13}\right)^2} & 0 \\
    0 & 0 & 0 & \frac{1}{c_{55}} \\
  \end{pmatrix}\\
  \bE_v &= \begin{pmatrix}
    \rho & 0 & \rho_f & 0 \\
    0 & \rho & 0 & \rho_f \\
    \rho_f & 0 & m_1 & 0 \\
    0 & \rho_f & 0 & m_3
  \end{pmatrix}.
  \end{aligned}
\end{equation}
Note that the individual blocks $\bE_s$ and $\bE_v$ are themselves
symmetric positive-definite matrices.  The parameters $c_{ij}$ are the
{\em drained} elastic stiffness constants.

This energy density allows us to define in a natural way a norm for
the state vector $\bQ$ that bypasses problems with relative scaling of
its components when they are expressed in typical unit systems.  We
define the energy inner product of two state vectors $\bQ_1$ and
$\bQ_2$ as
\begin{equation}
  \langle \bQ_1, \bQ_2 \rangle_E := \bQ_2^H \bE \bQ_1.
\end{equation}
(We use the Hermitian conjugate-transpose here to ensure that this
remains an inner product if it is extended to complex vectors.)  From
this inner product, we define the induced energy norm,
\begin{equation} \label{eq:energy-norm-def}
  \| \bQ \|_E := \sqrt{\langle \bQ, \bQ\rangle_E}.
\end{equation}

As discussed in~\refcartpaper{}, the matrix $\bE$ allows us to show
several useful properties of the system, including the
$\bE$-orthogonality of the eigenvectors of $\bA$, $\bB$, or any linear
combination of them, and the fact that the energy density
$\mathcal{E}$ is a strictly convex entropy function of the system.
The $\bE$-orthogonality of the eigenvectors also allows an easy proof
that the system is hyperbolic.

\subsection{Linear acoustics}

The partial differential equations governing linear acoustics are
well-known, and we will not re-derive them here.  We will, however,
state how they are incorporated into the same framework as
poroelasticity.

We use the same form of first-order linear system as
\eqref{eq:first-order-system} to model acoustic wave propagation in a
fluid, with the same state vector; however, in a fluid we define the
variables $\btau$ and $\vect{v}$ to be identically zero.  We use $p$
for the fluid pressure and $\vect{q}$ for its velocity.  (In fact the
total stress tensor in the fluid is $-p\B{I}$, but it is more
convenient to use the single pressure variable and ignore $\btau$ in
the fluid.)
The appropriate coefficient matrices have the same block form as for
poroelasticity, with blocks given by
\begin{align}
  \label{eq:A-blocks-acou} \bA_{sv} &= \begin{pmatrix}
    0 & 0 & K_f & 0 \\
    0 & 0 & 0 & 0 \\
    0 & 0 & 0 & 0 \\
    0 & 0 & 0 & 0
  \end{pmatrix} &
  \bA_{vs} &= \begin{pmatrix}
    0 & 0 & 0 & 0 \\
    0 & 0 & 0 & 0 \\
    \frac{1}{\rho_f} & 0 & 0 & 0 \\
    0 & 0 & 0 & 0
  \end{pmatrix} \\
  \label{eq:B-blocks-acou} \bB_{sv} &= \begin{pmatrix}
    0 & 0 & 0 & K_f \\
    0 & 0 & 0 & 0 \\
    0 & 0 & 0 & 0 \\
    0 & 0 & 0 & 0
  \end{pmatrix} &
  \bB_{vs} &= \begin{pmatrix}
    0 & 0 & 0 & 0 \\
    0 & 0 & 0 & 0 \\
    0 & 0 & 0 & 0 \\
    \frac{1}{\rho_f} & 0 & 0 & 0
  \end{pmatrix}.
\end{align}
The dissipation matrix $\bD$ is identically zero for a fluid.

Similarly to a poroelastic medium, we can also write a matrix $\bE$
such that the energy density in the fluid is $\mathcal{E} =
\frac{1}{2} \bQ^T \bE \bQ$.  This $\bE$ matrix has the same block
structure as for poroelasticity; its blocks are
\begin{equation} \label{eq:E-blocks-acou}
  \begin{aligned}
    \bE_s &= \begin{pmatrix}
      \frac{1}{K_f} & 0 & 0 & 0 \\
      0 & 0 & 0 & 0 \\
      0 & 0 & 0 & 0 \\
      0 & 0 & 0 & 0
    \end{pmatrix}, &
    \bE_v &= \begin{pmatrix}
      0 & 0 & 0 & 0 \\
      0 & 0 & 0 & 0 \\
      0 & 0 & \rho_f & 0 \\
      0 & 0 & 0 & \rho_f
    \end{pmatrix}.
  \end{aligned}
\end{equation}
Note that the $\bE$ defined this way is only positive-semidefinite,
not positive-definite as for poroelasticity.  However, the null space
of $\bE$ consists only of the variables that are defined to be
identically zero in the fluid, so it is essentially positive-definite,
and it is still sensible to use it to define an energy inner product
and norm in the fluid.

\subsection{Interface conditions}

Nontrivial conditions relating the state variables on either side
of an interface between distinct poroelastic materials have been
proposed by a number of authors.  These include Deresiewicz and
Skalak~\cite{deresiewicz-skalak:poro-uniqueness-1963}, who proposed an
imperfect hydraulic contact condition relating the pressure difference
across the interface to the normal fluid flow rate, and
Sharma~\cite{sharma:dissimilar-boundary-2008}, who formulated a loose
contact condition modeling lubricated slippage between the two sides.
For this work, we use Deresiewicz and Skalak's imperfect hydraulic
contact condition; Deresiewicz and Skalak showed that this condition
is sufficient to give a unique solution to Biot's equations in a
discontinuous medium, and Gurevich and
Schoenberg~\cite{gurevich-schoenberg:interface-1999} examined how such
a condition could arise asymptotically from a smoothly varying medium
as the region over which the material properties vary is shrunk to
zero thickness.

The imperfect hydraulic contact condition can be written in the form
\begin{equation} \label{eq:poro-poro-physical-base}
  \begin{aligned}
    \btau_l \cdot \vect{n} &= \btau_r \cdot \vect{n}\\
    \vect{v}_l &= \vect{v}_r\\
    \vect{q}_l \cdot \vect{n} &= \vect{q}_r \cdot \vect{n}\\
    p_l - p_r &= \frac{1}{\mathcal{K}} \widehat{\vect{q} \cdot \vect{n}}.
  \end{aligned}
\end{equation}
Here the subscripts $l$ and $r$ represent the left and right sides of
the interface, chosen arbitrarily, and $\vect{n}$ is the unit normal
to the interface, pointing from left to right.  These equations have
direct physical significance: the first is a statement of the
continuity of traction across the interface, the second states that
the materials stay connected to each other, the third requires that
all fluid entering the interface should exit the other side, and the
fourth relates the fluid flow rate across the interface to the
pressure difference forcing it across.  The parameter $\mathcal{K}$ is
a measure of the ability of the interface to conduct fluid, and ranges
from zero, representing a completely impermeable interface, to
$+\infty$, representing no impedance to fluid flow.  The quantity
$\widehat{\vect{q} \cdot \vect{n}}$ is the volume flow rate of fluid
across the interface, which is equal to both $\vect{q}_l \cdot
\vect{n}$ and $\vect{q}_r \cdot \vect{n}$ according to the third
equation of \eqref{eq:poro-poro-physical-base}.  We will revisit this
ambiguity when we discuss the implementation of this interface
condition.  Our previous work~\refcartpaper{}, which simply used the
wave structure of the system without any explicit interface condition,
was equivalent to this condition with $\mathcal{K} = +\infty$;
Gurevich and Schoenberg~\cite{gurevich-schoenberg:interface-1999}
showed that this is the most natural interface condition in the
context of Biot's equations as a PDE system, and we continue to use
the same approach where no discontinuities in the medium are
present.

Since both $\mathcal{K} = 0$ and $\mathcal{K} = +\infty$ are common
and important cases, and since infinite values are inconvenient in
numerical computations, we reparameterize the fourth equation of
\eqref{eq:poro-poro-physical-base}.  Noting that $1/\mathcal{K}$ has
the same units as acoustic impedance, we define
\begin{equation} \label{eq:discharge-efficiency-def}
\frac{1}{\mathcal{K}} =: Z_f \frac{1 - \eta_d}{\eta_d},
\end{equation}
where $Z_f$ is the impedance of the pore fluid in the left medium and
$\eta_d \in [0,1]$ is a new nondimensional parameter we term the {\em
  interface discharge efficiency}.  Setting $\mathcal{K} = +\infty$
now corresponds to setting $\eta_d = 1$, and $\mathcal{K} = 0$
corresponds to $\eta_d = 0$.  Substituting this into the fourth
equation of \eqref{eq:poro-poro-physical-base} and multiplying through
by $\eta_d$, we obtain
\begin{equation} \label{eq:poro-poro-physical}
  \begin{aligned}
    \btau_l \cdot \vect{n} &= \btau_r \cdot \vect{n}\\
    \vect{v}_l &= \vect{v}_r\\
    \vect{q}_l \cdot \vect{n} &= \vect{q}_r \cdot \vect{n}\\
    \eta_d (p_l - p_r) &= Z_f (1 - \eta_d) \widehat{\vect{q} \cdot \vect{n}},
  \end{aligned}
\end{equation}
which presents no special difficulty for any value of $\eta_d$,
including 0 and 1.  This
is the form of the interface condition that we use in our numerical
computations.

Between a poroelastic medium and a fluid, we use a similar condition,
which we write as
\begin{equation} \label{eq:poro-fluid-physical-base}
  \begin{aligned}
    \vect{q}_f \cdot \vect{n} &= (\vect{v}_p + \vect{q}_p)
    \cdot \vect{n}\\
    -p_f \vect{n} &= \btau_p \cdot \vect{n}\\
    p_p - p_f &= \frac{1}{\mathcal{K}} \vect{q}_p \cdot \vect{n}.
  \end{aligned}
\end{equation}
Here, the subscript $f$ indicates quantities in the fluid, while $p$
indicates quantities in the poroelastic medium.  We use $\vect{q}_f$
for the fluid velocity, and the unit interface normal $\vect{n}$ points
from the poroelastic medium into the fluid.  The meanings of the
equations are similar to before: the first states continuity of fluid
flow, the second, continuity of traction, and the third relates fluid
flow rate to pressure difference.  The parameter $\mathcal{K}$ has the same meaning
as before, though its value may be different.  This condition models
possibly imperfect hydraulic contact; it seems to have been first used
at a poroelastic-fluid interface by
Rosenbaum~\cite{rosenbaum:synthetic-microseismo}, and has since been
employed by other authors~\cite{bcz:acoustics-porous,
  feng-johnson:fluid-poro-hifreq-surface}.  Chiavassa and
Lombard~\cite{chiavassa-lombard:acou-poro} in particular use this
method in their numerical work, and also demonstrate that the
resulting coupled systems of PDEs on the fluid and poroelastic domains
are well-posed.  For easier implementation, we again
replace $1/\mathcal{K}$ with $Z_f (1 - \eta_d)/\eta_d$, where $Z_f$ is
the acoustic impedance of the fluid medium, obtaining the alternative
form
\begin{equation} \label{eq:poro-fluid-physical}
  \begin{aligned}
    \vect{q}_f \cdot \vect{n} &= (\vect{v}_p + \vect{q}_p)
    \cdot \vect{n}\\
    -p_f \vect{n} &= \btau_p \cdot \vect{n}\\
    \eta_d (p_p - p_f) &=  Z_f (1 - \eta_d) \vect{q}_p \cdot \vect{n}.
  \end{aligned}
\end{equation}

\section{Finite volume methods for poroelastic-fluid systems on logically rectangular mapped grids}
\label{sec:mapped-fvm}

In order to be able to model geometries that do not lend themselves to
the straight lines and 90 degree angles of a rectilinear grid, we use
logically rectangular mapped grids.  Compared to
unstructured grids, mapped grids have the
advantage of simpler data structures and less computational overhead,
although creating the desired mapping function is not always trivial.
They also combine well with finite volume methods, since the finite
volume solution tends to maintain good quality even in the face of a
mapping that severely distorts the grid.  However, when modeling
anisotropic poroelastic media, having to deal with cell interfaces
that may be oriented in any arbitrary direction means that the simple
compute-and-cache method used to solve the Riemann problem in our
previous work no longer suffices.  Section \ref{sec:riemann-solve}
discusses a Riemann solution process designed to function efficiently
in this context, and to also incorporate interface conditions such as
\eqref{eq:poro-poro-physical} and \eqref{eq:poro-fluid-physical}.

An additional difficulty we encounter is that at an interface between
a poroelastic medium and a fluid, the classical formulation of a
high-resolution finite volume method can be qualitatively incorrect
--- the classical transverse and second-order correction fluxes result
in poroelastic variables such as skeleton stress and solid velocity
being carried into the fluid, where they make no sense.
We discuss this further in section \ref{sec:fvm-mods}, and are able to produce
qualitatively correct solutions, though the best way to fully
reformulate the method, including generalizations of the classical
limiters and second-order correction fluxes, is still an open
question.

In the remaining parts of this section, we briefly discuss the
implementation of the source term in the poroelasticity system and the
\clawpack{} software framework in which we implement our numerical
solution.

\subsection{Mapped grids}
\label{sec:mapped-grid-basics}

\begin{figure}[t]
  \begin{center}
    \includegraphics[width=0.8\textwidth]{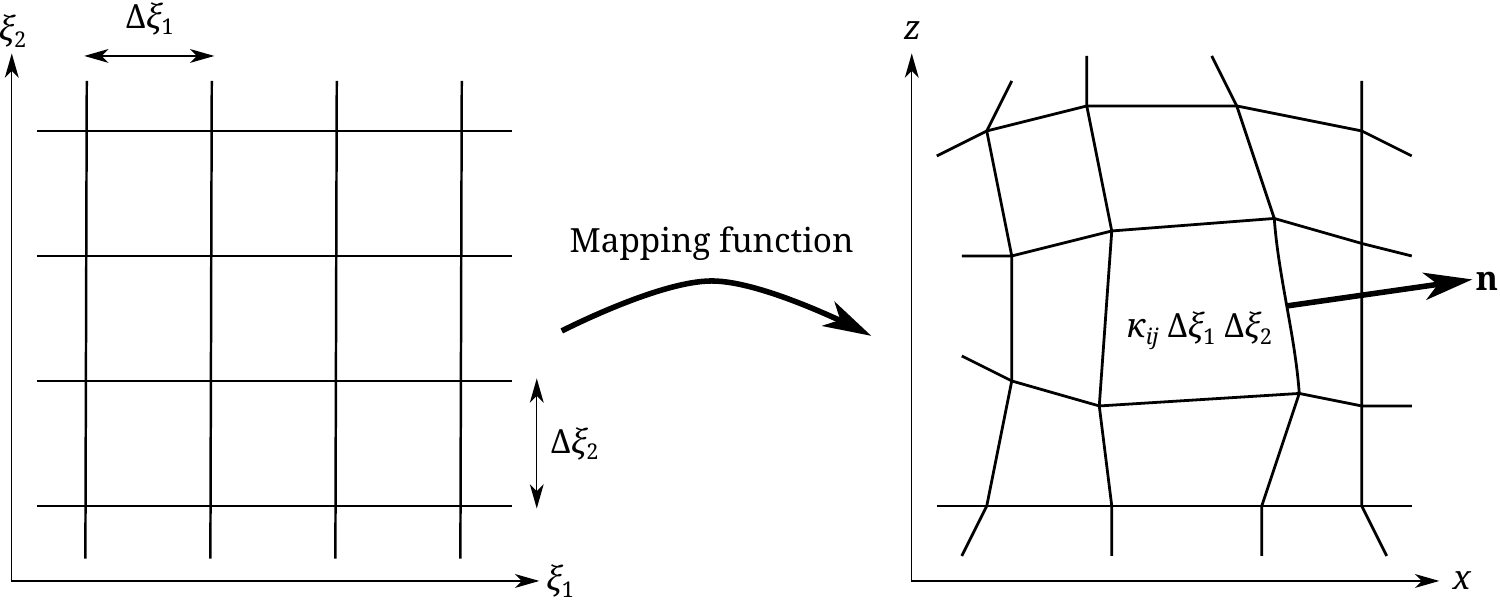}
    \caption{Example of a mapped grid, showing computational space at
      the left and physical space at the
      right. \label{fig:map-figure}}
  \end{center}
\end{figure}

For a mapped grid approach, we start with a uniform rectangular grid
in the computational coordinates, denoted $\xi_1$ and $\xi_2$.  We
then apply a mapping function $\B{X}(\xi_1,\xi_2)$ to obtain the grid
in physical coordinates $(x,z)$.  This mapping function is typically
chosen so that the grid boundaries or interior grid lines follow some
feature of interest in the problem, although once the grid mapping is
formulated, the actual cells are taken to be quadrilaterals with
straight sides, whose vertices are obtained using the mapping
function.  Each grid cell $ij$ has an associated capacity
$\kappa_{ij}$, which is the ratio of the area of the cell in physical
coordinates to its area $\Delta \xi_1\, \Delta \xi_2$ in computational
coordinates, and each cell interface has an associated unit normal
vector $\B{n}$ pointing in the positive grid direction.  Figure
\ref{fig:map-figure} shows an example of a mapped grid.  Finite volume
methods on mapped grids are discussed in greater detail in Chapter 23
of~\cite{rjl:fvm-book}.

\subsection{Riemann problems on mapped grids}
\label{sec:riemann-solve}

One of the most critical parts of any high-resolution finite volume
code is the solution of the Riemann problem.  Solving the Riemann
problem for poroelasticity at first appears challenging due to the
size and complexity of the hyperbolic part of the system
\eqref{eq:first-order-system}.  However, by taking advantage of the
structure of the system, a straightforward and fairly efficient
solution is possible.

\subsubsection{Eigenvalues and eigenvectors}

The first order of business in solving the Riemann problem is finding
the eigenvectors of $\breve{\bA} = n_x \bA + n_z \bB$ corresponding to
the propagating waves, and the corresponding eigenvalues giving the
wave speeds.  For linear acoustics, this eigensystem is simple: the
eigenvalues are $\pm \sqrt{K_f/\rho_f}$, and the eigenvectors may be
easily verified as
\begin{equation}
  \begin{aligned}
    \B{r}_{\text{acoustic, left}} &= \begin{pmatrix}
      -Z_f & 0 & 0 & 0 & 0 & 0 & n_x & n_z
    \end{pmatrix}^T && \text{(left-going wave)}\\
    \B{r}_{\text{acoustic, right}} &= \begin{pmatrix}
      Z_f & 0 & 0 & 0 & 0 & 0 & n_x & n_z
    \end{pmatrix}^T && \text{(right-going wave)}.
  \end{aligned}
\end{equation}
For poroelasticity the eigensystem is substantially more complicated,
but it is possible to exploit the block structure, symmetry, and
simple null space of the poroelastic system to reduce the eigenproblem
from an $8\times 8$ non-symmetric to a $3\times 3$ real symmetric one;
see Appendix \ref{sec:riemann-solver} for details.  As a note on implementation, our code caches
the eigensystems computed for each Riemann problem, and checks whether
the next Riemann problem has the same input data to the eigensolves;
if so, the previous eigensystem is re-used.  This speeds up execution
in the common special case of large sections of the domain that have
straight grid lines and uniform composition.

\subsubsection{Solution between identical materials}

The solution of the Riemann problem for linear acoustics with
identical media on either side of the interface is simple and
well-known, and will not be repeated here.  The interested reader may
refer to~\cite{rjl:fvm-book}, as one of many sources.

Solution of the Riemann problem for poroelasticity with identical
materials on either side of the interface is more complex, but it can
be simplified greatly by taking advantage of the structure
of the system.  In
particular, as shown in~\refcartpaper{}, the eigenvectors of the
poroelasticity system are orthogonal with respect to the energy matrix
$\bE$, and in fact our eigensolution process detailed in Appendix
\ref{sec:riemann-solver} produces
eigenvectors that are orthonormal with respect to $\bE$.  Since we do
not apply any special interface condition at cell interfaces within a
homogeneous medium, to solve the Riemann
problem we simply seek a vector of wave strengths $\bbeta$ such that
\begin{equation}
  \B{R} \bbeta = \bQ_r - \bQ_l,
\end{equation}
where $\B{R}$ is the matrix of eigenvectors of $\breve{\bA}$.  By the
orthonormality of the eigenvectors with respect to $\bE$, we can
multiply from the left by $\B{R}^T \bE$ and obtain
\begin{equation}
  \bbeta = \B{R}^T \bE (\bQ_r - \bQ_l).
\end{equation}

\subsubsection{Solution between different materials}
\label{sec:riemann-interface}

At an interface between different materials, the Riemann solution
process becomes more complex.  For a pair of fluids with different
properties, the solution process of the acoustic Riemann problem is
still well-known, and is covered for instance in~\cite{rjl:fvm-book}.
For a pair of distinct poroelastic media, or a poroelastic medium and
a fluid, we must satisfy an interface condition, either
\eqref{eq:poro-poro-physical} or \eqref{eq:poro-fluid-physical}.

To solve the Riemann problem in this context, we first note that both
interface conditions can be cast in the form
\begin{equation} \label{eq:interface-linear-general}
  \B{C}_l \lim_{x' \to 0^-}\bQ = \B{C}_r \lim_{x' \to 0^+} \bQ
\end{equation}
for some matrices $\B{C}_l$ and $\B{C}_r$.  Here $x'$ is the signed
normal distance from the interface, and is positive on the right side.
The solution within the left and right media is still a set of
discontinuities in $\bQ$ proportional to the eigenvectors of the
$\breve{\bA}$ matrices in the respective media, propagating at speeds
equal to the corresponding eigenvalues, so we can relate $\lim_{x' \to
  0^-} \bQ$ and $\lim_{x' \to 0^+} \bQ$ to the left and right states
$\bQ_l$ and $\bQ_r$ by
\begin{equation} \label{eq:interface-states}
  \lim_{x' \to 0^-} \bQ = \bQ_l + \sum_{\text{left}}
  \beta_{li} \B{r}_{li}, \quad
  \lim_{x' \to 0^+} \bQ = \bQ_r - \sum_{\text{right}}
  \beta_{ri} \B{r}_{ri}.
\end{equation}
Here $\beta_{li}$ and $\beta_{ri}$ are the strengths of the left-going
and right-going waves, $\B{r}_{li}$ and $\B{r}_{ri}$ are the
corresponding eigenvectors, and the sums are over only the left-going
and right-going waves, respectively.  Substituting this into
\eqref{eq:interface-linear-general} and rearranging, we obtain a
linear system for the wave strengths,
\begin{equation} \label{eq:interface-linear-system}
  \begin{pmatrix} \B{C}_l \B{R}_l & \B{C}_r \B{R}_r \end{pmatrix}
  \begin{pmatrix} \bbeta_l \\ \bbeta{r} \end{pmatrix} =
  \B{C}_r \bQ_r - \B{C}_l \bQ_l,
\end{equation}
where $\B{R}_l$ and $\B{R}_r$ are the matrices of left-going and
right-going wave eigenvectors, respectively.  Note that while it is
possible to formulate this linear system, we have no {\em a priori}
guarantee that it has a solution; at a minimum, there must be exactly
as many equations in the interface condition
\eqref{eq:interface-linear-general} as there are propagating waves.
Solution of this system has, however, not been a problem for any of
the cases considered here.

Now all that remains is to explicitly write the matrices $\B{C}_l$ and
$\B{C}_r$ corresponding to \eqref{eq:poro-poro-physical} and
\eqref{eq:poro-fluid-physical}.  We address
\eqref{eq:poro-fluid-physical} first; if we take the left medium to be
poroelastic, a straightforward component-by-component accounting gives
\begin{equation} \label{eq:poro-fluid-interface-matrices}
  \begin{aligned}
    \B{C}_{l, \text{poro-fluid}} &= \begin{pmatrix}
      0 & 0 & 0 & 0 & n_x & n_z & n_x & n_z \\
      0 & n_x & 0 & n_z & 0 & 0 & 0 & 0 \\
      0 & 0 & n_z & n_x & 0 & 0 & 0 & 0 \\
      \eta_d & 0 & 0 & 0 & 0 & 0 & -Z_f (1 - \eta_d)
      n_x & -Z_f (1 - \eta_d) n_z
    \end{pmatrix}\\
    \B{C}_{r, \text{poro-fluid}} &= \begin{pmatrix}
      0 & 0 & 0 & 0 & 0 & 0 & n_x & n_z \\
      -n_x & 0 & 0 & 0 & 0 & 0 & 0 & 0 \\
      -n_z & 0 & 0 & 0 & 0 & 0 & 0 & 0 \\
      \eta_d & 0 & 0 & 0 & 0 & 0 & 0 & 0
    \end{pmatrix}.
  \end{aligned}
\end{equation}
If the poroelastic material is on the right side of the interface, we
can simply exchange the subscripts $l$ and $r$, and negate $\vect{n}$.

Writing the matrices $\B{C}_l$ and $\B{C}_r$ for
\eqref{eq:poro-poro-physical} is less straightforward, due to the
ambiguity in $\widehat{\vect{q} \cdot \vect{n}}$.  This quantity is
equal to both $\vect{q}_l \cdot \vect{n}$ and $\vect{q}_r \cdot
\vect{n}$; in exact arithmetic it is irrelevant which one we choose,
but it is not obvious how best to define $\widehat{\vect{q} \cdot
  \vect{n}}$ for numerical solution.  Faced with this ambiguity, we
let $\widehat{\vect{q} \cdot \vect{n}} = (1-\zeta) \vect{q}_l \cdot
\vect{n} + \zeta \vect{q}_r \cdot \vect{n}$, where $\zeta$ is a free
parameter used to improve the conditioning of
\eqref{eq:interface-linear-system}.  With this choice of
$\widehat{\vect{q} \cdot \vect{n}}$, $\B{C}_l$ and $\B{C}_r$ become
\begin{equation}
  \begin{aligned}
    \B{C}_{l, \text{poro-poro}} &= \begin{pmatrix}
      0 & n_x & 0 & n_z & 0 & 0 & 0 & 0 \\
      0 & 0 & n_z & n_x & 0 & 0 & 0 & 0 \\
      0 & 0 & 0 & 0 & 1 & 0 & 0 & 0 \\
      0 & 0 & 0 & 0 & 0 & 1 & 0 & 0 \\
      0 & 0 & 0 & 0 & 0 & 0 & n_x & n_z \\
      \eta_d & 0 & 0 & 0 & 0 & 0 & -(1-\zeta) Z_f (1 -
      \eta_d) n_x & -(1-\zeta) Z_f (1 - \eta_d) n_z
    \end{pmatrix}\\
    \B{C}_{r, \text{poro-poro}} &= \begin{pmatrix}
      0 & n_x & 0 & n_z & 0 & 0 & 0 & 0 \\
      0 & 0 & n_z & n_x & 0 & 0 & 0 & 0 \\
      0 & 0 & 0 & 0 & 1 & 0 & 0 & 0 \\
      0 & 0 & 0 & 0 & 0 & 1 & 0 & 0 \\
      0 & 0 & 0 & 0 & 0 & 0 & n_x & n_z \\
      \eta_d & 0 & 0 & 0 & 0 & 0 & \zeta Z_f (1 -
      \eta_d) n_x & \zeta Z_f (1 - \eta_d) n_z
    \end{pmatrix}.
  \end{aligned}
\end{equation}

\begin{figure}[t]
  \begin{center}
    \includegraphics[width=0.6\textwidth]{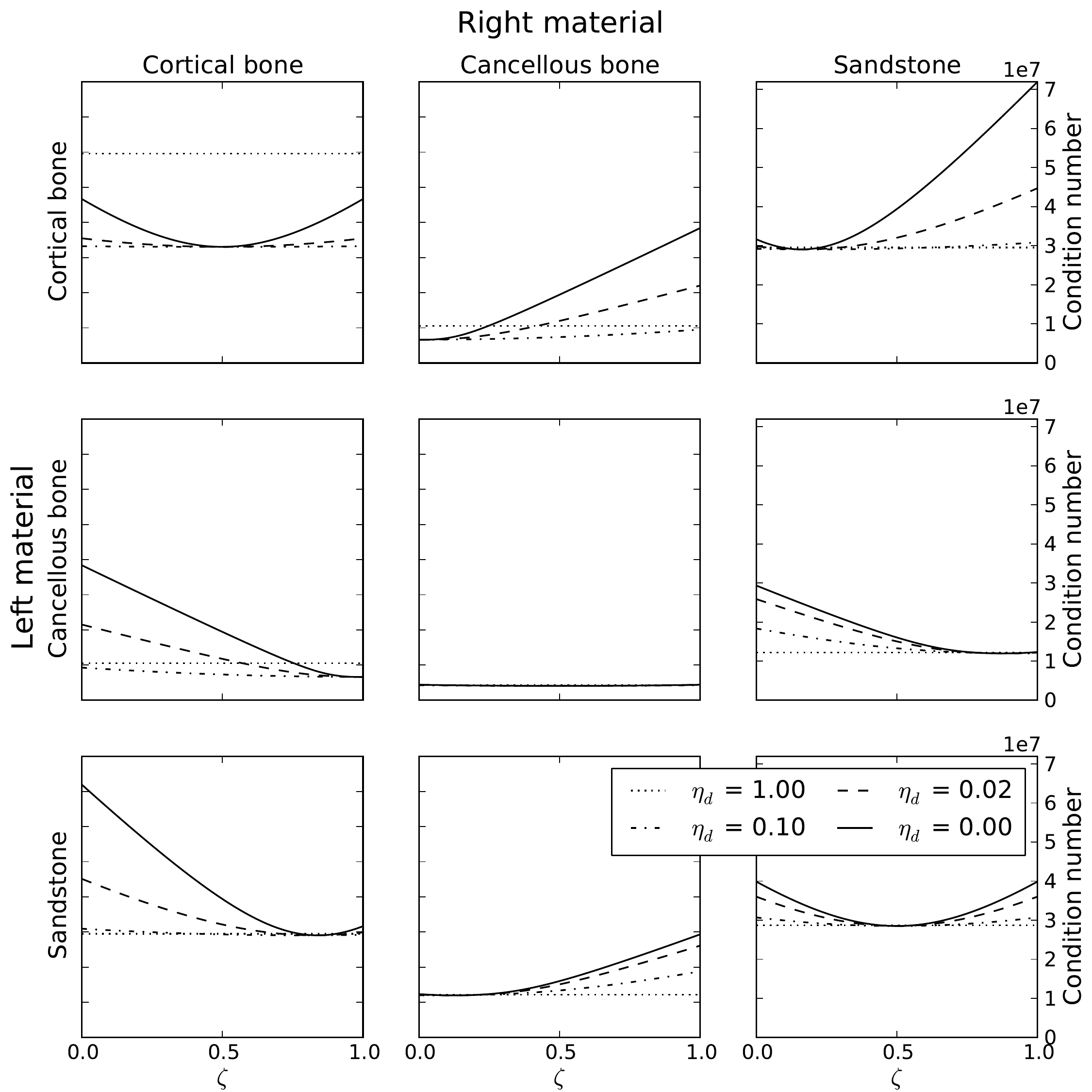}
    \caption{Variation of the condition number of the Riemann solution
      linear system \eqref{eq:interface-linear-system} for interfaces
      between various pairs of poroelastic materials, as a function of
      the free parameter $\zeta$.  Material properties are given in
      Table \ref{tab:matprops}.  The maximum condition number observed
      is around $7 \times 10^7$ --- not large enough to endanger
      accuracy in IEEE double-precision arithmetic. \label{fig:zeta}}
  \end{center}
\end{figure}

Figure \ref{fig:zeta} shows the variation of the condition number of
the coefficient matrix in \eqref{eq:interface-linear-system} as a
function of $\zeta$ for various $\eta_d$ values and various pairs of
materials.  The choice of $\zeta$ is essentially irrelevant except for
impermeable or nearly impermeable interfaces, and even for these
interfaces we see no worse than about a factor of five difference
between the lowest and highest condition number for a given material
combination --- not a significant variation, given that the highest
condition number observed is around $7 \times 10^7$.  For similar
materials, the minimum condition number occurs near $\zeta =
\frac{1}{2}$, and even for very dissimilar materials this value of
$\zeta$ yields a condition number substantially below the worst
possible.  Based on these results, for simplicity we set $\zeta =
\frac{1}{2}$ for all further computations, corresponding to equal
weighting of both sides in computing $\widehat{\vect{q} \cdot
  \vect{n}}$.

\subsection{Modifications to the high-resolution finite volume scheme}
\label{sec:fvm-mods}

The presence of interfaces between media governed by different PDEs
presents a special problem for implementing certain aspects of
high-resolution finite volume methods.  To explore this in more
detail, we first write the full update formula for cell $ij$ from
timestep $n$ to $n+1$ on a mapped grid, following (21.12) and (23.11)
from~\cite{rjl:fvm-book}:
\begin{multline} \label{eq:fvm-update-classic}
  \bQ_{ij}^{n+1} = \bQ_{ij} - \frac{\Delta t}{\kappa_{ij} \Delta \xi_1}
  \left( \mathcal{A}^+ \Delta \bQ_{i-1/2,j} + \mathcal{A}^- \Delta
  \bQ_{i+1/2,j} \right) - \frac{\Delta t}{\kappa_{ij} \Delta \xi_2}
  \left( \mathcal{B}^+ \Delta \bQ_{i,j-1/2} + \mathcal{B}^- \Delta
  \bQ_{i,j+1/2} \right) \\
  - \frac{\Delta t}{\kappa_{ij} \Delta \xi_1} \left(
  \tilde{\B{F}}_{i+1/2,j} - \tilde{\B{F}}_{i-1/2,j} \right)
  - \frac{\Delta t}{\kappa_{ij} \Delta \xi_2} \left(
  \tilde{\B{G}}_{i,j+1/2} - \tilde{\B{G}}_{i,j-1/2} \right).
\end{multline}
All quantities on the right-hand side are evaluated at timestep $n$.
The coordinates $\xi_1$ and $\xi_2$ are computational coordinates,
$\kappa_{ij}$ is the ratio of the area of cell $ij$ to $\Delta \xi_1\,
\Delta \xi_2$, $\mathcal{A}^+ \Delta \bQ_{i-1/2,j}$ etc. are the
fluctuations computed from the Riemann solutions at the interfaces
with neighboring cells, and $\tilde{\B{F}}_{i-1/2,j}$ etc. are the
correction fluxes.  These correction fluxes are a combination of
corrections coming from transverse Riemann solutions, which are
important for stability in multi-dimensional problems, and higher-order
correction fluxes, which allow the method to achieve second-order
accuracy.

It is these correction fluxes that pose a problem here.  Because at a
fluid-poroelastic interface the solutions in the two domains are
governed by different PDEs, it is not possible to construct a single
flux transferring quantities from one domain to the other in this form.  For
example, in the formulation of acoustics we use here, elements 2
through 6 of $\bQ$ are identically zero in a fluid, because they
correspond to the total stress tensor and solid velocity.  The state
of stress in a fluid is described completely by the pressure (since we
ignore viscosity in an all-fluid medium here), and there is no solid
component present to have a velocity, so a flux appropriate to the
poroelastic medium would produce nonsensical results if used to update
the solution in the fluid.  Similarly, a flux appropriate to the fluid
would only update the pressure and fluid flow rate in the poroelastic
medium, whereas we should in general expect all of the state variables
to be updated.  We must therefore reformulate
\eqref{eq:fvm-update-classic} to correctly handle a fluid-poroelastic
interface.

Because we cannot define a single flux across an interface between
media of different types, we instead define two correction
fluctuations, one on either side of the interface --- instead of a
correction flux $\tilde{\B{G}}_{i,j+1/2}$ acting on both cells $(i,j)$
and $(i,j+1)$, we define a $\tilde{\B{G}}^+_{i,j+1/2}$ acting on cell
$(i,j+1)$ and a $\tilde{\B{G}}^-_{i,j+1/2}$ acting on cell $(i,j)$.  These
fluctuations have a similar meaning to the first-order fluctuations
$\mathcal{A}^{\pm} \Delta \bQ$ --- they measure the rate of change of
cell averages caused by waves propagating from the interface.  The
correction fluctuations are not necessarily equal on either side of
the interface, but they are compatible in the sense that they arise
from the same Riemann problems and respect the underlying physics of
both media.  In terms of the correction
fluctuations, our new update formula derived from
\eqref{eq:fvm-update-classic} is
\begin{multline} \label{eq:fvm-update-new}
  \bQ_{ij}^{n+1} = \bQ_{ij} - \frac{\Delta t}{\kappa_{ij} \Delta \xi_1}
  \left( \mathcal{A}^+ \Delta \bQ_{i-1/2,j} + \mathcal{A}^- \Delta
  \bQ_{i+1/2,j} \right) - \frac{\Delta t}{\kappa_{ij} \Delta \xi_2}
  \left( \mathcal{B}^+ \Delta \bQ_{i,j-1/2} + \mathcal{B}^- \Delta
  \bQ_{i,j+1/2} \right) \\
  - \frac{\Delta t}{\kappa_{ij} \Delta \xi_1} \left(
  \tilde{\B{F}}^-_{i+1/2,j} - \tilde{\B{F}}^+_{i-1/2,j} \right)
  - \frac{\Delta t}{\kappa_{ij} \Delta \xi_2} \left(
  \tilde{\B{G}}^-_{i,j+1/2} - \tilde{\B{G}}^+_{i,j-1/2} \right).
\end{multline}
Note that the new correction fluctuations keep the same sign
convention as the old correction fluxes.  Also note that for a Riemann
problem with an interface condition of the form
\eqref{eq:interface-linear-general}, the solution is not necessarily a
function purely of the difference in cell states; we keep the
$\mathcal{A}^\pm \Delta \bQ_{i-1/2,j}$ notation for familiarity's
sake, but it should be interpreted as fluctuations arising from a
Riemann problem at cell interface $(i-1/2,j)$, not as an operator
applied to a difference of states.

The following subsections describe how we compute the
correction fluctuations.

\subsubsection{Transverse Riemann solution}
\label{sec:transverse-riemann-2flux}

\begin{figure}[t]
  \begin{center}
    \subfloat[Classical \label{fig:transverse-comparison-classical}]
             {\includegraphics[scale=0.85]{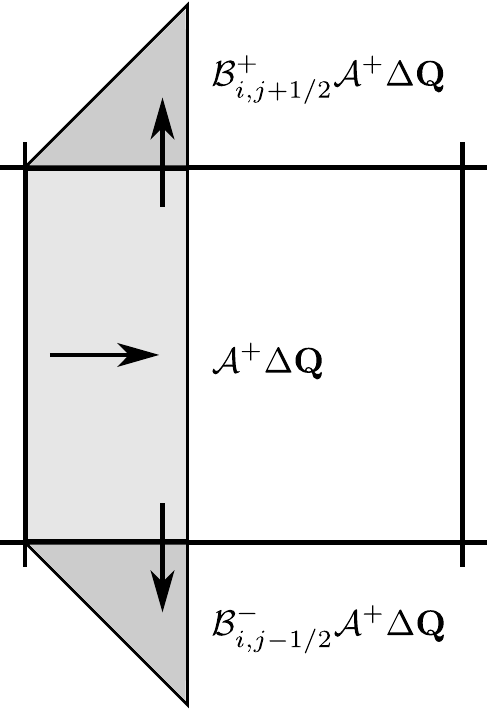}}
    \hspace{1in}
    \subfloat[Fluctuation-based \label{fig:transverse-comparison-fluctuation}]
             {\includegraphics[scale=0.85]{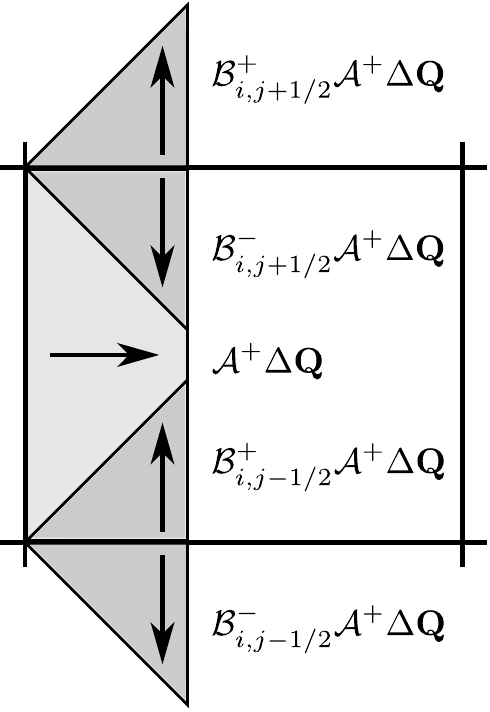}}
    \caption{Comparison of classical and fluctuation-based approaches
      to cell updates from transverse Riemann solves.  Classical
      approach: both cells on either side of the transverse interface
      are updated using the same transverse flux.  Fluctuation-based
      approach: each cell is updated using a separate fluctuation. \label{fig:transverse-comparison}}
  \end{center}
\end{figure}

The transverse Riemann solution process computes the contribution of
the solution on a cell in one timestep to solutions on the cells
diagonally adjacent to it in the next timestep, and is important for
stability in a dimensionally-unsplit high-resolution finite volume
scheme.  Suppose we are considering the contribution to
$\bQ_{ij}^{n+1}$ of the normal Riemann solution at left edge of cell
$(i,j)$.  Ordinarily, the transverse Riemann solution process would
use the right-going fluctuation $\mathcal{A}^+ \Delta \bQ_{i-1/2,j}$
to calculate up-going and down-going transverse fluctuations
$\mathcal{B}^+ \mathcal{A}^+ \Delta \bQ_{i-1/2,j}$ and $\mathcal{B}^-
\mathcal{A}^+ \Delta \bQ_{i-1/2,j}$.  This is illustrated in Figure
\ref{fig:transverse-comparison-classical}.  In the case of a linear
system such as ours, these transverse fluctuations are found by
decomposing $\mathcal{A}^+ \Delta \bQ_{i-1/2,j}$ into eigenvectors of
the flux Jacobians corresponding to the interfaces between cell
$(i,j)$ and cells $(i,j-1)$ and $(i,j+1)$, then multiplying by the
respective eigenvalues and by the ratio of the physical length of the
cell interface to $\Delta \xi_1$; the portion of the
eigendecomposition corresponding to positive eigenvalues is
$\mathcal{B}^+ \mathcal{A}^+ \Delta \bQ_{i-1/2,j}$, while the portion
corresponding to negative eigenvalues is $\mathcal{B}^- \mathcal{A}^+
\Delta \bQ_{i-1/2,j}$.  The correction fluxes
$\tilde{\B{G}}_{i,j+1/2}$ and $\tilde{\B{G}}_{i,j-1/2}$ are then
incremented by $-\frac{\Delta t}{2\Delta \xi_1} \mathcal{B}^+
\mathcal{A}^+ \Delta \bQ_{i-1/2,j}$ and $-\frac{\Delta t}{2\Delta
  \xi_1} \mathcal{B}^- \mathcal{A}^+ \Delta \bQ_{i-1/2,j}$,
respectively.

For our new
approach, illustrated in Figure
\ref{fig:transverse-comparison-fluctuation}, we compute the
contributions of the fluctuation $\mathcal{A}^+ \Delta \bQ_{i-1/2,j}$ to,
for
instance, the two distinct correction fluctuations $\tilde{\B{G}}^{\pm}_{i,j-1/2}$ by applying the Riemann
solver to the states in cells $(i,j)$ and $(i,j-1)$ modified by the
fluctuations from the normal Riemann solve.  (For our particular case,
we exploit linearity by separately computing the transverse
fluctuations from each normal fluctuation, which is the easiest approach to fit into the
existing \clawpack{} framework.  It would also be possible to reduce
computational effort by computing combined transverse fluctuations from the
normal fluctuations in cells in two adjacent rows taken together --- e.g., computing
$\mathcal{A}^+ \Delta \bQ_{i-1/2,j}$ and $\mathcal{A}^+ \Delta
\bQ_{i-1/2,j+1}$, then using them together to compute
$\mathcal{B}^\pm_{i,j+1/2}$ fluctuations by a process very similar to
solving a normal Riemann problem between states $\mathcal{A}^+ \Delta
\bQ_{i-1/2,j}$ and $\mathcal{A}^+ \Delta \bQ_{i-1/2,j+1}$.)
We then use the fluctuations $\mathcal{B}^{\pm} \mathcal{A}^+ \Delta
\bQ_{i-1/2,j}$ from the Riemann solution to increment
$\tilde{\B{G}}^-_{i,j-1/2}$ by $-\frac{\Delta t}{2\Delta \xi_1}
\mathcal{B}^- \mathcal{A}^+ \Delta \bQ_{i-1/2,j}$, and
$\tilde{\B{G}}^+_{i,j-1/2}$ by $+\frac{\Delta t}{2\Delta \xi_1}
\mathcal{B}^+ \mathcal{A}^+ \Delta \bQ_{i-1/2,j}$.  The other
correction fluctuations are updated analogously from transverse
Riemann solutions computed at each cell interface.

To compare these two approaches, we write out the full
contribution to the update of cell $(i,j)$ of the transverse solves based
on the right-going fluctuation from its left edge.  To make the
difference between the methods clear, we attach a subscript to
each transverse fluctuation operator identifying the interface at
which it was computed --- so, for example, the up-going transverse
fluctuation from the upper edge of cell $(i,j)$ is
$\mathcal{B}^+_{i,j+1/2} \mathcal{A}^+ \Delta \bQ_{i-1/2,j}$.  The
classical transverse solve approach increments cell $(i,j)$ by
\begin{equation} \label{eq:classical-transverse-increment}
  \Delta \bQ_{ij,\text{trans},\text{classical}} = -\frac{\Delta t^2}
         {2 \kappa_{ij} \Delta \xi_1\, \Delta \xi_2}
  \left( -\mathcal{B}^+_{i,j+1/2} \mathcal{A}^+ \Delta \bQ_{i-1/2,j}
  +
  \mathcal{B}^-_{i,j-1/2} \mathcal{A}^+ \Delta \bQ_{i-1/2,j}
  \right),
\end{equation}
while the new fluctuation-oriented approach increments the cell by
\begin{equation} \label{eq:fluctuation-transverse-increment}
  \Delta \bQ_{ij,\text{trans},\text{fluctuation}} = -\frac{\Delta t^2}
         {2 \kappa_{ij} \Delta \xi_1\, \Delta \xi_2}
  \left( \mathcal{B}^-_{i,j+1/2} \mathcal{A}^+ \Delta \bQ_{i-1/2,j}
  -
  \mathcal{B}^+_{i,j-1/2} \mathcal{A}^+ \Delta \bQ_{i-1/2,j}
  \right).
\end{equation}
The two approaches produce identical results if the $\mathcal{B}^\pm_{i,j+1/2}$
operators are the same as $\mathcal{B}^\pm_{i,j-1/2}$, but differ
otherwise.  Notice that \eqref{eq:fluctuation-transverse-increment}
updates the cell using fluctuations computed from waves entering it,
while \eqref{eq:classical-transverse-increment} uses fluctuations from
waves leaving it, which may not be appropriate.  Both approaches
produce the same increment to cell $ij$ from the transverse solves
applied to the normal fluctuations from the rows above and below,
$\mathcal{A}^+ \Delta \bQ_{i-1/2,j\pm1}$.

We make one final note here on the broader applicability of this way
of handling transverse propagation, regarding conservation.  While we
do not cast the governing equations of a fluid-poroelastic system
in conservation form, one of the great strengths of
high-resolution finite volume methods is their ability to maintain
conservation of physically relevant quantities such as mass, momentum,
or energy when they are applied to systems of conservation laws.
Going from a single transverse flux to a pair of transverse
fluctuations appears to jeopardize this --- it is not obvious that
$\mathcal{B}^+_{i,j+1/2} \mathcal{A}^+_{i-1/2,j} \Delta \bQ$ will
update cell $(i,j+1)$ by an equal and opposite amount to the update of
$\mathcal{B}^-_{i,j+1/2} \mathcal{A}^+_{i-1/2,j} \Delta \bQ$ to cell
$(i,j)$.  Ensuring this happens, though, is very similar to ensuring
that the two normal fluctuations $\mathcal{B}^{\pm} \Delta \bQ$ update
neighboring cells in a conservative fashion, so it should be no
barrier to practical use.  To demonstrate, suppose that we have some
general, possibly nonlinear problem.  Let $\mathcal{B}^\pm(\bQ_{i,j+1},
\bQ_{ij})$ denote the up-going and down-going fluctuations from the
normal Riemann solver applied to states $\bQ_{i,j+1}$ and $\bQ_{ij}$
in cells $(i,j+1)$ and $(i,j)$.  Assume that for all $i$ and $j$ these
fluctuations satisfy the condition
\begin{equation} \label{eq:conservation-condition}
  \mathcal{B}^+(\bQ_{i,j+1}, \bQ_{ij}) + \mathcal{B}^-(\bQ_{i,j+1},
  \bQ_{ij}) = \B{G}(\bQ_{i,j+1}) - \B{G}(\bQ_{ij}).
\end{equation}
This condition is sufficient to ensure
conservation~\cite{rjl:fvm-book}; the function $\B{G}$ will typically
be the vertical component of the flux vector of the underlying
conservation law, but this is not necessary.  In that
case, a sufficient condition for the new transverse solution scheme to
maintain conservation is to require the output of the transverse
Riemann solve to be related to that of the normal Riemann solve by
\begin{multline} \label{eq:transverse-conservation-construction}
  \mathcal{B}^\pm_{i,j+1/2} \mathcal{A}^+_{i-1/2,j} \Delta \bQ +
  \mathcal{B}^\pm_{i,j+1/2} \mathcal{A}^+_{i-1/2,j+1} \Delta \bQ \\
  = \mathcal{B}^\pm(\bQ_{i,j+1} + \mathcal{A}^+_{i-1/2,j+1}\Delta \bQ,
  \bQ_{ij} + \mathcal{A}^+_{i-1/2,j} \Delta \bQ) - \mathcal{B}^\pm(\bQ_{i,j+1}, \bQ_{ij}).
\end{multline}
Essentially, this condition states that the total of the transverse
fluctuations should be equal to the difference between the up- and
down-going normal fluctuations from the states behind the right-going
waves, and the states in front of them.

If condition \eqref{eq:transverse-conservation-construction} holds,
then summing the transverse fluctuations at interface $(i,j+1/2)$ and
rearranging a bit gives us
\begin{multline}
    \mathcal{B}^+_{i,j+1/2} \mathcal{A}^+_{i-1/2,j} \Delta \bQ +
    \mathcal{B}^+_{i,j+1/2} \mathcal{A}^+_{i-1/2,j+1} \Delta \bQ +
    \mathcal{B}^-_{i,j+1/2} \mathcal{A}^+_{i-1/2,j} \Delta \bQ +
    \mathcal{B}^-_{i,j+1/2} \mathcal{A}^+_{i-1/2,j+1} \Delta \bQ \\
    = (\B{G}(\bQ_{i,j+1} + \mathcal{A}^+_{i-1/2,j+1} \Delta \bQ) - \B{G}(\bQ_{i,j+1}))
    - (\B{G}(\bQ_{ij} + \mathcal{A}^+_{i-1/2,j} \Delta \bQ) - \B{G}(\bQ_{ij})) 
    = \B{K}_{i,j+1} - \B{K}_{i,j},
\end{multline}
where $\B{K}_{ij} := \B{G}(\bQ_{ij} + \mathcal{A}^+_{i-1/2,j} \Delta
\bQ) - \B{G}(\bQ_{ij})$.  This is precisely analogous to
\eqref{eq:conservation-condition}, and ensures that the transverse
fluctuations will update the solution in a conservative fashion.
Condition \eqref{eq:transverse-conservation-construction} is typically
satisfied by transverse solvers for linear problems (including
our own transverse solvers away from interfaces between distinct
media, which is the only situation where conservation of the elements
of $\bQ$ makes sense for our problem), and
\eqref{eq:transverse-conservation-construction} may be taken as a
recipe for defining conservative transverse fluctuations for problems
for which it is not immediately obvious how to do so.

\subsubsection{Second-order correction term}

These terms modify the basic Godunov-type update of the first line of
\eqref{eq:fvm-update-classic} with an additional antidiffusive flux.
Mathematically, they provide the second-order terms in the Taylor
expansion of the solution in time; qualitatively, they remove the
diffusive error from the Godunov update, sharpening the solution.  In
a classical high-resolution finite volume method, these second-order
correction fluxes are computed as linear combinations of the waves
propagating from each interface, with {\em limiters} applied to each
wave.  These limiters compare the strength of a wave at an interface
to that of the wave in the same family at the neighboring interface in
the upwind direction; based on its relative strength the wave may be
scaled downward in magnitude to avoid overshoot, or amplified to
sharpen the solution.  For a more thorough discussion, see chapter 6
of LeVeque~\cite{rjl:fvm-book}.

There are two problems with implementing appropriate second-order
correction fluxes here.  The first is in formulating the appropriate
numerical flux function on either side of an interface where a
condition of the form \eqref{eq:interface-linear-general} holds.  This
is not trivial, but can be accomplished with, for instance, an
approach based on the Immersed Interface Method~\cite{li-ito:iim-book,
  zhang:elastic, rjl-zhang:acou}.  The second, harder problem is in
creating an appropriate limiting scheme when the waves in the upwind
direction are in a different medium and may be in no way analogous to
the waves we wish to limit.

Because of these problems, we omit the second-order correction term at
interfaces between different materials.  While this reduces the
accuracy of our solution to first order locally, we note that
classical high-resolution finite volume methods also lose formal
accuracy at such interfaces, even if the correction term is included
(Section 9.12 of~\cite{rjl:fvm-book}), and that it is only at interfaces between
different materials that we incur higher error.  At cell interfaces
between identical materials, we use the standard second-order
correction flux.

\subsection{Source term}

As in our previous work~\refcartpaper{}, we handle the source term
$\bD\bQ$ in \eqref{eq:first-order-system} using operator splitting.
Based on the our experience with Godunov and Strang splitting, we use
only Strang splitting here.  The matrix $\bD$ has a simple
eigenstructure, so we use the exact solution operator $\exp(\bD \Delta
t)$ to advance the solution under the action of the dissipation term;
this has the advantages of being the most accurate method available to
apply the source term, and of being unconditinally stable for all
$\Delta t$.  We found in~\refcartpaper{} that even Strang splitting
degrades to first-order accuracy in the stiff regime where the
timestep is longer than the characteristic dissipation time, but for
short timesteps it displays the expected second-order behavior.

\subsection{Numerical software}

We implemented the numerical solution techniques described here using
the \clawpack{} finite volume method package, version
4.6~\cite{claw.org}.  Normally, writing a \clawpack{} application
would require only writing a few plug-in subroutines, such as the
Riemann solver; because of our modified formulation, we also had to
modify some of its internal subroutines, but using \clawpack{} still
greatly reduced the time and effort required for coding compared to
starting from scratch.  The package supports operator splitting for
source terms, such as the dissipative term here, by means of a
user-supplied subroutine that advances the system by a specified time
step under the action of the source term.  Both Godunov and Strang
splitting are available, though we use only Strang splitting here.
Block-structured Berger-Colella-Oliger adaptive mesh refinement (AMR)
is available from the \amrclaw{} package~\cite{mjb-rjl:amrclaw};
\amrclaw{} can also run in parallel on shared-memory systems using
OpenMP.

\section{Results for rectilinear grids}
\label{sec:results-interface}

In this section we test our code's ability to correctly model
interfaces between different poroelastic materials, and between a
poroelastic material and a fluid, on a rectilinear grid.  Note that
although the grid lines are straight, these results still use the
mapped-grid solver; the grid mapping just happens to be particularly
simple.  For all cases in this section, we choose the simulation time
step such that the global maximum CFL number is 0.9.  We also use no
limiters anywhere in these simulations --- while limiters improve
solution accuracy on typical grids, they obscure the convergence
behavior of the underlying wave propagation method, which is what we
seek to observe here.  Since the solutions used in this section are
smooth except at interfaces between different media, we will not
encounter serious trouble from the numerical dispersion that limiters
are designed to suppress.  Table \ref{tab:matprops} lists the properties
of all the materials used here.

\begin{table}[t]
\caption{Properties of the poroelastic media used in test cases.
  Stone properties are taken from de la Puente et
  al.~\cite{delapuente-dumbser-kaser-igel:poro-dg}; bone properties
  are from~\cite{cowin:bone-poro,
    hosokawa-otani:ultrasound-bovine-cancellous,
    hughes-leighton-white:aniso-cancellous,
    smit-huyghe-cowin:cortical-params}.  Cortical bone
    properties refer to vascular pore space.  The sandstone properties
    are orthotropic, but the other materials are isotropic.  Wave
    speeds are correct in the high-frequency limit.  $c_{pf}$ is the
    fast P wave speed, $c_s$ is the S wave speed, $c_{ps}$ is the slow
    P wave speed, and $\tau_d$ is the time constant for dissipation.
    Subscript numbers indicate principal directions.}
\label{tab:matprops}
\begin{center}
\begin{tabular}{rlllll}
  \toprule
  & & Sandstone & Shale & Cortical bone & Cancellous bone \\
  \midrule
  \multicolumn{6}{l}{Base properties} \\
  \midrule
  $K_s$ & (GPa) & 80 & 7.6 & 14 & 18.5 \\
  $\rho_s$ & (kg/m$^3$) & 2500 & 2210 & 1960 & 1960 \\
  $c_{11}$ & (GPa) & 71.8 & 11.9 & 20.6 & 5.2 \\
  $c_{12}$ & (GPa) & 3.2 & 3.96 & 10.6 & 2.4 \\
  $c_{13}$ & (GPa) & 1.2 & 3.96 & 10.6 & 2.4 \\
  $c_{33}$ & (GPa) & 53.4 & 11.9 & 20.6 & 5.2 \\
  $c_{55}$ & (GPa) & 26.1 & 3.96 & 5 & 1.38 \\
  $\phi$ & & 0.2 & 0.16 & 0.04 & 0.75 \\
  $\kappa_1$ & ($10^{-15}$ m$^2$) & 600 & 100 & 630 & $7\times 10^6$ \\
  $\kappa_3$ & ($10^{-15}$ m$^2$) & 100 & 100 & 630 & $7\times 10^6$ \\
  $T_1$ & & 2 & 2 & 2 & 1 \\
  $T_3$ & & 3.6 & 2 & 2 & 1 \\
  $K_f$ & (GPa) & 2.5 & 2.5 & 2.3 & 2.2 \\
  $\rho_f$ & (kg/m$^3$) & 1040 & 1040 & 1060 & 990 \\
  $\eta$ & ($10^{-3}$ kg/m$\cdot$s) & 1 & 1 & 1 & 40 \\ \midrule
  \multicolumn{6}{l}{Derived quantites} \\
  \midrule
  $c_{pf1}$ & (m/s) & 6000 & 2480 & 3290 & 3260 \\
  $c_{pf3}$ & (m/s) & 5260 & 2480 & 3290 & 3260 \\
  $c_{s1}$ & (m/s) & 3480 & 1430 & 1620 & 1680 \\
  $c_{s3}$ & (m/s) & 3520 & 1430 & 1620 & 1680 \\
  $c_{ps1}$ & (m/s) & 1030 & 1130 & 1123 & 1480 \\
  $c_{ps3}$ & (m/s) & 746 & 1130 & 1123 & 1480 \\
  $\tau_{d1}$ & ($\mu$s) & 5.95 & 1.25 & 33 & 92 \\
  $\tau_{d3}$ & ($\mu$s) & 1.82 & 1.25 & 33 & 92 \\
  \bottomrule
\end{tabular}
\end{center}
\end{table}

The results we wish to show in this section are error values and
convergence rates.  Where we measure error relative to the true
solution, we do so using grid energy norms: we take the energy norm
(defined by \eqref{eq:energy-norm-def}) of the difference between the
numerical solution for each grid cell and the true solution evaluated
at the cell centroid.  We then define the 1-norm and max-norm errors
as
\begin{equation} \label{eq:err-norm-rectilinear-def}
  \begin{aligned}
  \text{error 1-norm} &= \frac{1}{N_1 N_2} \sum_{i=1}^{N_1}
  \sum_{j=1}^{N_2} \| \bQ_{ij,\text{numerical}} - \bQ_{ij,\text{true}}
  \|_E\\
 \text{error max-norm} &= \max_{i,j} \| \bQ_{ij,\text{numerical}} -
 \bQ_{ij,\text{true}} \|_E
 \end{aligned}
\end{equation}
where $N_1$ and $N_2$ are the grid dimensions in the two computational
axes.
The
incident waves are all scaled to unit peak energy density, so these
grid energy error norms provide an overall measure of relative error
in the numerical solution, with the various components of the solution
scaled appropriately relative to each other.  All cases in this
section test the ability of the code to correctly evolve a known
analytical solution: we initialize all simulations by evaluating the
true solution at cell centroids, and boundary conditions are
implemented using the standard ghost cell approach, with ghost cells
filled using the true solution evaluated at ghost cell centroids.  In
addition, all the solutions considered in this section are periodic in
time.  We have observed that in such cases the error in the numerical
solution displays a periodic component, which gives it local minima at
integer multiples of the period after the starting time, so to avoid
artificially picking the best result we evaluate the solution error
after 1.25 periods.  Since we run a very large number of simulations,
we do not display convergence using plots of error versus grid size or
spacing.  Instead, for each true solution we perform a linear
least-squares fit of the logarithm of error in each norm with respect
to the logarithm of grid dimension; we report the best, worst, and
mean convergence rates across all solutions, as measured by the slopes
of these fit lines, along with the worst $R^2$ value for the fit, and
the best and worst errors in each norm on the finest grid used.  For
all test cases in this section, we assess convergence using grid
dimensions of $100 \times 100$, $200 \times 200$, $400 \times 400$,
and $800 \times 800$ cells.

Before we treat material interfaces on rectilinear grids, we will
first explore the consequences of omitting the second-order correction
term along a line of cell interfaces.

\subsection{Effect of omitting the second-order correction term}
\label{sec:no2nd}

To see the results of omitting the second-order correction term, we
use time-harmonic simple plane waves propagating in homogeneous media.
The correction term is omitted for the vertical-direction fluxes along
one grid line passing horizontally through the center of the domain.
We examine the results both for poroelastic waves of all three
families --- though with viscosity omitted, to isolate the effect of
the missing term --- and for acoustic waves; the material used is the
orthotropic sandstone of Table \ref{tab:matprops} for the poroelastic
waves, or just the brine contained within it for acoustic waves.  To
test grid alignment and transverse solve effects, we vary the wave
propagation direction from straight down in the $-z$ direction to
$7.5^\circ$ clockwise of the $+x$ direction, in $7.5^\circ$
increments.  Because the sample poroelastic material is orthotropic,
we also vary its principal $1$-direction from horizontal to
$165^\circ$ counterclockwise of horizontal, in $15^\circ$ increments.
We omit the viscous dissipation term for the poroelastic waves, since
we wish to focus on error in the wave-propagation part of the
algorithm; since these tests are inviscid, the period of the wave and
time required to cross the domain are the only time scales present, so
we choose an angular frequency of 1~rad/s.  For all cases, the
computational domain is a square whose side length is two wavelengths
of a wave in the incident family propagating in the material principal
$1$-direction.

\begin{figure}[t]
  \begin{center}
    \subfloat[1-norm \label{fig:rt-no2nd-compare-1norm}]
             {\includegraphics[width=0.4\textwidth]{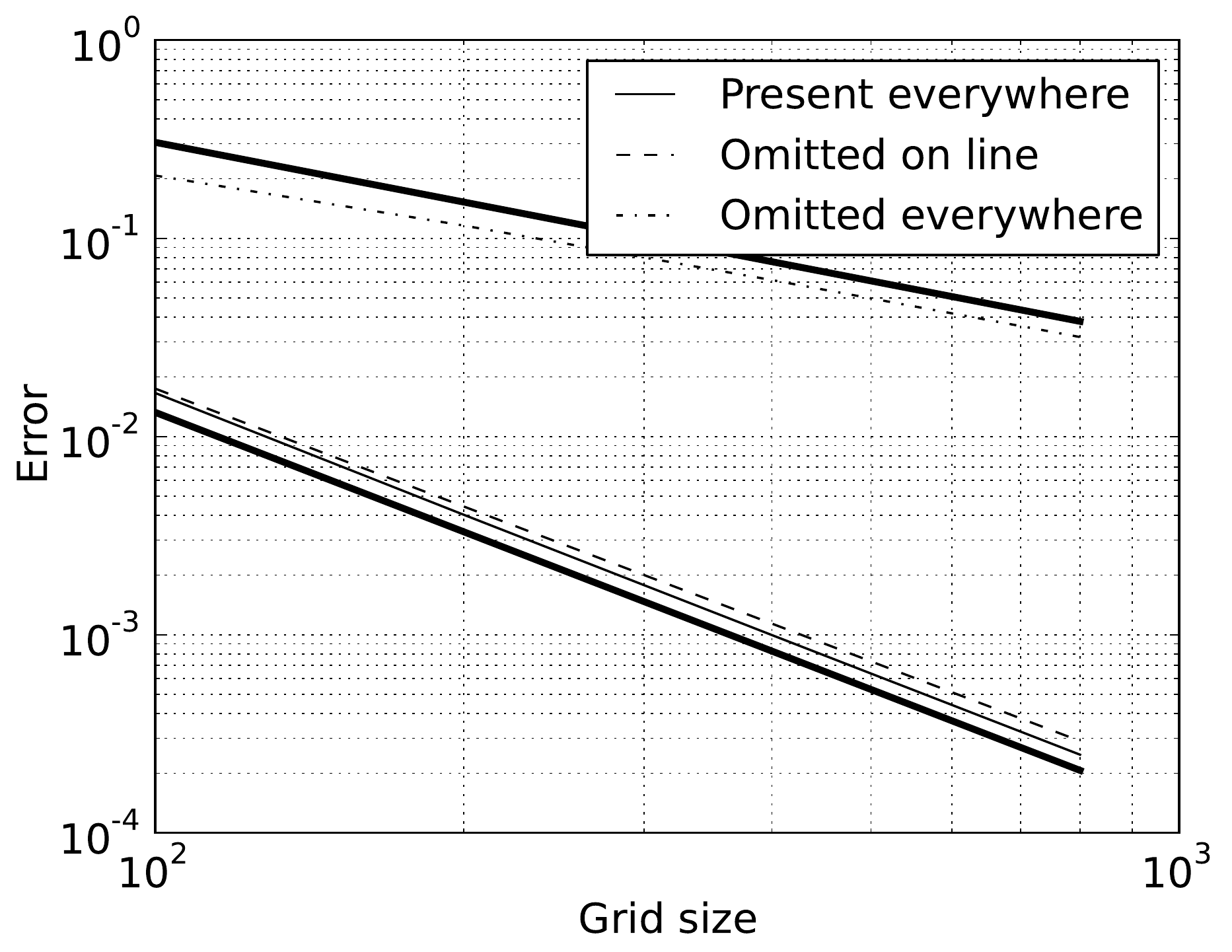}}
    \hspace{0.5in}
    \subfloat[Max-norm \label{fig:rt-no2nd-compare-maxnorm}]
             {\includegraphics[width=0.4\textwidth]{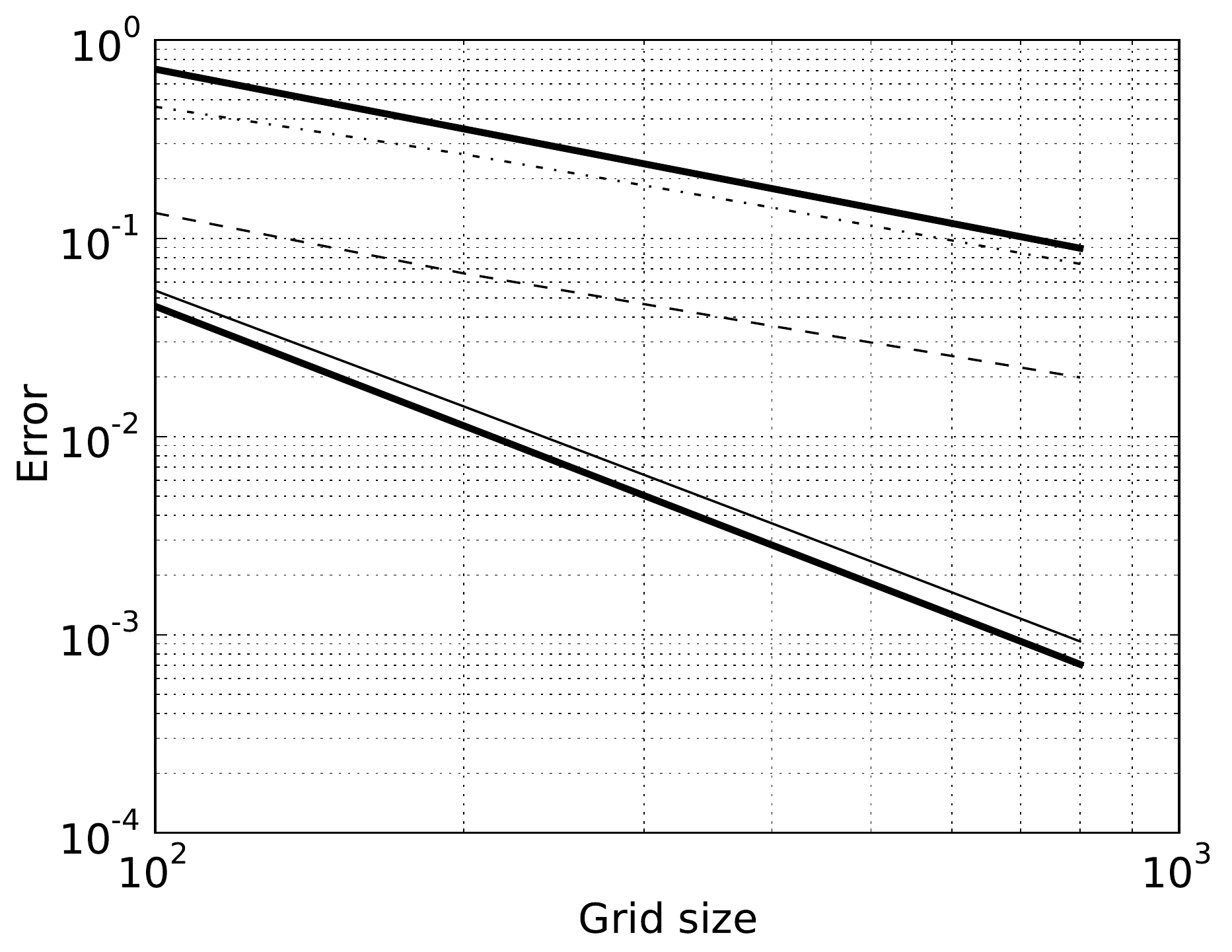}}
  \end{center}
  \caption{Detailed examination of the effect of omitting the
    second-order correction term for the case of Table
    \ref{tab:rt-results-no2nd} that had the highest error on the
    finest grid.  Results of omitting the second-order correction term
    on a line are compared to results if the term is present or
    omitted everywhere.  The thick black lines are first- and
    second-order reference lines. \label{fig:rt-no2nd-compare}}
\end{figure}

Table \ref{tab:rt-results-no2nd} shows the results of this convergence
study.  Even with the second-order correction term omitted along a
line, we still obtain second-order convergence in the 1-norm.
Convergence degrades to first-order (typically) or somewhat below (in
the worst case) for the max-norm.  These results are about as good as
could be expected from omitting the second-order term.  Figure
\ref{fig:rt-no2nd-compare} examines in greater
detail the case that gives the worst error on the $800 \times 800$
grid --- a slow P wave propagating vertically downward, with the
principal direction of the poroelastic sandstone medium oriented
$15^\circ$ counterclockwise of horizontal.  In the 1-norm, results are
almost as good as if the second-order correction term were present
everywhere, while in the max-norm, even though the convergence rate is
reduced to first-order when the term is omitted along a line, the
magnitude of the error is almost four times less than if it were
omitted everywhere.  This suggests that we can expect results of
reasonable accuracy even with the second-order correction term omitted
at interfaces between different media.

\begin{table}[t]
\caption{Convergence results for acoustic and inviscid poroelastic
  waves with the second-order correction term omitted along a line of
  cell interfaces.}
\label{tab:rt-results-no2nd}
\begin{center}
\begin{tabular}{cccccccc}
  \toprule
  & & \multicolumn{3}{c}{Convergence rate} & & \multicolumn{2}{c}{Error on $800 \times 800$ grid}\\
  \cmidrule(r){3-5} \cmidrule(r){7-8}
  & Error norm & Best & Worst & Mean & Worst $R^2$ value & Best & Worst \\
  \midrule
  \multirow{2}{*}{Acoustic} & 1-norm & 2.01 & 1.99 & 2.00 & 0.99957 & $2.57 \times 10^{-5}$ & $4.09 \times 10^{-5}$ \\
  & Max-norm & 1.49 & 1.03 & 1.15 & 0.95837 & $2.76 \times 10^{-4}$ & $2.51 \times 10^{-3}$ \\
  \midrule
  \multirow{2}{*}{Fast P} & 1-norm & 2.02 & 1.97 & 1.99 & 0.99883 & $2.75 \times 10^{-5}$ & $7.83 \times 10^{-5}$ \\
  & Max-norm & 1.45 & 0.91 & 1.13 & 0.86779 & $4.62 \times 10^{-4}$ & $3.74 \times 10^{-3}$ \\
  \midrule
  \multirow{2}{*}{S} & 1-norm & 2.01 & 1.98 & 1.99 & 0.99999 & $4.18 \times 10^{-5}$ & $9.49 \times 10^{-5}$ \\
  & Max-norm & 1.76 & 0.95 & 1.15 & 0.96214 & $4.05 \times 10^{-4}$ & $5.06 \times 10^{-3}$ \\
  \midrule
  \multirow{2}{*}{Slow P} & 1-norm & 2.01 & 1.95 & 1.98 & 0.99999 & $8.21 \times 10^{-5}$ & $3.10 \times 10^{-4}$ \\
  & Max-norm & 1.85 & 0.76 & 1.09 & 0.98524 & $9.91 \times 10^{-4}$ & $1.99 \times 10^{-2}$ \\
  \bottomrule
\end{tabular}
\end{center}
\end{table}

\subsection{Reflected and transmitted waves at a fluid-poroelastic
  interface}
\label{sec:rt-fluid-poro}

\begin{figure}[t]
  \begin{center}
    \subfloat[Fluid-poroelastic \label{fig:rt-cartoon-fluid-poro}]
             {\includegraphics{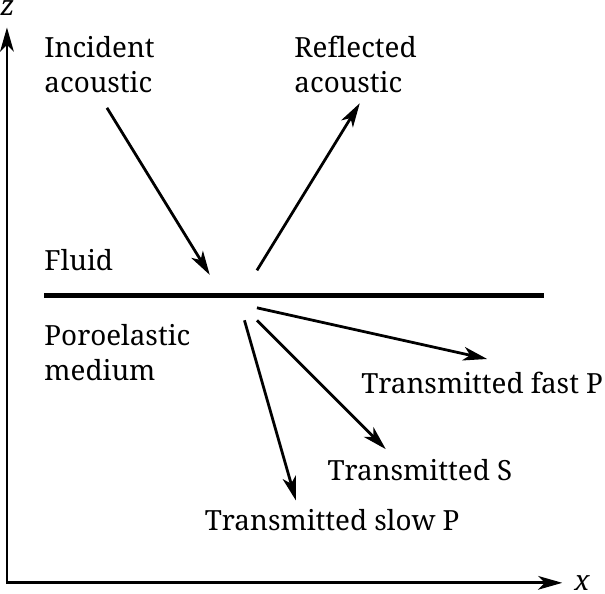}}
    \hspace{1in}
    \subfloat[Poroelastic-poroelastic \label{fig:rt-cartoon-poro-poro}]
             {\includegraphics{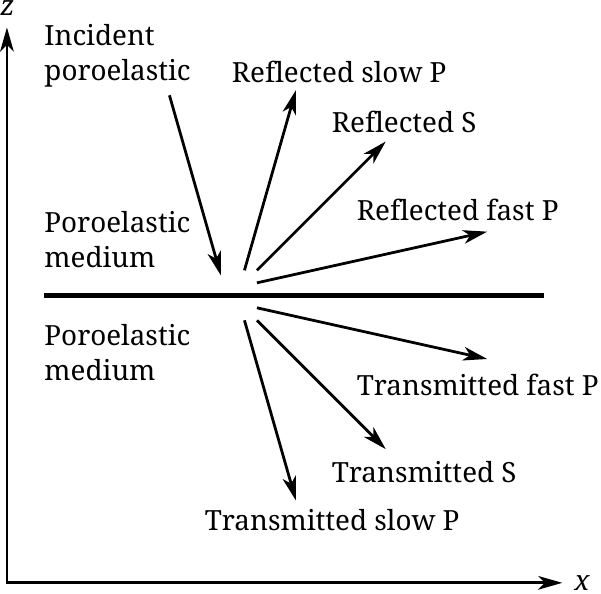}}
    \caption{Simple sketches of reflection/transmission problems used
      to test convergence.  The incident, reflected, and transmitted
      waves are labeled by type; an incident poroelastic wave may be
      of any of the three types. \label{fig:rt-cartoon}}
  \end{center}
\end{figure}

To test our code's ability to correctly handle material interfaces, we
first examine the case of a time-harmonic plane acoustic wave striking
a flat interface between a poroelastic medium and a fluid.  In all
cases, the interface between the two media is horizontal, and the
incident wave impinges from the top side.  Figure
\ref{fig:rt-cartoon-fluid-poro} shows a simple sketch of the problem.
We generate analytical solutions for these cases using the procedure
described in Appendix \ref{sec:reflect-transmit}.  The poroelastic medium used is the brine-saturated
sandstone of Table \ref{tab:matprops}, and the fluid medium is the
same brine contained in the sandstone.  We again perform a
combinatorial sweep over the relevant parameters.  In this case, we
vary the direction of the incident wave from $7.5^\circ$ below
horizontal to straight down in $7.5^\circ$ increments, and the
principal material $1$-direction from horizontal to $165^\circ$
counterclockwise of horizontal in $15^\circ$ increments.  For each
combination of incident wave and material principal directions, we
simulate with interface discharge efficiencies $\eta_d$ of $0$ (sealed
pores), $0.5$ (imperfect hydraulic contact), and $1$ (open pores).
The dimensions of the domain are two wavelengths of an acoustic wave
at the chosen frequency.  We perform these tests both with viscosity
ignored, in order to investigate the performance of the hyperbolic
solver by itself, and with viscosity included.  For the inviscid
tests, since the period of the wave and time required for it to cross the
domain are the only time scales present, the wave frequency is not
directly relevant to the solution error, so we chose an angular
frequency of $1$~rad/s.  With viscosity included, the dissipation term
has its own intrinsic time scale indepdendent of the wave behavior, so
choice of frequency is important.  For these viscous tests, we
restrict ourselves to a high frequency (10~kHz; the maximum frequency
for validity of Biot theory in the sandstone medium is 25~kHz),
because at low frequencies the slow P wave dissipates over distances
much shorter than a wavelength.  In fact, the characteristic decay
length is much shorter than the grid spacing on an otherwise
reasonably resolved grid, so at low frequencies the transmitted slow P
wave could not be resolved without refining the grid to impractical
levels.  This also results in grid sizes and corresponding time steps
short enough that these models are outside the stiff regime identified
in~\refcartpaper{}.  For all cases, the computational domain is again
a square with side length equal to two wavelengths of incident
acoustic wave.

Tables \ref{tab:rt-results-basic-invisc} and
\ref{tab:rt-results-basic-visc-hifreq} show the results of these
convergence studies for the inviscid and viscous high-frequency tests
respectively.  We obtain second-order convergence for all cases in the
1-norm, though results are degraded to first-order in the max-norm due
to the omission of the second-order correction term at the interface,
as expected from the results of the previous section.  All three
interface conditions (open and sealed pores, and imperfect hydraulic
contact) show similar levels of error, and the difference between best
and worst error obtained in the finest grid is also not large.
Convergence rates are slightly worse for the viscous cases; they are
still typically quite close to the expected second-order and
first-order results, but some outlying cases degrade as far as order
0.7, possibly due to the compounding of first-order error at the
interface with operator splitting error.

\begin{table}[t]
\caption{Convergence results for an acoustic wave impinging on
  sandstone, with viscosity omitted.}
\label{tab:rt-results-basic-invisc}
\begin{center}
\begin{tabular}{cccccccc}
  \toprule
  & & \multicolumn{3}{c}{Convergence rate} & &
  \multicolumn{2}{c}{Error on $800 \times 800$ grid}\\
  \cmidrule(r){3-5} \cmidrule(r){7-8}
  & Error norm & Best & Worst & Mean & Worst $R^2$ value & Best & Worst \\
  \midrule
  \multirow{2}{*}{$\eta_d=0$} & 1-norm & 2.00 & 1.98 & 1.99 & 1.00000 & $7.93 \times 10^{-5}$ & $2.09 \times 10^{-4}$ \\
  & Max-norm & 1.50 & 0.88 & 1.12 & 0.96350 & $2.17 \times 10^{-3}$ & $1.25 \times 10^{-2}$ \\
  \midrule
  \multirow{2}{*}{$\eta_d=0.5$} & 1-norm & 2.01 & 1.97 & 1.99 & 0.99999 & $1.03 \times 10^{-4}$ & $2.94 \times 10^{-4}$ \\
  & Max-norm & 1.47 & 0.88 & 1.14 & 0.95713 & $2.12 \times 10^{-3}$ & $1.24 \times 10^{-2}$ \\
  \midrule
  \multirow{2}{*}{$\eta_d=1$} & 1-norm & 2.01 & 1.97 & 2.00 & 0.99999 & $1.07 \times 10^{-4}$ & $3.14 \times 10^{-4}$ \\
  & Max-norm & 1.52 & 0.88 & 1.15 & 0.95981 & $2.16 \times 10^{-3}$ & $1.24 \times 10^{-2}$ \\
  \bottomrule
\end{tabular}
\end{center}
\end{table}

\begin{table}[t]
\caption{Convergence results for a $10$ kHz acoustic wave impinging on
  sandstone, with viscosity included.}
\label{tab:rt-results-basic-visc-hifreq}
\begin{center}
\begin{tabular}{cccccccc}
  \toprule
  & & \multicolumn{3}{c}{Convergence rate} & &
  \multicolumn{2}{c}{Error on $800 \times 800$ grid}\\
  \cmidrule(r){3-5} \cmidrule(r){7-8}
  & Error norm & Best & Worst & Mean & Worst $R^2$ value &
  Best & Worst \\
  \midrule
  \multirow{2}{*}{$\eta_d=0$} & 1-norm & 1.99 & 1.96 & 1.98 & 1.00000 & $6.93 \times 10^{-5}$ & $1.21 \times 10^{-4}$ \\
  & Max-norm & 1.44 & 0.88 & 1.09 & 0.96397 & $2.46 \times 10^{-3}$ & $1.25 \times 10^{-2}$ \\
  \midrule
  \multirow{2}{*}{$\eta_d=0.5$} & 1-norm & 1.97 & 1.90 & 1.95 & 0.99999 & $6.84 \times 10^{-5}$ & $1.79 \times 10^{-4}$ \\
  & Max-norm & 1.30 & 0.72 & 0.98 & 0.97612 & $2.96 \times 10^{-3}$ & $2.06 \times 10^{-2}$ \\
  \midrule
  \multirow{2}{*}{$\eta_d=1$} & 1-norm & 1.97 & 1.90 & 1.94 & 0.99999 & $6.75 \times 10^{-5}$ & $1.83 \times 10^{-4}$ \\
  & Max-norm & 1.29 & 0.70 & 0.98 & 0.97672 & $3.02 \times 10^{-3}$ & $2.19 \times 10^{-2}$ \\
  \bottomrule
\end{tabular}
\end{center}
\end{table}

\subsection{Reflected and transmitted waves at an interface between
  two poroelastic materials}
\label{sec:rt-poro-poro}

We also examine plane waves reflected and transmitted at an interface
between two poroelastic media.  For these cases the upper medium is
the isotropic, brine-saturated shale of Table \ref{tab:matprops}, and
the lower medium is again orthotropic sandstone; Figure
\ref{fig:rt-cartoon-poro-poro} shows a simple sketch of the reflected
and transmitted waves in this case.  These cases are
quite similar to the fluid-poroelastic cases, with the exception that
we send in incident waves in multiple different poroelastic wave
families.  To reduce the total number of cases, we run a smaller
number of incident wave directions --- we vary the angle of incidence
from $7.5^\circ$ below horizontal to $82.5^\circ$ in steps of
$15^\circ$, which halves the number of angles of incidence and drops
the relatively uninteresting normal incidence, while retaining the
grazing $7.5^\circ$ angle.  We use the same set of principal material
directions for the sandstone as in the fluid-poroelastic cases;
principal material direction is irrelevant for the shale because it is
isotropic.  We also use the same sets of interface discharge
efficiencies and grid dimensions.  For all cases the domain size is
two wavelengths of the fast P wave in shale at the chosen frequency.
The inviscid cases use waves with an angular frequency of 1~rad/s,
while the viscous cases use 20~kHz waves, a frequency high enough to
be out of the stiff regime (or just at the edge of it in shale for a
$100 \times 100$ grid) but low enough for low-frequency Biot theory to
be valid for both materials.  We simulate incident waves in all three
families for the inviscid cases in order to exercise all different
possibilities in the solution code, but with viscosity present we do
not include incident slow P waves --- their decay rate is high enough
that their amplitude is reduced by a factor of $10^{20}$ over the width
of the domain, making these cases both intractable to simulate and
unrealistic.

Tables \ref{tab:rt-results-poro-invisc} and
\ref{tab:rt-results-poro-visc} show the results of these convergence
studies.  Similarly to the acoustic-poroelastic cases, we obtain
consistent second-order convergence in the 1-norm, and first-order or
better convergence in the max-norm.  This occurs for all incident wave
types, even in the face of the operator splitting error present in the
viscous cases.  In fact, the worst-case convergence rates are better
than for the acoustic-poroelastic cases, perhaps due to a more
tractable interface condition.

\begin{table}[t]
\caption{Convergence results for poroelastic waves in shale impinging
  on sandstone, with viscosity omitted.}
\label{tab:rt-results-poro-invisc}
\begin{center}
\begin{tabular}{ccccccccc}
  \toprule
  & \multirow{2}{*}[-3.5pt]{\parbox{0.5in}{\centering Incident wave}} & & \multicolumn{3}{c}{Convergence rate} & \multirow{2}{*}[-3.25pt]{\parbox{0.6in}{\centering Worst $R^2$ value}} & \multicolumn{2}{c}{Error on $800 \times 800$ grid}\\
  \cmidrule(r){4-6} \cmidrule(r){8-9}
  & & Error norm & Best & Worst & Mean & & Best & Worst \\
  \midrule
  \multirow{6}{*}[-4.5pt]{$\eta_d=0$} & \multirow{2}{*}{Fast P} & 1-norm & 2.01 & 1.99 & 2.00 & 1.00000 & $8.81 \times 10^{-5}$ & $3.26 \times 10^{-4}$ \\
  & & Max-norm & 1.49 & 1.06 & 1.17 & 0.97570 & $2.01 \times 10^{-3}$ & $5.62 \times 10^{-3}$ \\
  \cmidrule(r){2-9}
  & \multirow{2}{*}{S} & 1-norm & 2.01 & 1.98 & 1.99 & 0.99999 & $2.19 \times 10^{-4}$ & $4.27 \times 10^{-4}$ \\
  & & Max-norm & 1.86 & 1.03 & 1.25 & 0.98292 & $2.39 \times 10^{-3}$ & $1.18 \times 10^{-2}$ \\
  \cmidrule(r){2-9}
  & \multirow{2}{*}{Slow P} & 1-norm & 2.00 & 1.99 & 2.00 & 0.99999 & $3.34 \times 10^{-4}$ & $6.34 \times 10^{-4}$ \\
  & & Max-norm & 1.92 & 1.08 & 1.40 & 0.99555 & $3.90 \times 10^{-3}$ & $2.13 \times 10^{-2}$ \\
  \midrule
  \multirow{6}{*}[-4.5pt]{$\eta_d=0.5$} & \multirow{2}{*}{Fast P} & 1-norm & 2.00 & 1.98 & 1.99 & 1.00000 & $7.65 \times 10^{-5}$ & $2.00 \times 10^{-4}$ \\
  & & Max-norm & 1.22 & 1.01 & 1.13 & 0.99179 & $1.99 \times 10^{-3}$ & $5.63 \times 10^{-3}$ \\
  \cmidrule(r){2-9}
  & \multirow{2}{*}{S} & 1-norm & 2.01 & 1.98 & 1.99 & 0.99999 & $2.13 \times 10^{-4}$ & $3.49 \times 10^{-4}$ \\
  & & Max-norm & 1.86 & 0.99 & 1.24 & 0.98301 & $2.38 \times 10^{-3}$ & $1.18 \times 10^{-2}$ \\
  \cmidrule(r){2-9}
  & \multirow{2}{*}{Slow P} & 1-norm & 2.01 & 1.99 & 2.00 & 0.99999 & $4.70 \times 10^{-4}$ & $1.44 \times 10^{-3}$ \\
  & & Max-norm & 1.99 & 1.06 & 1.40 & 0.96976 & $3.81 \times 10^{-3}$ & $1.96 \times 10^{-2}$ \\
  \midrule
  \multirow{6}{*}[-4.5pt]{$\eta_d=1$}& \multirow{2}{*}{Fast P} & 1-norm & 2.00 & 1.99 & 1.99 & 1.00000 & $7.66 \times 10^{-5}$ & $1.99 \times 10^{-4}$ \\
  & & Max-norm & 1.22 & 1.00 & 1.13 & 0.99181 & $1.99 \times 10^{-3}$ & $5.63 \times 10^{-3}$ \\
  \cmidrule(r){2-9}
  & \multirow{2}{*}{S} & 1-norm & 2.01 & 1.98 & 1.99 & 0.99999 & $2.11 \times 10^{-4}$ & $3.49 \times 10^{-4}$ \\
  & & Max-norm & 1.86 & 0.98 & 1.24 & 0.98301 & $2.38 \times 10^{-3}$ & $1.18 \times 10^{-2}$ \\
  \cmidrule(r){2-9}
  & \multirow{2}{*}{Slow P} & 1-norm & 2.01 & 1.99 & 2.00 & 0.99999 & $4.86 \times 10^{-4}$ & $1.50 \times 10^{-3}$ \\
  & & Max-norm & 1.99 & 1.05 & 1.41 & 0.97068 & $3.86 \times 10^{-3}$ & $1.98 \times 10^{-2}$ \\
  \bottomrule
\end{tabular}
\end{center}
\end{table}

\begin{table}[t]
\caption{Convergence results for 20~kHz poroelastic waves in shale
  impinging on sandstone, with viscosity included.}
\label{tab:rt-results-poro-visc}
\begin{center}
\begin{tabular}{ccccccccc}
  \toprule
  & \multirow{2}{*}[-3.5pt]{\parbox{0.5in}{\centering Incident wave}} & & \multicolumn{3}{c}{Convergence rate} & \multirow{2}{*}[-3.25pt]{\parbox{0.6in}{\centering Worst $R^2$ value}} & \multicolumn{2}{c}{Error on $800 \times 800$ grid}\\
  \cmidrule(r){4-6} \cmidrule(r){8-9}
  & & Error norm & Best & Worst & Mean & & Best & Worst \\
  \midrule
  \multirow{4}{*}{$\eta_d = 0$} & \multirow{2}{*}{Fast P} & 1-norm & 2.00 & 1.96 & 1.99 & 1.00000 & $6.97 \times 10^{-5}$ & $1.80 \times 10^{-4}$ \\
  & & Max-norm & 1.21 & 1.01 & 1.12 & 0.99164 & $1.98 \times 10^{-3}$ & $5.88 \times 10^{-3}$ \\
  \cmidrule(r){2-9}
  & \multirow{2}{*}{S} & 1-norm & 2.01 & 1.98 & 1.99 & 1.00000 & $2.10 \times 10^{-4}$ & $4.40 \times 10^{-4}$ \\
  & & Max-norm & 1.86 & 1.02 & 1.25 & 0.97550 & $2.56 \times 10^{-3}$ & $1.20 \times 10^{-2}$ \\
  \midrule
  \multirow{4}{*}{$\eta_d = 0.5$} & \multirow{2}{*}{Fast P} & 1-norm & 2.00 & 1.97 & 1.99 & 1.00000 & $6.99 \times 10^{-5}$ & $1.81 \times 10^{-4}$ \\
  & & Max-norm & 1.21 & 1.01 & 1.12 & 0.99168 & $1.98 \times 10^{-3}$ & $5.78 \times 10^{-3}$ \\
  \cmidrule(r){2-9}
  & \multirow{2}{*}{S} & 1-norm & 2.01 & 1.99 & 1.99 & 0.99999 & $2.09 \times 10^{-4}$ & $4.43 \times 10^{-4}$ \\
  & & Max-norm & 1.86 & 1.02 & 1.25 & 0.97509 & $2.53 \times 10^{-3}$ & $1.20 \times 10^{-2}$ \\
  \midrule
  \multirow{4}{*}{$\eta_d = 1$} & \multirow{2}{*}{Fast P} & 1-norm & 2.00 & 1.97 & 1.99 & 1.00000 & $7.00 \times 10^{-5}$ & $1.81 \times 10^{-4}$ \\
  & & Max-norm & 1.21 & 1.01 & 1.12 & 0.99169 & $1.98 \times 10^{-3}$ & $5.78 \times 10^{-3}$ \\
  \cmidrule(r){2-9}
  & \multirow{2}{*}{S} & 1-norm & 2.01 & 1.99 & 1.99 & 0.99999 & $2.09 \times 10^{-4}$ & $4.43 \times 10^{-4}$ \\
  & & Max-norm & 1.86 & 1.02 & 1.25 & 0.97510 & $2.53 \times 10^{-3}$ & $1.20 \times 10^{-2}$ \\
  \bottomrule
\end{tabular}
\end{center}
\end{table}

\section{Results for curved interfaces and mapped grids}
\label{sec:results-mapped}

Now that we have examined the convergence behavior of our code for
problems with interfaces on rectilinear grids, we turn to curvilinear
mapped grids.  Here we treat two types of problem.  First, we perform a
convergence study for a time-harmonic plane wave scattering off an
isotropic poroelastic cylinder --- a case for which we have an
analytic solution --- and second, we simulate an acoustic pulse
striking a simplified model of a human femur bone.

\subsection{Cylindrical scatterer}
\label{sec:results-scatterer}

We model a cylindrical scatterer composed of the isotropic shale of
Table \ref{tab:matprops}.  Our reference solution is
for a time-harmonic acoustic plane wave scattering off a cylindrical
isotropic poroelastic body; we obtained it by transforming the
poroelastic-acoustic system into a set of coupled Helmholtz equations,
then expressing their solutions as bi-infinite series of Bessel and
Hankel function modes.  This solution is very similar to that of
Laperre and Thys~\cite{laperre-thys:cylindrical-scatterer}, although
we extend it to imperfect hydraulic contact in addition to open or
closed pores.  It has a variety of interesting properties, including
resonance-like behavior at certain frequencies; we plan a companion
publication exploring its properties in more detail.  For now, we
confine ourselves to a specific model --- a shale cylinder 5~cm in
diameter, in a bath of fluid identical to its pore fluid.

\begin{figure}[t]
  \begin{center}
    \includegraphics[width=0.6\textwidth]{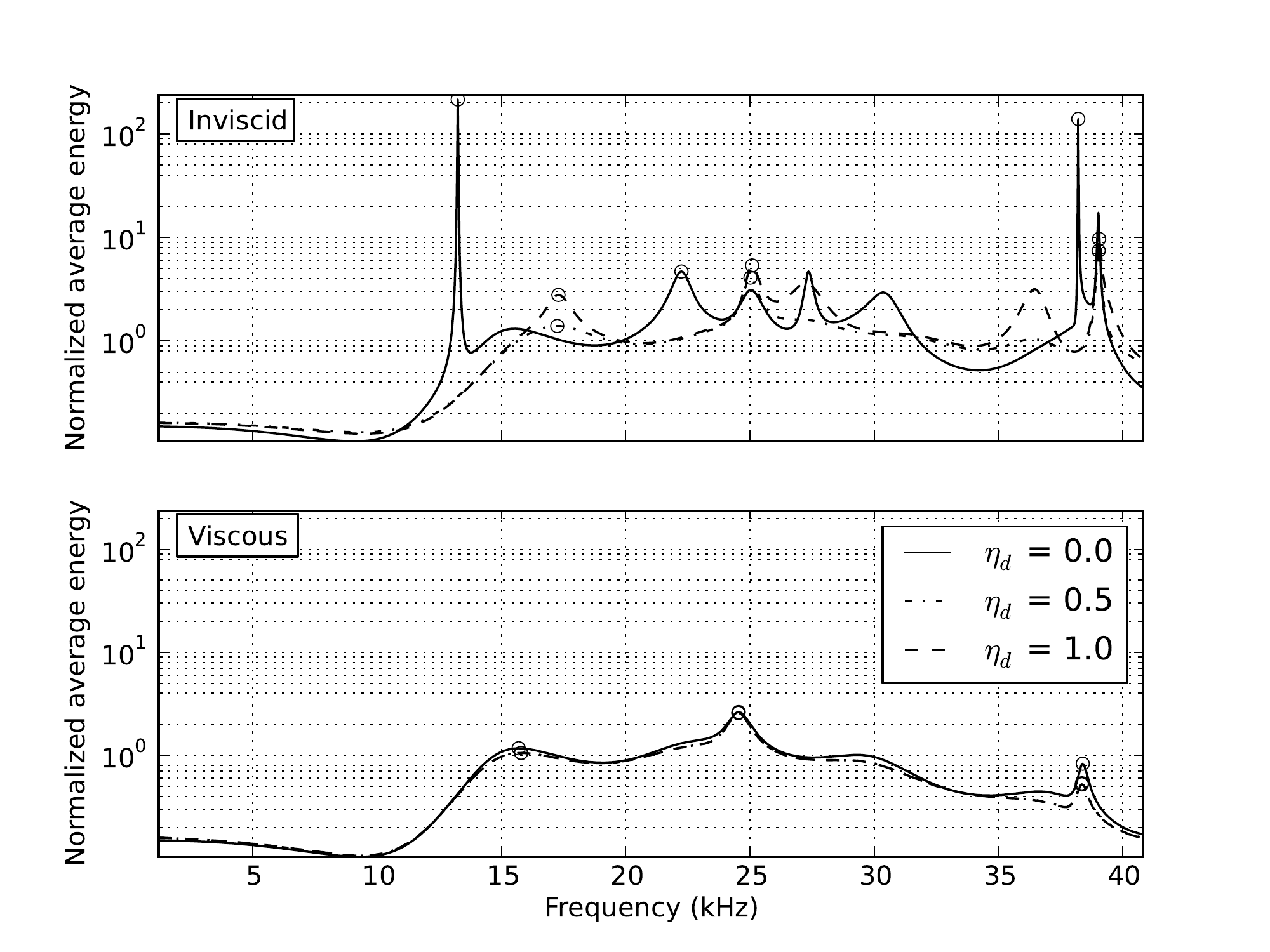}
    \caption{Frequency response plot for 5~cm diameter cylindrical
      scatterer composed of shale, with and without viscosity, for
      several interface discharge efficiencies $\eta_d$.  The
      vertical axis is ratio of the average energy contained within
      the cylinder over one cycle to the average energy contained
      within an equal volume of the undisturbed incident wave.  The
      circles indicate cases chosen for convergence evaluation of the
      mapped-grid poroelastic/fluid code.  This plot goes up to a
      maximum frequency of one-third of the cutoff frequency for
      validity of low-frequency Biot theory; there is no barrier to
      modeling all the way up to the cutoff frequency, but the
      response curve becomes increasingly complicated beyond the
      portion shown here. \label{fig:sc-freq-response}}
  \end{center}
\end{figure}

\begin{figure}[t]
  \begin{center}
    \subfloat[Inviscid, $\eta_d = 0$, 13.25 kHz, analytical]
             {\includegraphics[width=0.4\textwidth]{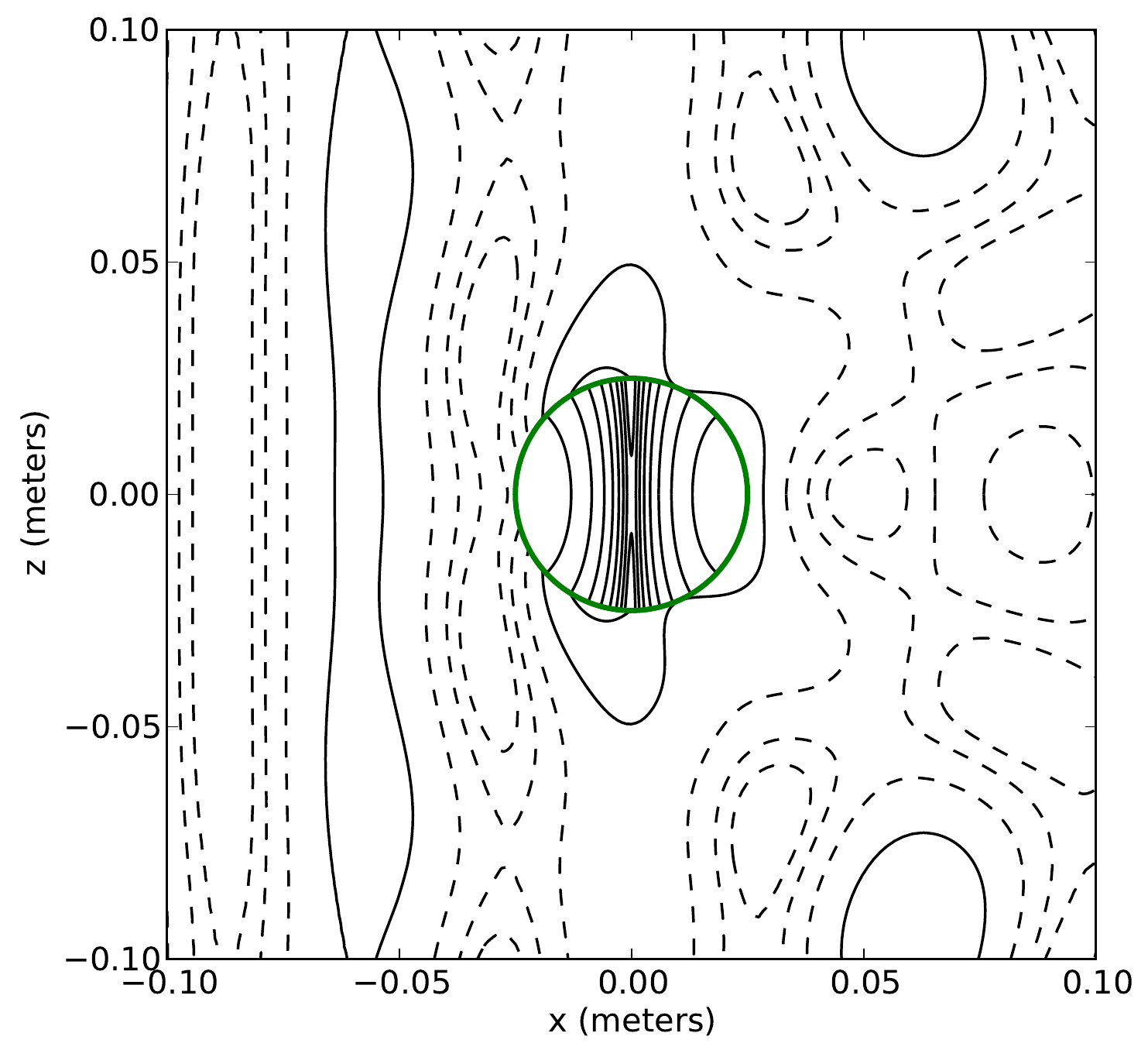}}
    \hspace{0.5in}
    \subfloat[Viscous, $\eta_d = 1$, 38.35 kHz, analytical]
             {\includegraphics[width=0.4\textwidth]{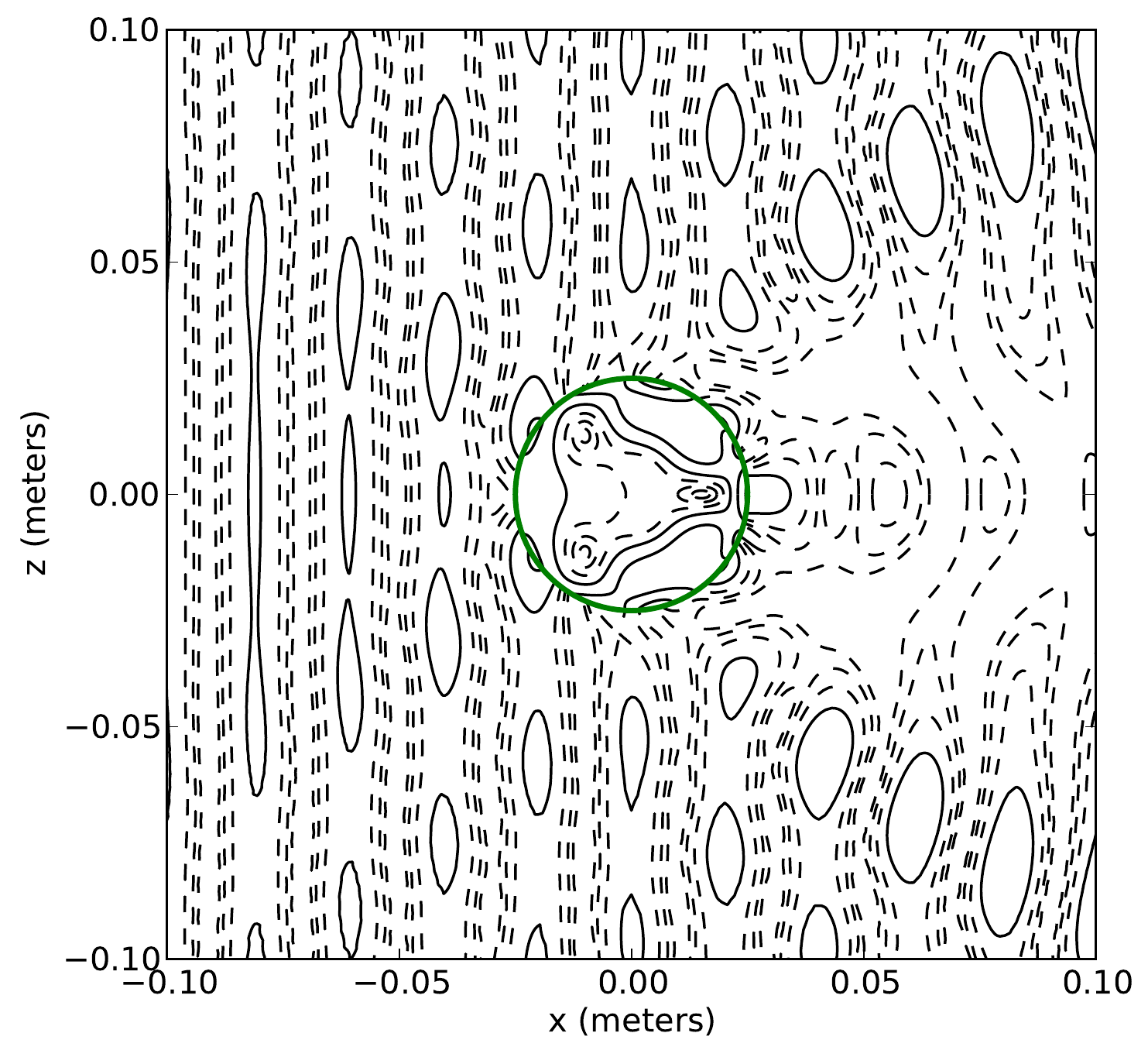}}\\
    \subfloat[Inviscid, $\eta_d = 0$, 13.25 kHz, numerical]
             {\includegraphics[width=0.4\textwidth]{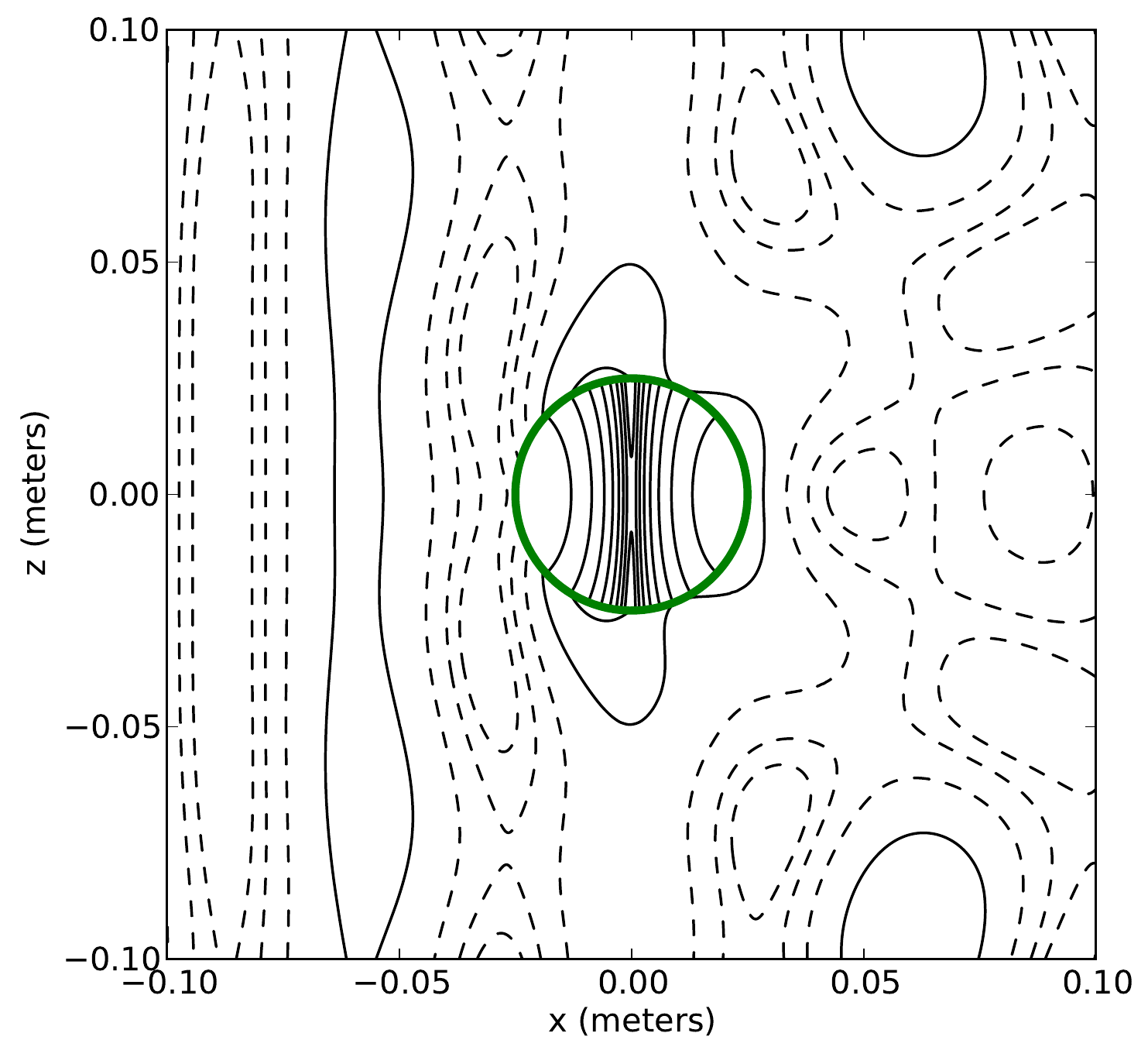}}
    \hspace{0.5in}
    \subfloat[Viscous, $\eta_d = 1$, 38.35 kHz, numerical]
             {\includegraphics[width=0.4\textwidth]{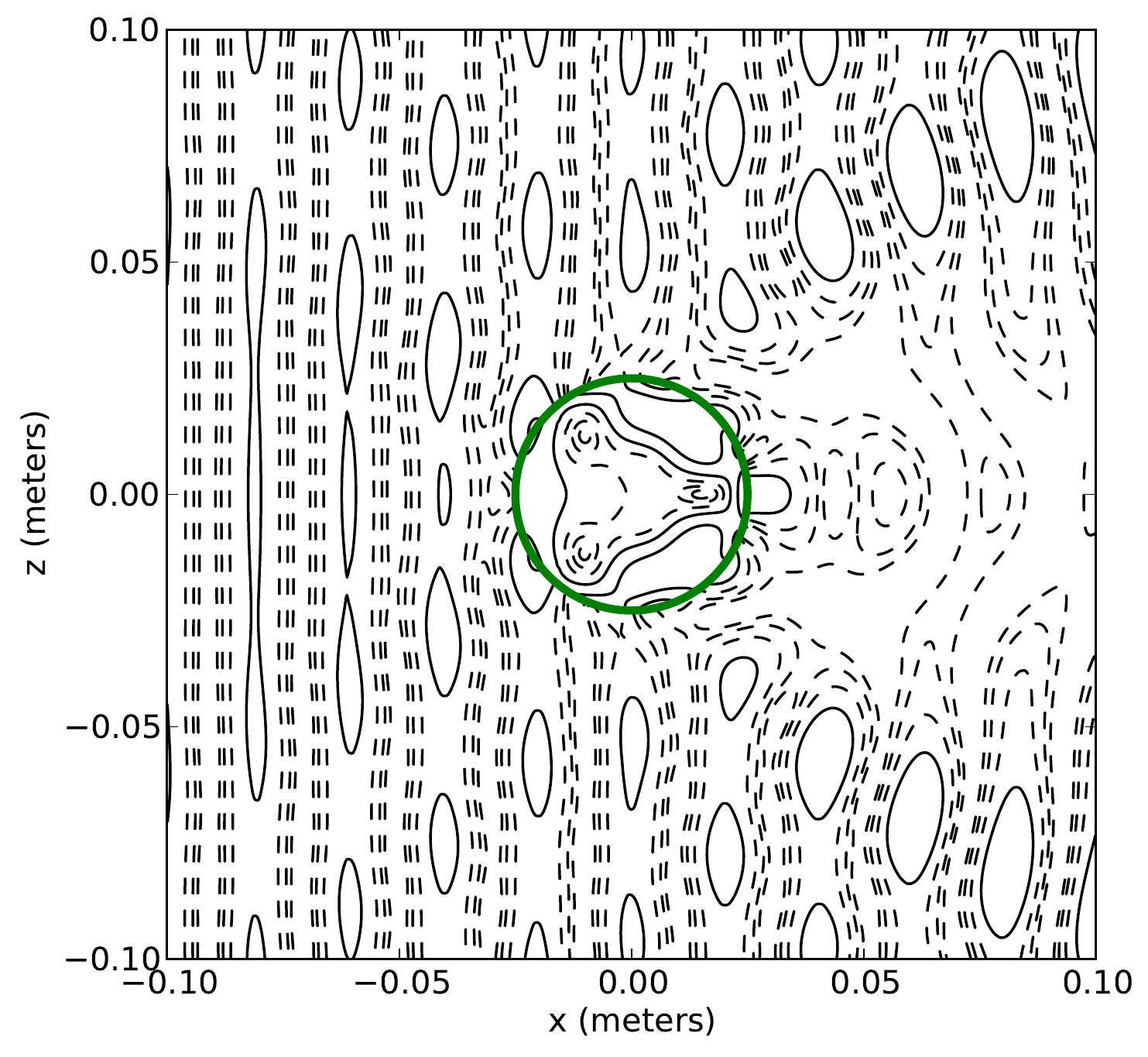}}\\
    \caption{Contour plots of energy density for the analytical
      solutions of selected test cases from Table
      \ref{tab:sc-case-list} at time $t = 0$, and the numerical
      solution on a $1024 \times 1024$ grid after one cycle.  Values
      are normalized by peak energy density of incident wave.
      Contours are placed at every power of 2, with dashed contours
      indicating negative powers; the lowest solid contour is at unit
      normalized energy density, with the adjacent dashed contour at
      $\frac{1}{2}$ and the adjacent solid contour at 2.  The thick
      circle is the boundary of the cylinder.  Maximum normalized
      energy densities inside the cylinder in the upper-left plot are
      967 at the left edge, and 929 at the right.  The numerical
      solutions closely match the analytical ones.  Similar plots for
      the analytical solutions for all cases can be found in Appendix
      \ref{sec:sc-gallery}. \label{fig:sc-case-portraits}}
  \end{center}
\end{figure}

In order to choose sensible input frequencies in the face of the
complex behaivor of the analytical solution, we first examine the
average energy contained in the cylinder over one cycle as a function
of frequency, normalized by the average energy contained in the
incident wave over the same volume and time.  The required integral of
the solution in the radial direction is not readily computable
analytically, so we use Romberg quadrature instead, applied
recursively until the relative difference between the highest-order
iterates at successive depths is less than $10^{-9}$.  We also
truncate the series solution to 35 terms in each direction, by which
index the series terms have decreased to $10^{-9}$ or less of their
typical magnitudes for low index at all frequencies below the Biot
low-frequency validity cutoff, and are decaying rapidly.  Figure
\ref{fig:sc-freq-response} shows this normalized average energy
content with and without viscosity, for open, sealed, and imperfect
pore conditions.  Based on this frequency response plot, we selected
three frequencies that showed strong response for each combination of
viscosity and interface permeability to test the convergence of our
code.  We chose these strongly-responding frequencies because they
would be a more rigorous test of our numerical model: the solution is
large within the cylinder, coinciding with the most distorted grid
cells, and accurate handling of the transfer of energy and momentum
across the interface is likely to be more important.  Table
\ref{tab:sc-case-list} lists the cases chosen, and Figure
\ref{fig:sc-case-portraits} plots the energy density at the initial
time for two selected cases.  Similar plots for all cylindrical
scatterer cases can be found in Appendix \ref{sec:sc-gallery}.

\begin{table}[t]
  \caption{Cylindrical scatterer cases chosen for convergence
    analysis.  Note that the reference to ``high'' frequency here is
    relative to the scale of Figure \ref{fig:sc-freq-response}; all
    cases are well below the cutoff frequency for validity of
    low-frequency Biot theory.}
  \label{tab:sc-case-list}
  \begin{center}
    \begin{tabular}{cC{0.5in}C{0.5in}C{0.5in}p{0.2in}C{0.5in}C{0.5in}C{0.5in}}    
      \toprule
      & \multicolumn{3}{c}{Inviscid case frequency (kHz)} & &
      \multicolumn{3}{c}{Viscous case frequency (kHz)}\\
      \cmidrule(r){2-4} \cmidrule(r){6-8}
      & Low & Mid & High & & Low & Mid & High\\
      \midrule
      $\eta_d = 0$   & 13.25 & 22.25 & 38.20 & & 15.70 & 24.55 & 38.39 \\
      $\eta_d = 0.5$ & 17.25 & 25.02 & 39.03 & & 15.80 & 24.54 & 38.35 \\
      $\eta_d = 1$   & 17.30 & 25.09 & 39.04 & & 15.80 & 24.54 & 38.35 \\
      \bottomrule
    \end{tabular}
  \end{center}
\end{table}

\begin{figure}[t]
  \begin{center}
    \subfloat[Cylindrical scatterer \label{fig:sc-map-figure}]
             {\includegraphics[width=0.5\textwidth]{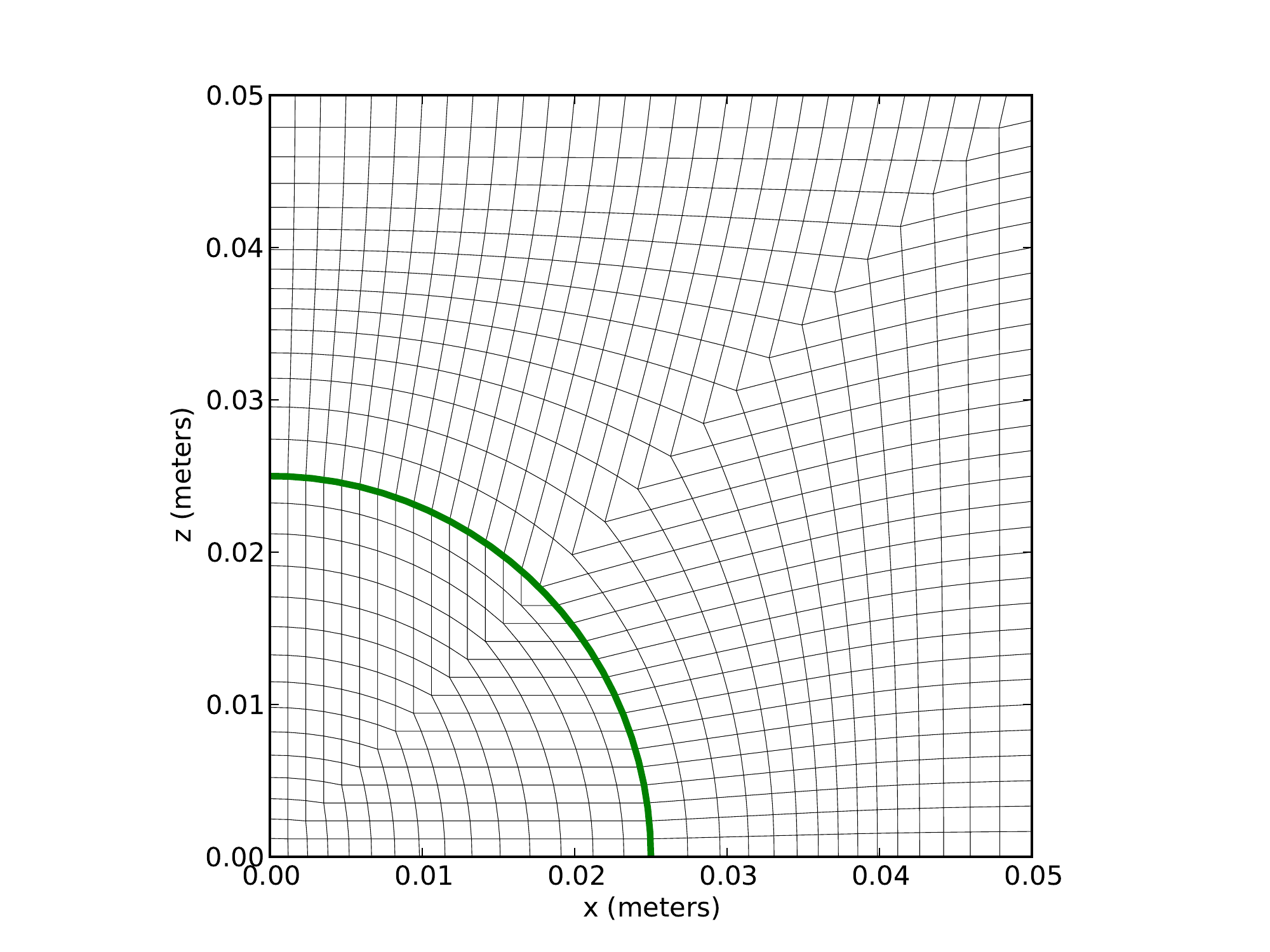}}
    \subfloat[Femur model \label{fig:femur-grid-map}]
             {\includegraphics[width=0.5\textwidth]{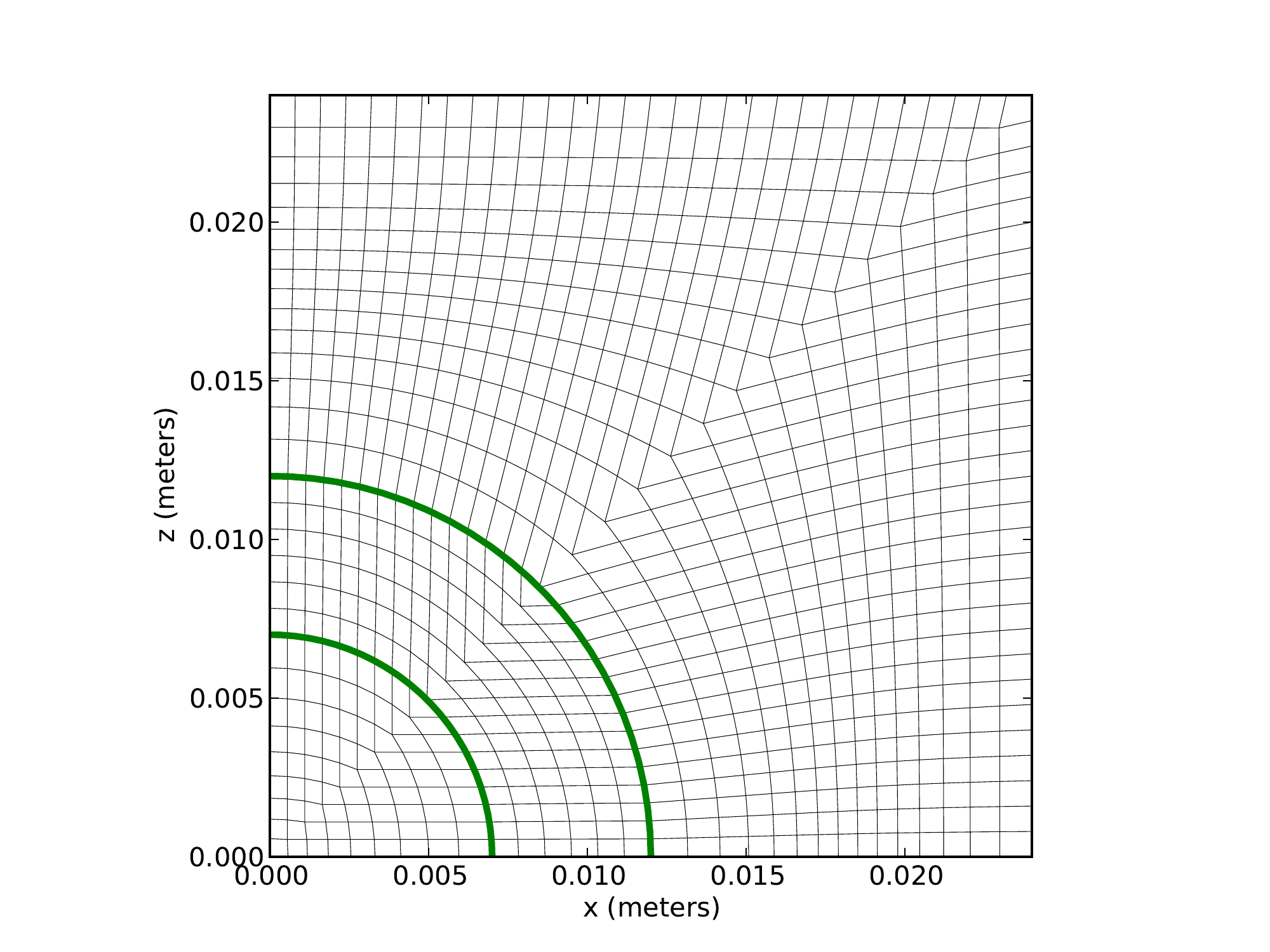}}
    \caption{Grid mappings used for cylindrical scatterer and femur
      models, illustrated on a grid $120 \times 120$ cells in total
      size.  The thick green circles denote material boundaries.  Only
      the upper-right quadrant of the distorted region is shown; the
      grid is rectilinear above and to the right of the region shown,
      while the parts below and to the left are reflections of the
      upper-right quadrant. \label{fig:mapped-grids}}
  \end{center}
\end{figure}

The grid mapping used here deserves some discussion.  It is closely
related to the square-to-circle mappings of Calhoun, Helzel, and
LeVeque~\cite{cal-hel-rjl:circles}, but the function defining the
radius of curvature of the grid lines has been modified to improve
solution quality.  Based on experimentation, having concentric grid
lines near the surface of the scatterer seems to reduce error, as does
having as little grid line curvature as possible and as even a cell
size as possible in the interior of the scatterer.  Using these
considerations as a guide, in the notation of Section 3.2
of~\cite{cal-hel-rjl:circles} our grid mapping is defined by $D(d) =
r_1 d/\sqrt{2}$, $R(d) = r_1 \left( \frac{9}{10} + d^{19} -
\frac{9}{10} d^{20} \right)$.  This mapping has $R'(1) = r_1$ and
$R''(1) = 0$, giving very nearly concentric grid lines near the
scatterer surface, and maintains a grid line radius of curvature of
$\frac{9}{10} r_1$ or greater throughout the interior of the
scatterer; the small variation in radius of curvature also helps keep
the cell size relatively uniform.  Note that even with these
modifications, any grid of this type must contain some highly
distorted cells, so it should be viewed as one of the worst reasonable
cases that might be encountered.  Figure \ref{fig:sc-map-figure} shows
the mapping in action on a coarse grid.

For each case in Table \ref{tab:sc-case-list}, we simulate for 1.25
cycles of the solution, with initial and boundary conditions
set using the true solution as in the previous section.  The
time step is again chosen so that the global maximum CFL number is
0.9.  We continue to measure error using the grid energy 1- and
max-norms, but since the grid is not uniform we use an area-weighted
grid energy 1-norm, computed as
\begin{equation}
  \text{error 1-norm} = \frac{\sum_{i=1}^{N_1} \sum_{j=1}^{N_2}
    \kappa_{ij} \| \bQ_{ij,\text{numerical}} - \bQ_{ij,\text{true}} \|_E}
  {\sum_{i=1}^{N_1} \sum_{j=1}^{N_2} \kappa_{ij}},
\end{equation}
where $N_1$ and $N_2$ are the numbers of grid cells in the
computational $\xi_1$ and $\xi_2$ directions, and $\kappa_{ij}$ is the
cell area ratio discussed in Section \ref{sec:mapped-grid-basics}.
This is meant to mimic the spatial average of the energy norm of the error.
The simulation domain is a square 20~cm on a side, with the cylinder
placed at the center.  Since the grid mapping is simplest if we use
grid dimensions that are a multiple of eight (the ratio of total
domain size to cylinder radius), and since the smaller number of cases
encourages us to devote more computing resources to each one, we
increase the grid sizes used to test convergence to $128 \times 128$,
$256 \times 256$, $512 \times 512$, and $1024 \times 1024$.

\begin{figure}[t]
  \begin{center}
    \subfloat[1-norm error for inviscid cases \label{fig:sc-error-invisc-1}]
             {\includegraphics[width=0.4\textwidth]{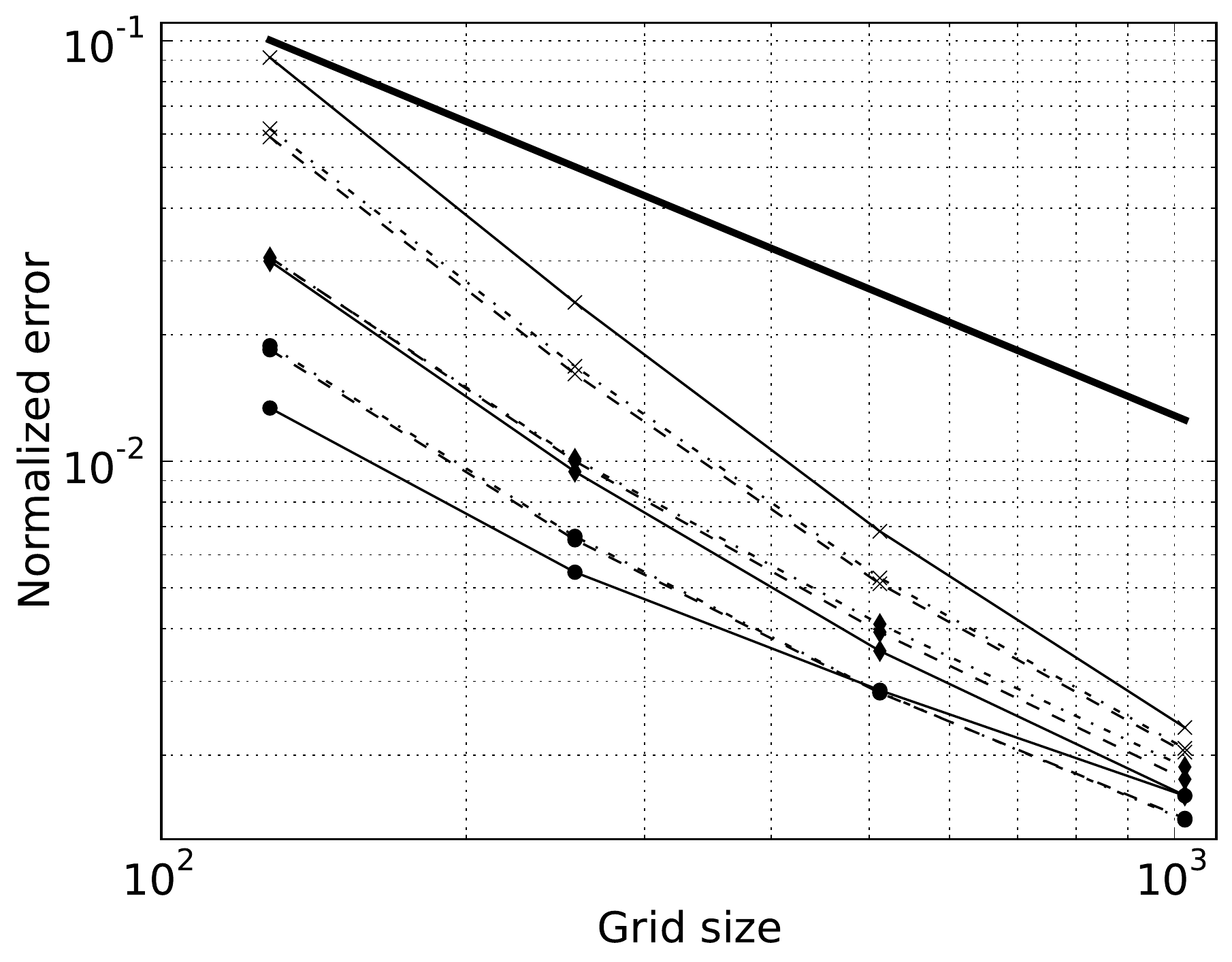}}
    \hspace{0.5in}
    \subfloat[Max-norm error for inviscid cases \label{fig:sc-error-invisc-max}]
             {\includegraphics[width=0.4\textwidth]{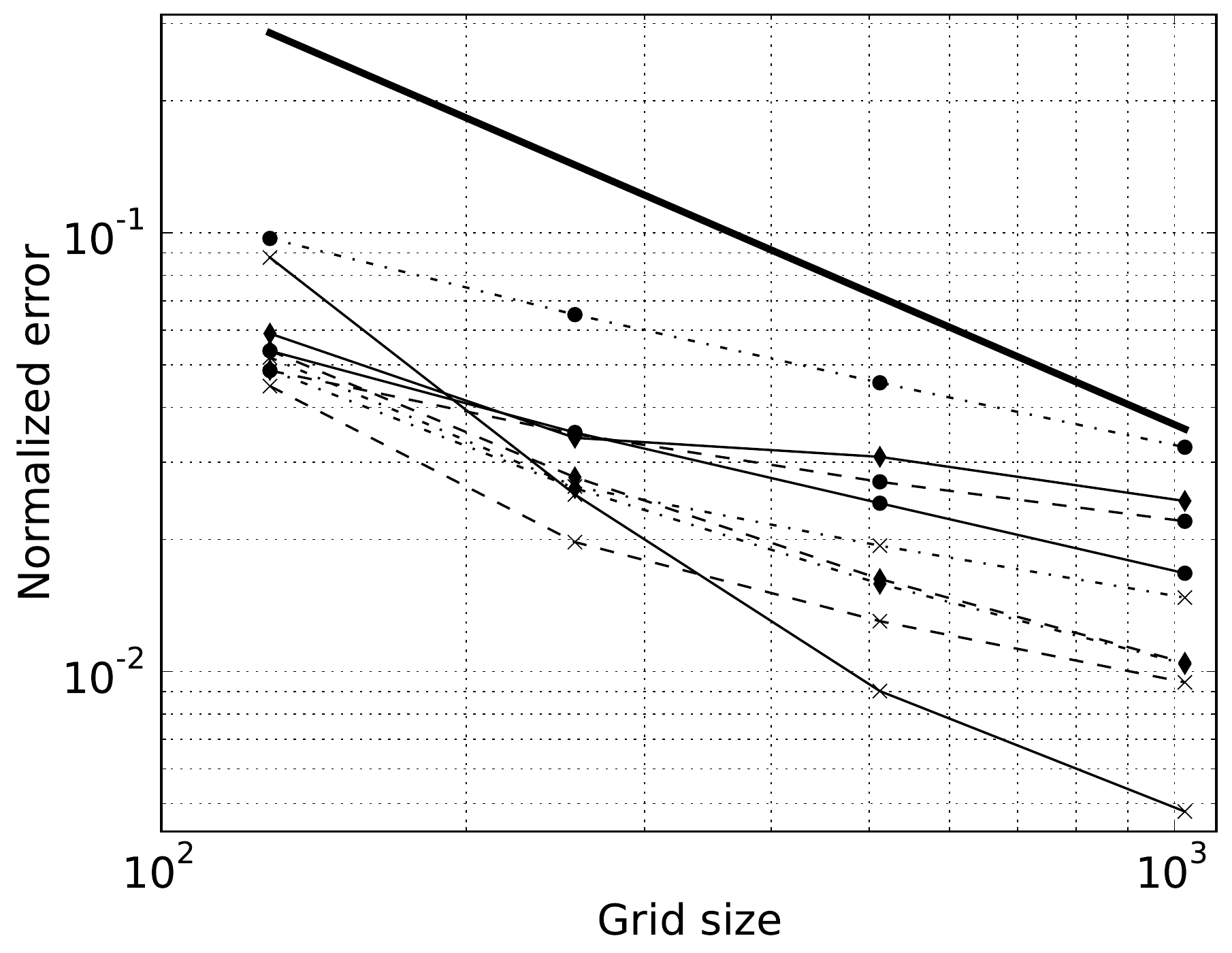}}\\
    \subfloat[1-norm error for viscous cases \label{fig:sc-error-visc-1}]
             {\includegraphics[width=0.4\textwidth]{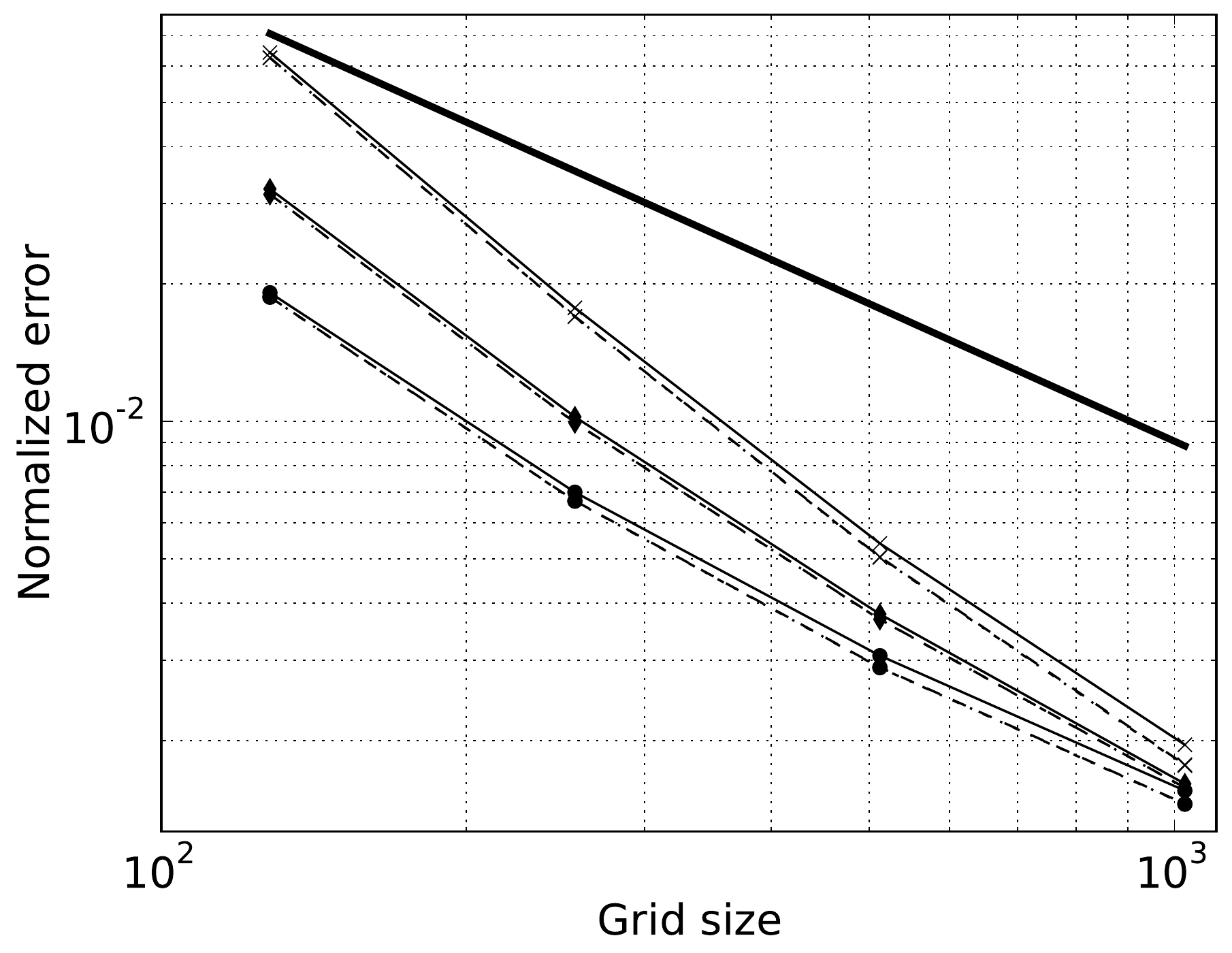}}
    \hspace{0.5in}
    \subfloat[Max-norm error for viscous cases \label{fig:sc-error-visc-max}]
             {\includegraphics[width=0.4\textwidth]{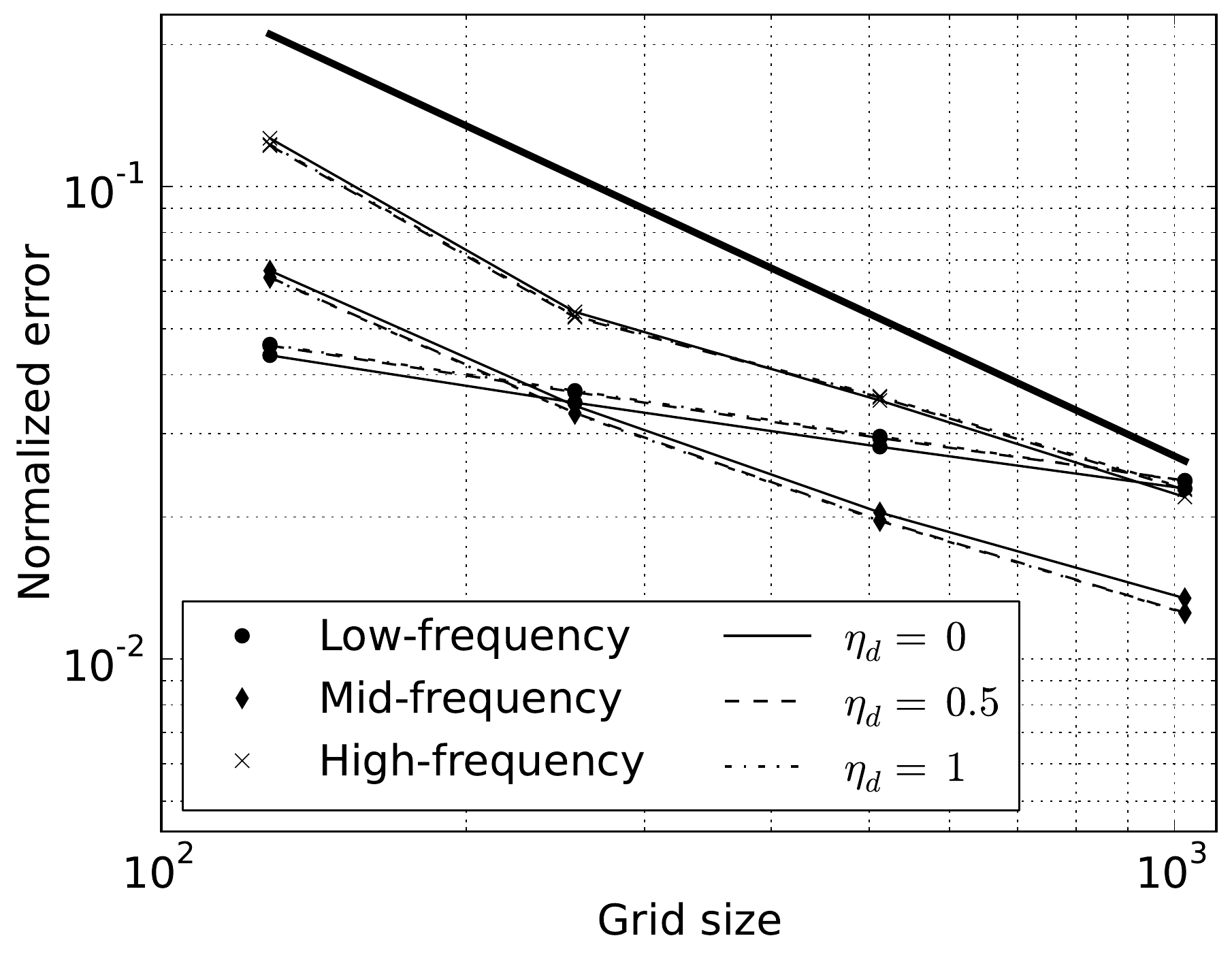}}
    \caption{Convergence behavior of cylindrical scatterer test cases.
      The thick black line is a first-order reference
      line.  Error is normalized by the corresponding norm of the
      true solution on the computational grid. \label{fig:sc-error}}
  \end{center}
\end{figure}

Figure \ref{fig:sc-error} shows the results of these convergence
tests.  Convergence rates are degraded relative to the rectilinear
grid results of Section \ref{sec:rt-fluid-poro}.  We achieve roughly
first-order convergence at large grid sizes in the 1-norm, with faster
convergence on coarser grids; in the max-norm, due to the compounding
of error due to the nonsmooth grid mapping with operator splitting
error and the omission of the second-order correction term at the
surface of the scatterer, we typically achieve below first-order
convergence, in some cases as low as order $\frac{1}{2}$.  Relative
errors on the finest grid are typically a few tenths of a percent in
the 1-norm, and a few percent in the max-norm.  The culprits for these
poor convergence results are the highly-distorted, nearly triangular
grid cells where the $45^\circ$ diagonals intersect the surface of the
scatterer.  These cells have very high error compared to all other
cells in the model; moreover, any logically rectangular grid mapping
of this type must produce such cells when it maps a rectangle on the
computational grid to a circle on the physical grid.  A polar grid
could be used instead, but it would introduce difficulties at the
pole, and would be difficult to generalize to model multiple
scatterers or other additional geometry features.  In future work,
these problems might be ameliorated somewhat by aggressive use of
adaptive mesh refinement --- even though the convergence rate is slow,
the error on these corner cells still decreases as the grid is
refined, and heavy grid refinement could still have a reasonable
computational cost if it is restricted to the vicinity of a few
problem cells.  Despite these problems, however, the bottom half of
Figure \ref{fig:sc-case-portraits} shows that on a fine grid, after
one cycle of the periodic solution, the numerical solution returns
very close to its proper value, without any significant artifacts
visible.  (Note that the numerical solutions shown in the figure were
computed using the MC limiter, as opposed to the convergence results,
which were computed without any limiter.  The limiter was included
because its use would be more typical for solving a problem whose
solution is not already known.)

\subsection{Femur model}

With the results of the previous sections providing reason to be
confident in the accuracy of our code, we now turn to a model of a
biological system.  Specifically, we examine a highly simplified model
of an adult human femur bone.  This model is composed of an inner
cylinder of marrow-filled cancellous bone of radius 7~mm, surrounded
by a concentric cylindrical shell of cortical bone 5~mm thick.  The
outer shell is in turn surrounded by a fluid bath, with both the inner
and outer walls of the shell taken to be impermeable to fluid flow.
The properties used for both types of bone material are given in Table
\ref{tab:matprops}; the fluid bath is taken to be water with bulk
modulus $2.25\,\text{GPa}$ and density $1000\,\text{kg}/\text{m}^3$.
The initial condition is an acoustic pulse traveling in the $+x$
direction toward the bone.  The pulse has a Gaussian profile with
frequency width 100~kHz, and starts with its peak 15~mm away from the
surface.  While the pulse frequency width is well beyond the
low-frequency Biot cutoff frequency, we intend this as a demonstration
example only, as opposed to an accurate quantitative model.  Since our
model is linear, we give the incoming pulse a peak pressure of 1~Pa
for simplicity.  We model this concentric cylindrical geometry using a
modification of the grid mapping used for the cylindrical scatterer,
with the grid lines chosen to be concentric circles between the inner
and outer radius of the cortical bone shell.  Figure
\ref{fig:femur-grid-map} shows the mapping on a coarse grid.  The
actual computational grid is $800 \times 800$ cells.  We use simple
zero-order extrapolation to implement nonreflecting boundary
conditions; the simulation time is too short for significant
interaction between the reflected and transmitted waves from the bone
and the boundaries of the domain in any case.  Since, unlike the other
test cases in this work, we are not attempting to isolate the
convergence behavior of the wave-propagation algorithm here, we used
the monotonized centered limiter for all waves in this computation.

\begin{figure}[t]
  \begin{center}
    \subfloat[Normalized energy density \label{fig:femur-energy}]
             {\includegraphics[width=0.4\textwidth]{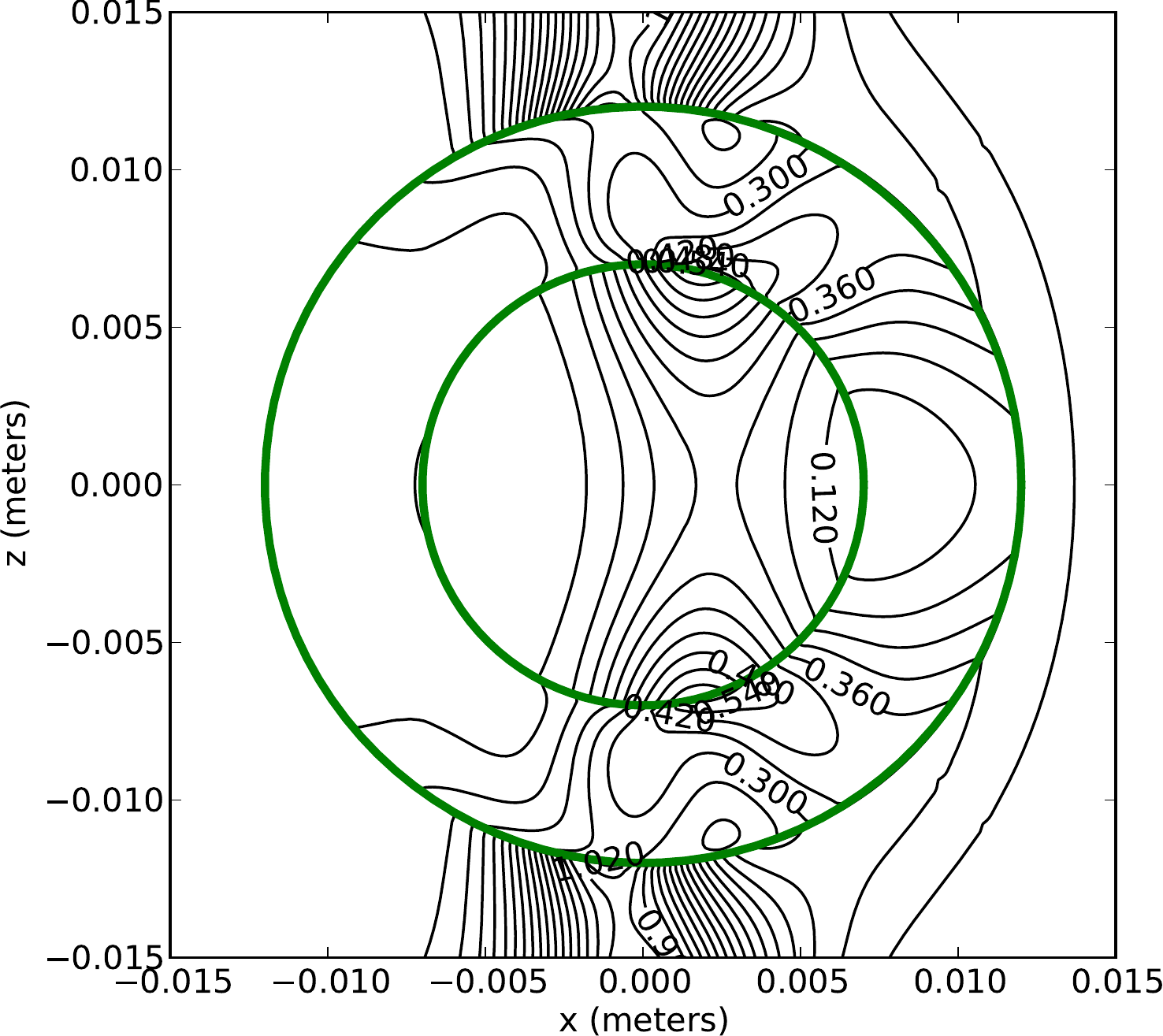}}
    \hspace{0.5in}
    \subfloat[Maximum in-plane shear stress (Pa) \label{fig:femur-shear}]
             {\includegraphics[width=0.4\textwidth]{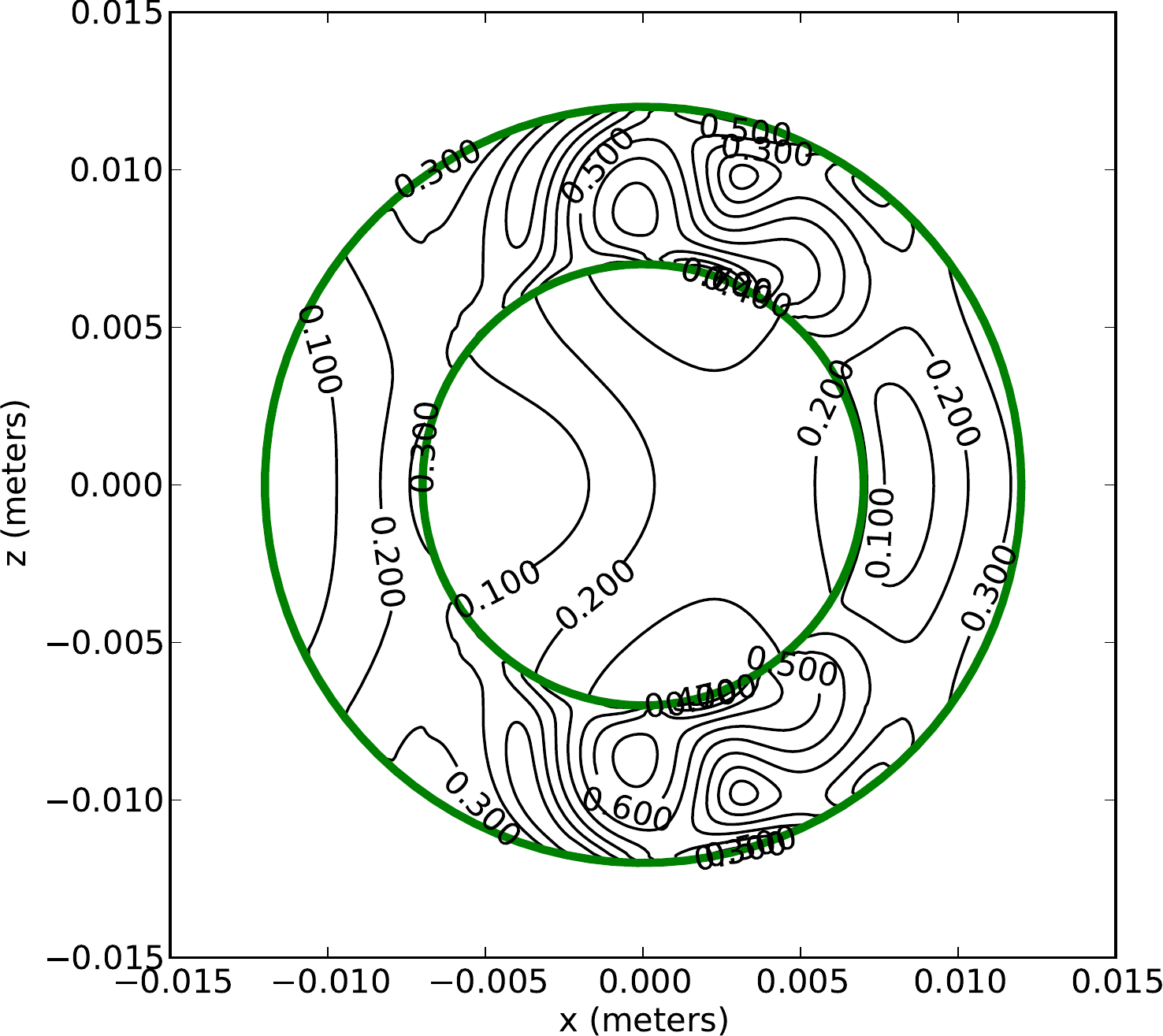}}\\
    \subfloat[Fluid flow speed ($\mu$m/s) \label{fig:femur-flowspeed}]
             {\includegraphics[width=0.4\textwidth]{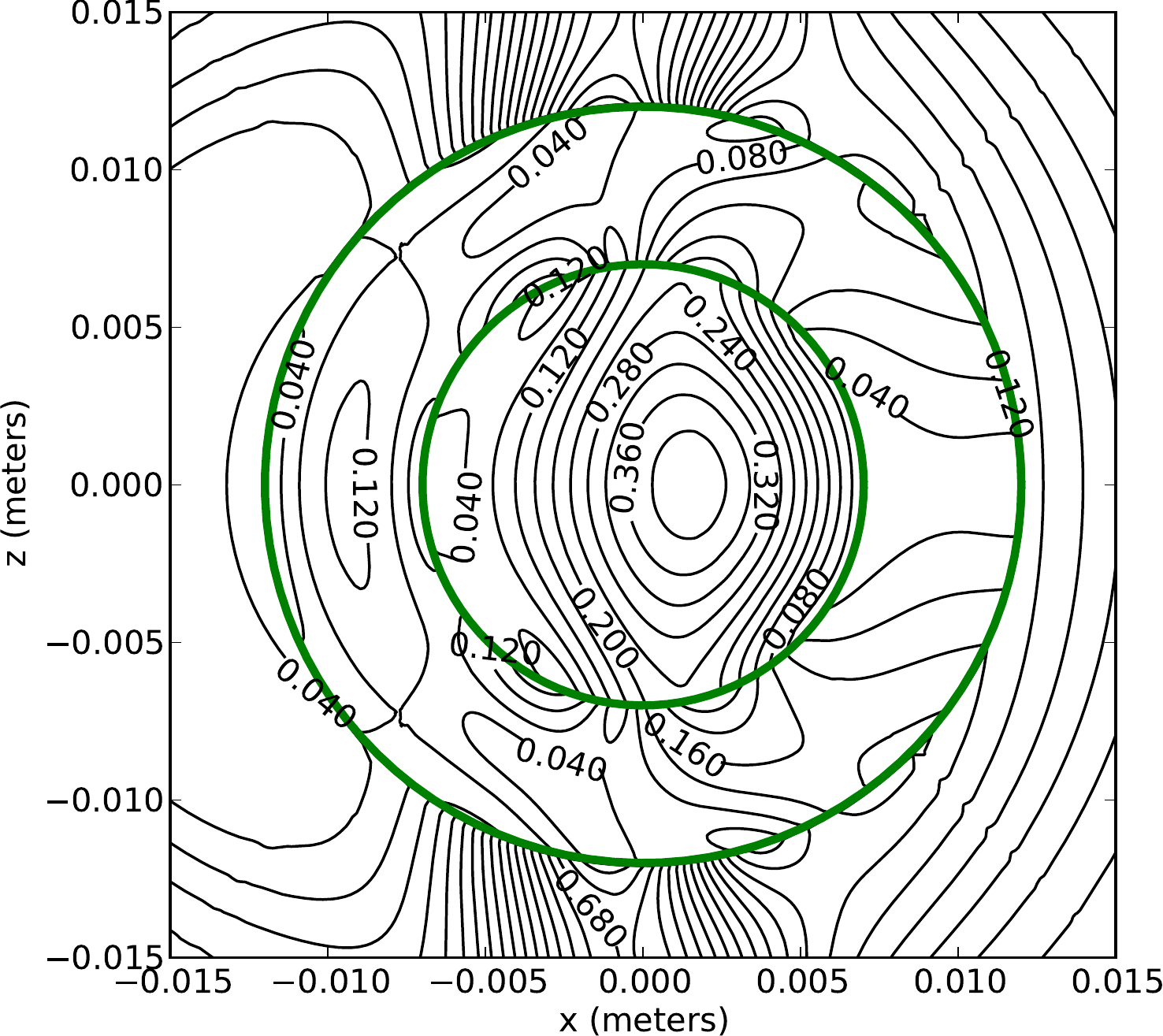}}
    \hspace{0.5in}
    \subfloat[Pressure (Pa) \label{fig:femur-pressure}]
             {\includegraphics[width=0.4\textwidth]{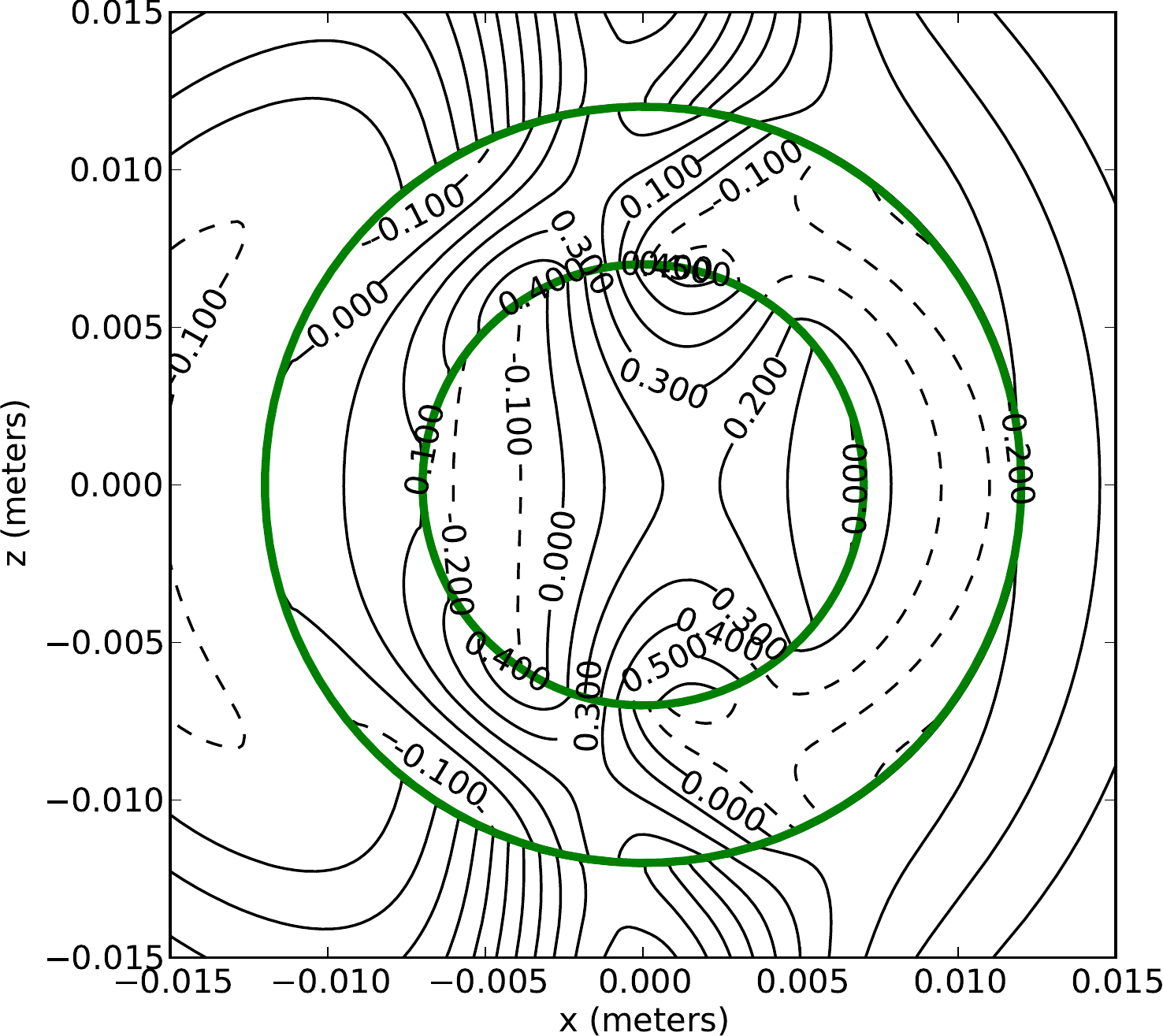}}
    \caption{Contour plots of various solution quantities in the femur
      model, 18~$\mu$s after the start of the simulation.  The thick
      circles mark the boundaries between different
      media. \label{fig:femur-results}}
  \end{center}
\end{figure}

Figure \ref{fig:femur-results} shows a few relevant quantities from
this simulation, at a snapshot 18~$\mu$s after the start.  We plot
energy density normalized by the peak energy density of the incident
pulse, for comparison with the cylindrical scatterer plots, as well
as the maximum in-plane shear stress, as a general measure of the
deformational load on the bone.  We also plot some quantities that may
be biologically relevant: the fluid pressure in the bone, and the
magnitude of fluid flow velocity (the fluid volume flow rate
$\vect{q}$ divided by porosity).  The inclusion of the latter two is
motivated by the review paper of Hillsley and
Frangos~\cite{hillsley-frangos:bone-tissue}, who identify interstitial
fluid flow as relevant to osteogenesis.  While the plots Figure
\ref{fig:femur-results} show some artifacts due to the highly
distorted grid mapping --- namely, the kinks in some of the contours
along the lines of nearly triangular cells extending diagonally from
the center of the circles --- the solution looks qualitatively good
overall.  At the point shown in the solution, the leading fast P
wavefront has passed entirely through the bone.  Meanwhile, separate
fast P and S waves from the initial impact of the pulse, visible in
Figures \ref{fig:femur-energy} and \ref{fig:femur-shear} as the lobes
coming off vertically and to the right from the top and bottom of the
cancellous bone core, are still propagating through the cortical bone.
There are local regions of high fluid velocity (Figure
\ref{fig:femur-flowspeed}), particularly in the cancellous bone but
also associated with the fast P wave propagating around the cortical
bone; there are also local regions of high or low fluid pressure,
associated with a slow P wave to the left of the top and bottom of the
cancellous core in Figure \ref{fig:femur-pressure}, and with the
reflection of the fast P wave off the right exterior surface of the
cortical bone.  The response of the pore fluids is hampered, though,
by the impermeability of the bone surfaces, which prevents strong slow
P waves from being excited by the incident wave.

\section{Summary and future work}
\label{sec:conclusion}

We have developed a high-resolution finite volume method code to model
wave propagation in two dimensions in systems of multiple orthotropic
media and/or fluids.  Our code is capable of using mapped grids to
model curved interfaces, and can also model several standard interface
conditions used for hydraulic contact with a poroelastic medium
--- specifically, open pores, sealed pores, or imperfect hydraulic
contact.  In order to correctly model interfaces between poroelastic
materials and fluids, we have introduced a new transverse Riemann
solution scheme that more carefully handles transverse updates to
media goverened by different systems of equations.  We have also found
that in order to avoid having the poroelastic variables contaminate
the solution in the fluid, at these interfaces we must omit the
standard second-order correction term used in high-resolution finite
volume methods.  While it seems possible to generalize the
second-order correction term to these cases, for now we accept the
reduction in accuracy in the interface cells.

In order to verify our code, we examined its convergence for a variety
of analytical solutions.  We first investigated the effect of omitting
the second-order correction term at an interface using simple plane
wave solutions, by omitting the correction along a line through our
computational domain.  In the 1-norm we observed second-order
convergence with error nearly as low as if the second-order term were
present everywhere; in the max-norm, convergence was degraded to
first-order, but the magnitude of the error was substantially lower
than for a wholly first-order method.  Following this, we tested the
ability of our code to model interfaces between media by comparing its
results to analytical solutions for a time-harmonic plane wave being
reflected off and transmitted through a material interface --- either
between a fluid and a poroelastic medium or between two different
poroelastic media.  We achieved the convergence results that were to
be expected from the plane wave tests: second-order convergence in
1-norm and first-order in the max-norm, demonstrating that our code
can model wave reflection and transmission at an interface with
reasonable accuracy.

Next, we evaluated the ability of our code to handle curved interfaces
and mapped grids by applying it to a time-harmonic acoustic plane wave
scattering off a poroelastic cylinder.  Convergence results were
substantially degraded, to not much better than first-order in the
1-norm, and to roughly order $\frac{1}{2}$ in the max-norm.  The
reason for this was the highly-distorted cells where the mapping
function mapped the corners of a square in the computational domain to
the surface of a circle in the physical domain.  This type of problem
is unavoidable when making a logically rectangular grid of this type,
but these error results could be ameliorated somewhat in the future by
aggressive use of adaptive mesh refinement at the problem cells; in
any case, the numerical solution on a fine grid after one cycle
closely matched the true solution.  We also simulated a simplified
model of a human femur bone in a water bath being struck by a pressure
pulse; all three poroelastic waves contribute noticeably to the
response of the bone, and localized regions of high stress, pressure,
and fluid flow rate occur, which may be of interest biologically.

From here, there are several open directions for future work.  The
most obvious is development of a generalized second-order correction
term that can be used at a material interface; an approach based on
the Immersed Interface Method~\cite{li-ito:iim-book,
  zhang:elastic, rjl-zhang:acou} seems effective, but formulating an appropriate
limiting scheme for this case is not trivial.  Another obvious
extension is to three-dimensional models, in order to handle more
interesting and realistic systems.  The move to three dimensions is
straightforward, but as always the increase in computational work from
two dimensions to three is quite large.  In particular, the number of
transverse Riemann solves grows dramatically in three
dimensions~\cite{langseth-rjl:3d-transverse}, which may necessitate
use of dimensional splitting instead.

There are also opportunities for expansion of the frequency range in
which our code is valid, and in which it achieves its optimal order of
accuracy.  While this paper has avoided the stiff regime identified in
our previous work~\refcartpaper{}, it still represents an impediment
to second-order convergence in cases where an otherwise reasonable
grid and CFL number result in a timestep much longer than the
characteristic dissipation time for the material; this might be
alleviated by an improved Riemann solution process that takes into
account the effect of the dissipation term on the evolution of the
waves, possibly using the work of Hittinger and
Roe~\cite{hittinger-roe:asymptotic-riemann} on Riemann solutions of
hyperbolic systems with relaxation terms, or by using a
semidiscrete numerical scheme with an exponential
time-integrator~\cite{hochbruck-ostermann:exponential} that more
explicitly models the effect of the stiff relaxation term.  Finally,
while we model low-frequency Biot theory here, our code could be
extended straightforwardly to higher-frequency poroelasticity models
by the inclusion of additional memory variables to model the
frequency-dependent kernel used for the generalization of Darcy's law
to high frequencies~\cite{lu-hanyga:poro-memory-drag}.

To aid in reproducibility, we provide all of the code used to generate
the results here at\\ \archivelink{}.

\section{Acknowledgements}

This work has benefited greatly from the direct input and advice of
Prof. Randall J. LeVeque of the Department of Applied Mathematics,
University of Washington, as well as from the \clawpack{} simulation
framework.

\bibliographystyle{siam}
\bibliography{poro}

\appendix

\section{Solution of the poroelastic eigenproblem}
\label{sec:riemann-solver}

\subsection{Equipartition of energy}

A useful property that can be derived from the block structure of the
system is that a simple plane wave of the homogeneous hyperbolic
equation \eqref{eq:first-order-homogeneous} carries equal amounts of
kinetic and potential energy.  To see this, note that a simple
traveling plane wave solution has the form $\bQ(x, z, t) = \B{r} f(n_x
x + n_z z - \lambda t)$, where $\B{r}$ and $\lambda$ satisfy the
eigenproblem
\begin{equation} \label{eq:simple-wave-eigen}
  \breve{\bA} \B{r} = \lambda \B{r}, \quad \breve{\bA} := n_x \bA +
  n_z \bB.
\end{equation}
For a traveling wave, we further assume $\lambda$ is nonzero.  Dividing
$\B{r}$ into stress and velocity parts according to the block
partitioning of the system, $\B{r} =: \begin{pmatrix} \B{r}_s^T &
  \B{r}_v^T \end{pmatrix}^T$, this becomes a pair of equations,
\begin{equation} \label{eq:stress-velocity-eigen}
  \breve{\bA}_{sv} \B{r}_v = \lambda \B{r}_s, \quad
  \breve{\bA}_{vs} \B{r}_s = \lambda \B{r}_v,
\end{equation}
or, multiplying through by $\bE$,
\begin{equation} \label{eq:stress-velocity-geneig}
  \bE_s \breve{\bA}_{sv} \B{r}_v = \lambda \bE_s \B{r}_s, \quad
  \bE_v \breve{\bA}_{vs} \B{r}_s = \lambda \bE_v \B{r}_v.
\end{equation}

Now, multiply the first equation of \eqref{eq:stress-velocity-geneig}
from the left by $\B{r}_s^T$, giving
\begin{equation}
  \B{r}_s^T \bE_s \breve{\bA}_{sv} \B{r}_v = \lambda \B{r}_s^T \bE_s
  \B{r}_s,
\end{equation}
and notice the symmetry relation between blocks (due to $\bE \bA$ and
$\bE \bB$ being symmetric as discussed in~\refcartpaper{}, and
remaining symmetric when the state variable are permuted to the order
used in this work),
\begin{equation} \label{eq:block-symmetry-relations}
  \bE_s \bA_{sv} = (\bE_v \bA_{vs})^T, \quad
  \bE_s \bB_{sv} = (\bE_v \bB_{vs})^T.
\end{equation}
Using this symmetry relation and the eigenproblem
\eqref{eq:stress-velocity-geneig}, $\B{r}_s^T \bE_s \breve{\bA}_{sv} =
\B{r}_s^T (\bE_v \breve{\bA}_{vs})^T = (\bE_v \breve{\bA}_{vs}
\B{r}_s)^T = (\lambda \bE_v \B{r}_v)^T$.  Since we took $\lambda \ne
0$ for a traveling wave, we may divide it out; remembering that
$\bE_v$ is symmetric, we obtain
\begin{equation} \label{eq:energy-equipartition}
  \B{r}_v^T \bE_v \B{r}_v = \B{r}_s^T \bE_s \B{r}_s.
\end{equation}
Note that there is no need to distinguish between the ordinary
transpose and the Hermitian here, since as demonstrated
in~\refcartpaper{} the eigenvectors are real-valued.  Note also that
there was no specific reference to the physics of poroelastic medium
here, only to general algebraic properties of the system --- it is
easy to verify that the matrices for linear acoustics satisfy the same
block symmetry relations, so the same $\bE$-orthogonality and energy
equipartition results hold for acoustics as well.

\subsection{Solution of the eigenproblem}

In order to solve the Riemann problem, we first need the eigenvalues
and eigenvectors of the $\breve{\bA}$ matrix on either side of the
interface.
For our poroelastic material of interest,
let $(n_1, n_3)$ be the components of the unit normal to the Riemann
solve interface, measured in the material principal coordinates.  Then
if we write the first-order system \eqref{eq:first-order-system} in
principal coordinates, $\breve{\bA}$ is
\begin{equation}
  \breve{\bA} = n_1 \bA + n_3 \bB,
\end{equation}
where $\bA$ and $\bB$ are exactly as given in
\eqref{eq:system-matrices-block}, \eqref{eq:A-blocks}, and
\eqref{eq:B-blocks}.  Furthermore, from the properties of $\bA$ and
$\bB$, we know that $\bE \breve{\bA}$ is symmetric, and if we define
$\breve{\bA}_{sv} := n_1 \bA_{sv} + n_3 \bB_{sv}$ and $\breve{\bA}_{vs} :=
n_1 \bA_{vs} + n_3 \bB_{vs}$, then $\bE_s \breve{\bA}_{sv} = (\bE_v
\breve{\bA}_{vs})^T$.

In order to see how this helps us, let $\lambda$ and $\B{r}$ be an
eigenpair of $\breve{\bA}$, so that $\breve{\bA} \B{r} = \lambda
\B{r}$.  Since $\bE$ is nonsingular, this is equivalent to the
generalized eigenproblem
\begin{equation} \label{eq:rp-poro-geneig}
  \bE \breve{\bA} \B{r} = \lambda \bE \B{r}.
\end{equation}
Dividing the eigenvector into stress and velocity blocks, $\B{r}
=: \begin{pmatrix} \B{r}_s^T & \B{r}_v^T \end{pmatrix}^T$, and
referring back to section \ref{sec:poro-review}, we can rewrite the
second equation of \eqref{eq:stress-velocity-geneig} as $(\bE_s
\breve{\bA}_{sv})^T \B{r}_s = \breve{\bA}_{sv}^T \bE_s \B{r}_s =
\lambda \bE_v \B{r}_v$.  If we then multiply the first equation on the
left by $\breve{\bA}_{sv}^T$, we then get a $4\times 4$
symmetric-definite eigenproblem for $\lambda^2$ and $\B{r}_v$:
\begin{equation} \label{eq:rp-geneig}
  \breve{\bA}_{sv}^T \bE_s \breve{\bA}_{sv} \B{r}_v = \lambda^2 \bE_v
  \B{r}_v.
\end{equation}
Aside from reducing the dimension of the eigenproblem, this also
illustrates that the eigenvalues of this system come in positive and
negative pairs -- if $\lambda$ is an eigenvalue, then so is
$-\lambda$.

While the reduction from an $8 \times 8$ unsymmetric ordinary
eigenproblem to a $4 \times 4$ symmetric-definite generalized
eigenproblem is already a substantial gain, further improvement is
possible.  First, factorize $\bE_v$ as $\bE_v = \B{L}\B{L}^T$, where
$\B{L}$ is upper-triangular; such an $\B{L}$ can be found
straightforwardly as
\begin{equation} \label{eq:L-matrix}
  \B{L} = \begin{pmatrix}
    \sqrt{\frac{\Delta_1}{m_1}} & 0 & \frac{\rho_f}{\sqrt{m_1}} & 0 \\
    0 & \sqrt{\frac{\Delta_3}{m_3}} & 0 & \frac{\rho_f}{\sqrt{m_3}} \\
    0 & 0 & \sqrt{m_1} & 0 \\
    0 & 0 & 0 & \sqrt{m_3}
  \end{pmatrix}.
\end{equation}
(We could choose many factorizations of $\bE_v$, but the
upper-triangular structure of $\B{L}$ will be useful later.  The
$\Delta$ quantities are defined as $\Delta_i := \rho m_i - \rho_f^2$.)
Making the variable substitution $\B{r}_v = \B{L}^{-T} \B{y}$, we
arrive at the $4 \times 4$ symmetric eigenproblem
\begin{equation} \label{eq:rp-geneig-M4}
  \B{M}_4 \B{y} = \lambda^2 \B{y}, \quad
  \B{M}_4 := \B{L}^{-1} \breve{\bA}_{sv}^T \bE_s \breve{\bA}_{sv} \B{L}^{-T}.
\end{equation}

As an aside, the way the matrix $\B{M}_4$ is defined offers an
additional opportunity to reduce the amount of computation performed
during the Riemann solve.  Substituting $\breve{\bA}_{sv} = n_1
\bA_{sv} + n_3 \bB_{sv}$ into the definition of $\B{M}_4$, we can
reduce it to a quadratic form in $n_1$ and $n_3$ with matrix
coefficients,
\begin{equation}
  \B{M}_4 = n_1^2 \B{M}_{411} + n_1 n_3 \B{M}_{413} + n_3^2 \B{M}_{433},
\end{equation}
where
\begin{equation}
  \begin{aligned}
  \B{M}_{411} &= \B{L}^{-1} \bA_{sv}^T \bE_s \bA_{sv} \B{L}^{-T}\\
  \B{M}_{413} &= \B{L}^{-1} \bA_{sv}^T \bE_s \bB_{sv} \B{L}^{-T} +
  \B{L}^{-1} \bB_{sv}^T \bE_s \bA_{sv} \B{L}^{-T}\\
  \B{M}_{433} &= \B{L}^{-1} \bB_{sv}^T \bE_s \bB_{sv} \B{L}^{-T}.
  \end{aligned}
\end{equation}
Since the matrices $\B{M}_{411}$, $\B{M}_{413}$, and $\B{M}_{433}$ do
not depend on the interface direction but only on the properties of
the medium, they can be computed once for each material and retrieved
when needed to form the full $\B{M}_4$.

Finally, we can reduce the dimension of the eigenproblem even further
by finding one of the eigenvectors explicitly and removing it from the
computation.  It is easy to verify that $\B{r}_{v0} = \begin{pmatrix}
  0 & 0 & -n_3 & n_1 \end{pmatrix}^T$ satisfies $\breve{\bA}_{sv}
\B{r}_{v0} = 0$; physically, $\B{r}_{v0}$ corresponds to fluid flow
parallel to the Riemann problem interface, which produces no
propagating waves in invisicd Biot theory.  This is in fact the only
vector in the null space of $\breve{\bA}_{sv}$, since the other three
eigenvectors correspond to the three propagating waves with nonzero
wave speeds $\lambda$.  Thus $\B{y}_0 := \B{L}^T \B{r}_{v0} / \|
\B{L}^T \B{r}_{v0} \|_2$ is a null vector of $\B{M}_4$.  Now construct
a $4\times 3$ matrix $\B{Y}_3$ whose columns are orthonormal and all
orthogonal to $\B{y}_0$.  (This is where the upper-triangular
structure of $\B{L}$ comes in handy, since it means that the first two
components of $\B{y}_0$ are zero, which makes this matrix easy to
construct.)  The matrix $\B{Y} := \begin{pmatrix} \B{Y}_3 &
  \B{y}_0 \end{pmatrix}$ is then an orthonormal matrix.  Making the
new variable substitution $\B{y} = \B{Y} \begin{pmatrix} \B{u}^T &
  u_0 \end{pmatrix}^T$ in \eqref{eq:rp-geneig-M4} and multiplying from
the left by $\B{Y}^T$, we then get
\begin{equation}
  \B{Y}^T \B{M}_4 \B{Y} \begin{pmatrix} \B{u} \\ u_0 \end{pmatrix}
  = \begin{pmatrix}
    \B{Y}_3^T \B{M}_4 \B{Y}_3 & 0 \\ 0 & 0
  \end{pmatrix}
  \begin{pmatrix} \B{u} \\ u_0 \end{pmatrix}
  = \lambda^2 \begin{pmatrix} \B{u} \\ u_0 \end{pmatrix}.
\end{equation}
Since we are primarily concerned with propagating waves in the Riemann
solution, we can ignore the possibility of $\lambda = 0$, and need
only find eigenvalues and eigenvectors of the $3\times 3$ matrix
$\B{M}_3 := \B{Y}_3^T \B{M}_4 \B{Y}_3$.  We use the classic QR
algorithm with Wilkinson shifts for this, following the implementation
outlined in Golub and van Loan~\cite{gvl:matrix-computations}, which
has the advantages of being fairly straightforward to implement and
computing all the eigenvectors at once.

The QR algorithm also has another beneficial property for our Riemann
solution process.  It naturally produces eigenvectors with unit
2-norm, which is useful in the orthonormalization of the eigenvectors.
Backing out the variable transformations above, we get $\B{r}_v =
\B{L}^{-T} \B{Y}_3 \B{u}$, so $\B{r}_v^T \bE_v \B{r}_v = \B{u}^T
\B{u}$.  If we divide each eigenvector $\B{u}$ by $\sqrt{2}$ after
computing it with the QR algorithm, we get $\B{r}_v^T \bE_v \B{r}_v =
\frac{1}{2}$, and by the energy equipartition result
\eqref{eq:energy-equipartition}, $\B{r}^T \bE \B{r} = \B{r}_s^T \bE_s
\B{r}_s + \B{r}_v^T \bE_v \B{r}_v = 2 \B{r}_v^T \bE_v \B{r}_v = 1$.
This energy orthonormalization is helpful for calculating wave
strengths in the common special case of a Riemann problem between
identical poroelastic materials.

Once we have $\lambda^2$ and $\B{r}_v$, we get $\lambda$ as either the
positive or negative square root of $\lambda^2$; negative $\lambda$
values correspond to left-going waves, while positive $\lambda$
indicate right-going waves.  Since we have excluded the null space of
$\breve{\bA}_{sv}$ when calculating these eigensolutions, no $\lambda$
thus obtained will be zero.  Finally, we can use
\eqref{eq:stress-velocity-eigen} to obtain $\B{r}_s = \breve{\bA}_{sv}
\B{r}_v / \lambda$, giving all the components of the eigenvector
$\B{r}$ measured in the principal axes of the material; we then
transform $\B{r}$ into the global computational axes to solve the
Riemann problem.

\section{Analytic solution procedure for a time-harmonic plane wave train
  striking a flat interface}
\label{sec:reflect-transmit}

To obtain an analytic solution for a sinusoidal plane wave train
striking a planar interface between two poroelastic media, or a
poroelastic medium and a fluid, we first specify the angular frequency
$\omega$ of the solution, the unit vector $\vect{p}$ in the
propagation direction of the incident wave, and the wave family
desired if the incident wave is in a poroelastic medium.  Calculations
will be done with a complex-valued solution for convenience; for actual
use, we take the real part of this solution.

We assume an ansatz for the incident wave of the form
\begin{equation} \label{eq:rt-incoming-wave}
 \bQ_{\text{in}} (\vect{x}, t) = \B{V}_{\text{in}}
 \exp\left(i(k_{\text{in}}\vect{x}\cdot\vect{p} - \omega t)\right),
\end{equation}
Substuting this into the PDE \eqref{eq:first-order-system}, which
models both poroelasticity and acoustics depending on the choice of
coefficient matrices, we obtain the eigenvalue problem
\begin{equation} \label{eq:rt-incoming-eigen}
  -i\omega \B{V}_{\text{in}} + i k_{\text{in}} (p_x \bA_R + p_z \bB_R)
  \B{V}_{\text{in}} = \bD_R \B{V}_{\text{in}}.
\end{equation}
Note that here $\bA_R$ and $\bB_R$ are taken to be in the global
$x$-$z$ axes, rather than the material principal axes.  The subscript
$R$ denotes the side of the interface corresponding to the incident
and reflected waves, as opposed to $T$, which will denote the side
corresponding to the transmitted waves.  We solve this eigenproblem
for $k_{\text{in}}$ and $\B{V}_{\text{in}}$, and select the
appropriate solution for the desired wave family, if applicable.  We
typically scale the eigenvector $\B{V}_{\text{in}}$ to have unit
energy norm.

With the incident wave solution in hand, we turn to the reflected and
transmitted waves.  We expect a reflected or transmitted wave in every
wave family in each medium, although some of these may be evanescent
waves if the angle of incidence is shallow enough.  While evanescent
waves decay rapidly away from the interface, they are important near
it and must be included in the solution if they occur.  We assume that
each reflected and transmitted wave has a similar complex exponential
form to the incoming wave, but with wavevectors whose magnitude and
direction are both unknown.  Specifically, we take the ansatz
\begin{equation} \label{eq:rt-outgoing-ansatz}
  \bQ_{(R,T)j}(\vect{x}, t) = \B{V}_{(R,T)j} \exp \left( i
  (\vect{k}_{(R,T)j} \cdot \vect{x} - \omega t) \right),
\end{equation}
where $j$ indexes the wave family and $(R,T)$ may be either $R$ or $T$
--- this equation holds for both the reflected and the transmitted
waves.  We then take the total solution field $\bQ(\vect{x},t)$ to be
\begin{equation} \label{eq:rt-total-soln}
  \bQ(\vect{x},t) = \begin{cases}
    \bQ_{\text{in}}(\vect{x},t) + \sum_{j=1}^{N_R}
    \bQ_{Rj}(\vect{x},t), \quad & \text{$\vect{x}$ on incident side of
      interface,}\\
    \sum_{j=1}^{N_T} \bQ_{Tj}(\vect{x},t), \quad & \text{$\vect{x}$ on
      outgoing side of interface,}
  \end{cases}
\end{equation}
where $N_R$ and $N_T$ are the numbers of reflected and transmitted
waves.

At first glance, there seem to be too many unknowns to find a unique
solution for each wave, but it is also necessary to satisfy the
appropriate interface condition; denoting the incident and outgoing
sides of the interface by the subscripts ``in'' and ``out,'' this
condition becomes $\bC_{\text{in}} \bQ_{\text{in, total}} =
\bC_{\text{out}} \bQ_{\text{out}}$, where $\bQ_{\text{in, total}}$ and
$\bQ_{\text{out}}$ are the limits of the state vector approaching the
interface from the incoming and outgoing sides, and $\bC_{\text{in}}$
and $\bC_{\text{out}}$ are the interface condition matrices of section
\ref{sec:riemann-interface}.  Writing $\bQ_{\text{in, total}}$ and
$\bQ_{\text{out}}$ in terms of individual waves, we get
\begin{equation} \label{eq:rt-interface-1}
  \B{C}_{\text{in}} \left. \left( \bQ_{\text{in}}(\vect{x},t) +
  \sum_{j=1}^{N_R} \bQ_{Rj}(\vect{x},t) \right) \right|_{\text{at
      interface}} = \B{C}_{\text{out}} \left. \sum_{j=1}^{N_T}
  \bQ_{Tj}(\vect{x},t) \right|_{\text{at interface}}.
\end{equation}

Now, let $\vect{t}$ be the unit tangent vector along the interface;
without loss of generality assume that the interface runs through the
origin of coordinates, and let $\xi$ measure distance along it from
the origin.  Factoring out the $\exp(-i\omega t)$ time dependence, the
interface condition \eqref{eq:rt-interface-1} becomes
\begin{equation} \label{eq:rt-interface-2}
  \B{C}_{\text{in}} \left( \B{V}_{\text{in}} \exp(i k_{\text{in}}
  \vect{p}\cdot \vect{t}\, \xi) + \sum_{j=1}^{N_R} \B{V}_{Rj} \exp(i
  \vect{k}_{Rj}\cdot \vect{t}\, \xi) \right) = \B{C}_{\text{out}}
  \sum_{j=1}^{N_T} \B{V}_{Tj} \exp(i \vect{k}_{Tj}\cdot \vect{t}\, \xi).
\end{equation}
For this condition to hold at all points on the interface, it is
necessary to have $\vect{k}_{(R,T)j}\cdot \vect{t} = k_{\text{in}}
\vect{p} \cdot \vect{t}$ for every reflected and transmitted wave; we
denote this common tangential component of the wavevector as $k_t$,
and note that it is entirely specified by the direction and frequency
of the incident wave.

Expressing \eqref{eq:first-order-system} in coordinates tangential and
normal to the interface, then applying it to the ansatz
\eqref{eq:rt-outgoing-ansatz} and dividing out the exponential time
and space depedence, we obtain for each outgoing wave
\begin{equation}
  -i\omega \B{V}_{(R,T)j} + i k_t \breve{\bA}_{(R,T)} \B{V}_{(R,T)j} + i
  k_{(R,T)nj} \breve{\bB}_{(R,T)} \B{V}_{(R,T)j} = \bD_{(R,T)} \B{V}_{(R,T)j}.
\end{equation}
Here $\breve{\bA}_{(R,T)} = t_x \bA_{(R,T)} + t_z \bB_{(R,T)}$ is the
flux Jacobian in the tangential direction, and $\breve{\bB}_{(R,T)} =
-t_z \bA_{(R,T)} + t_x \bB_{(R,T)}$ is the flux Jacobian in the
direction normal to the interface; $k_{(R,T)nj}$ is the component of
the wavevector in the interface normal direction.  The subscript
$(R,T)$ is included as a reminder that the system matrices will in
general be different on either side of the interface.  Multiplying
from the left by $\bE_{(R,T)}$ to symmetrize and reorganizing, we
obtain a complex symmetric generalized eigenproblem for each
$k_{(R,T)nj}$ and its corresponding $\B{V}_{(R,T)j}$:
\begin{equation} \label{eq:rt-outgoing-geneig}
  \left( \omega \bE_{(R,T)} - k_t \bE_{(R,T)} \breve{\bA}_{(R,T)} - i
  \bE_{(R,T)} \bD_{(R,T)} \right) \B{V}_{(R,T)j} = k_{(R,T)nj} \bE_{(R,T)}
  \breve{\bB}_{(R,T)} \B{V}_{(R,T)j}.
\end{equation}
Because these represent reflected and transmitted waves, we choose the
eigenvalues $k_{(R,T)nj}$ that correspond to waves propagating away
from the interface (for pure real eigenvalues, for which we determine
propagation direction from the sign of the normal energy flux $\bQ^H
\bE \breve{\bB} \bQ$), or that decay away from the interface (for
complex eigenvalues).  As a practical note, $\breve{\bB}$ is always
singular, and for our formulation of acoustics the matrix on the
left-hand side is also singular with some of its null space in common
with $\breve{\bB}$, so this eigenproblem presents numerical
difficulties (see, for example, section 7.7.3 of Golub and van
Loan~\cite{gvl:matrix-computations}).  To obtain an accurate solution,
we found it necessary to make an additional change of variables to
eliminate the null spaces of both matrices, then solve a smaller
generalized eigenproblem with both matrices nonsingular.  For the
actual eigensolution, we used the {\tt eig} command in
Scipy~\cite{scipy}, which in turn calls the
LAPACK~\cite{lapack-users-guide} routine {\tt ZGGEV}.

Combined with $k_t$, these eigensolutions fully define the reflected
and transmitted waves of \eqref{eq:rt-outgoing-ansatz} up to a scalar
factor; because we have not yet determined that factor, we denote the
eigenvectors found from \eqref{eq:rt-outgoing-geneig} by
$\B{v}_{(R,T)j}$, and let $\B{V}_{(R,T)j} = \beta_{(R,T)j}
\B{v}_{(R,T)j}$.  To find these scalar factors, we return to
\eqref{eq:rt-interface-2} and divide out the common factor of $\exp(i
k_t \xi)$ to obtain
\begin{equation} \label{eq:rt-interface-3}
  \B{C}_{\text{in}} \left( \B{V}_{\text{in}} + \sum_{j=1}^{N_R}
  \beta_{Rj} \B{v}_{Rj} \right) = \B{C}_{\text{out}} \sum_{j=1}^{N_T}
  \beta_{Tj} \B{v}_{Tj}.
\end{equation}
Rearranging and casting in matrix form, we obtain a linear system 
for the $\beta$ values,
\begin{equation}
  \begin{pmatrix}
    -\B{C}_{\text{in}} \B{v}_{R1} & \dots & -\B{C}_{\text{in}} \B{v}_{RN_R} &
    \B{C}_{\text{out}} \B{v}_{T1} & \dots & \B{C}_{\text{out}} \B{v}_{TN_T} 
  \end{pmatrix}
  \begin{pmatrix}
    \beta_{R1} \\ \vdots \\ \beta_{RN_R} \\ \beta_{T1} \\ \vdots
    \\ \beta_{TN_T}
  \end{pmatrix}
  = \B{C}_{\text{in}} \B{V}_{\text{in}}.
\end{equation}
Solving this system completes the information necessary to describe
the reflected and transmitted wave fields, and allows us to compute
the full solution field at any point and time.

\section{Gallery of energy density plots for cylindrical scatterer
  cases}
\label{sec:sc-gallery}

This section contains plots of the energy densities at the initial
time for all 18 cases selected for convergence investigation in
Section 5 of the main paper.  The plots here also include locations
and values of maxima of the energy density within the scatterer, which
were omitted from the main paper due to the difficulty of including
them legibly in a small plot.  The plotting is the same as in Figure
\ref{fig:sc-case-portraits} of the main paper: the plotted quantity is
the energy density, normalized by the peak energy density of the
incident wave, and contours are placed at powers of 2, with dashed
contours indicating negative powers and solid indicating positive
powers.  The blue $\times$ symbols indicate the locations of maxima,
with the adjacent labels giving the value at those maxima.  Since the
solutions are symmetric top-to-bottom, only the top or bottom maximum
of a pair is labeled, not both.

\begin{figure}[ht]
  \begin{center}
    \includegraphics[width=\textwidth]{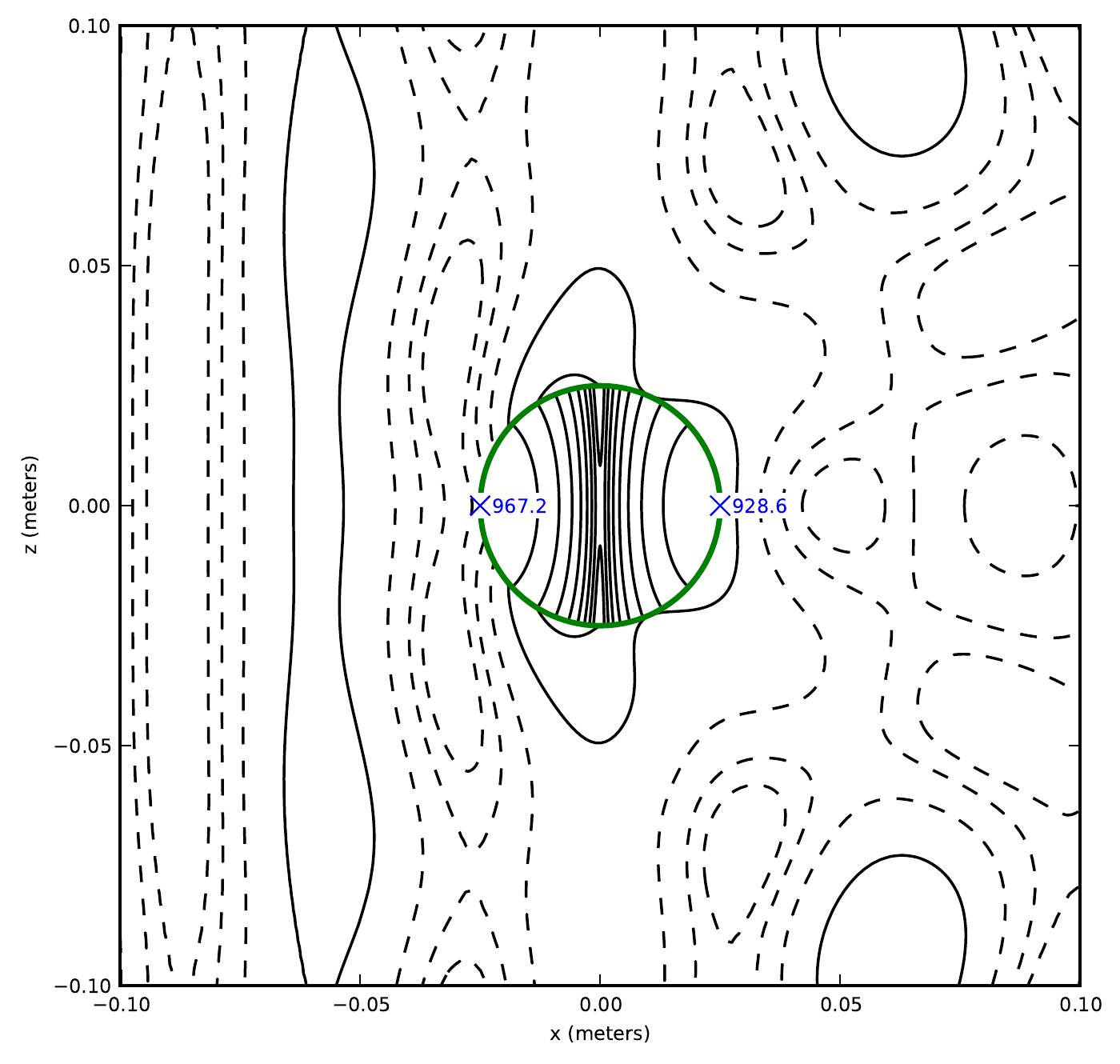}
    \caption{Inviscid, $\eta_d = 0$, 13.25 kHz}
  \end{center}
\end{figure}

\begin{figure}[ht]
  \begin{center}
    \includegraphics[width=\textwidth]{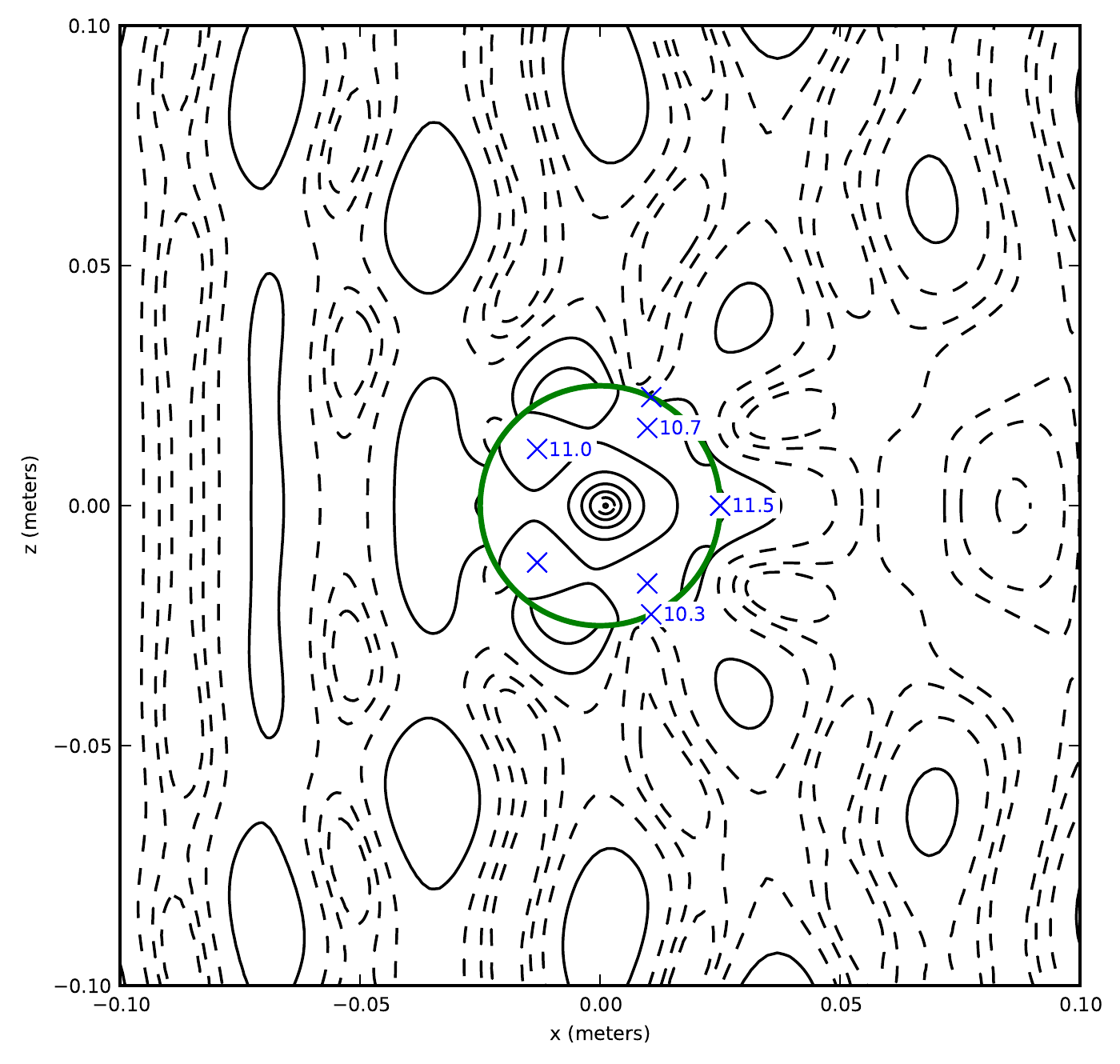}
    \caption{Inviscid, $\eta_d = 0$, 22.25 kHz}
  \end{center}
\end{figure}

\begin{figure}[ht]
  \begin{center}
    \includegraphics[width=\textwidth]{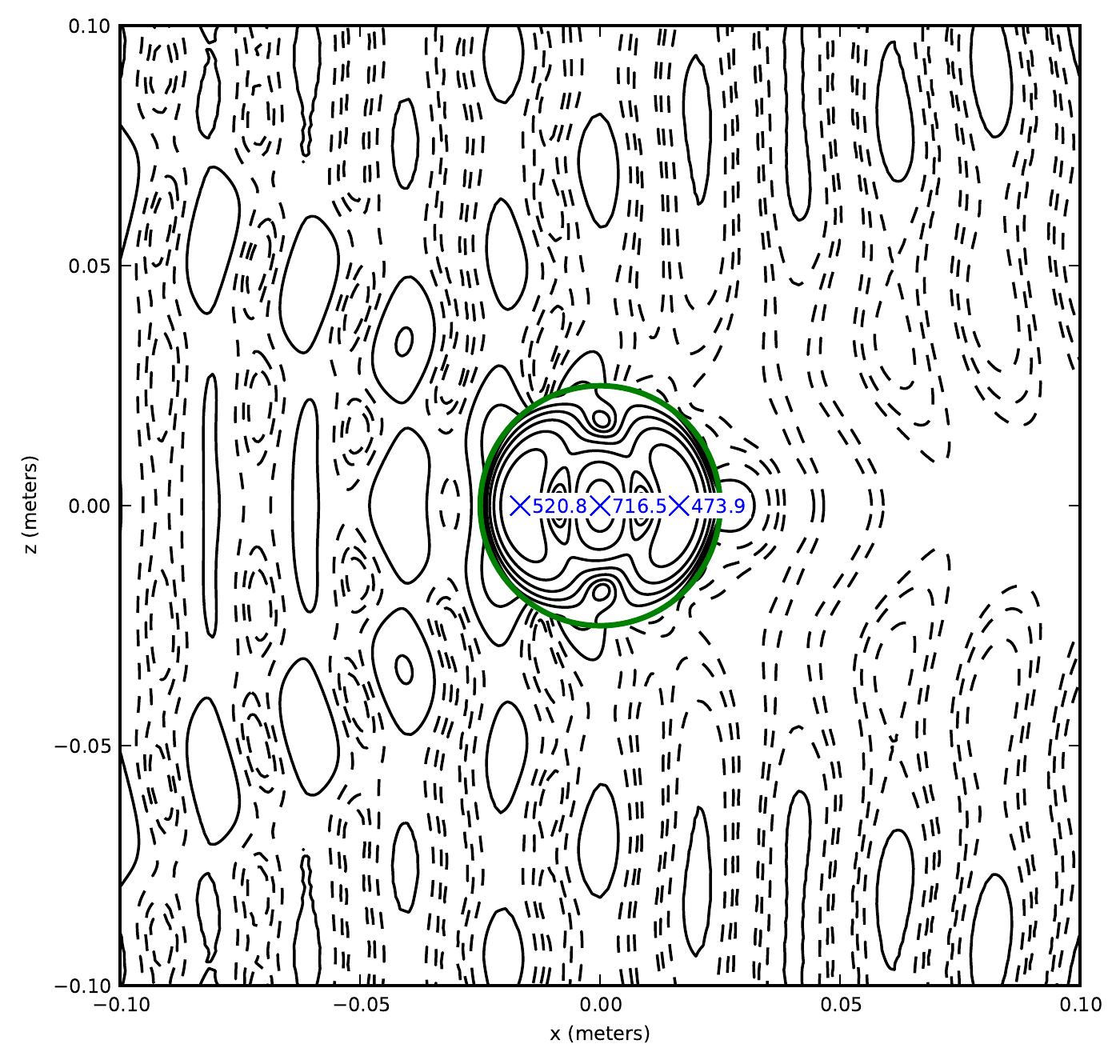}
    \caption{Inviscid, $\eta_d = 0$, 38.20 kHz}
  \end{center}
\end{figure}

\begin{figure}[ht]
  \begin{center}
    \includegraphics[width=\textwidth]{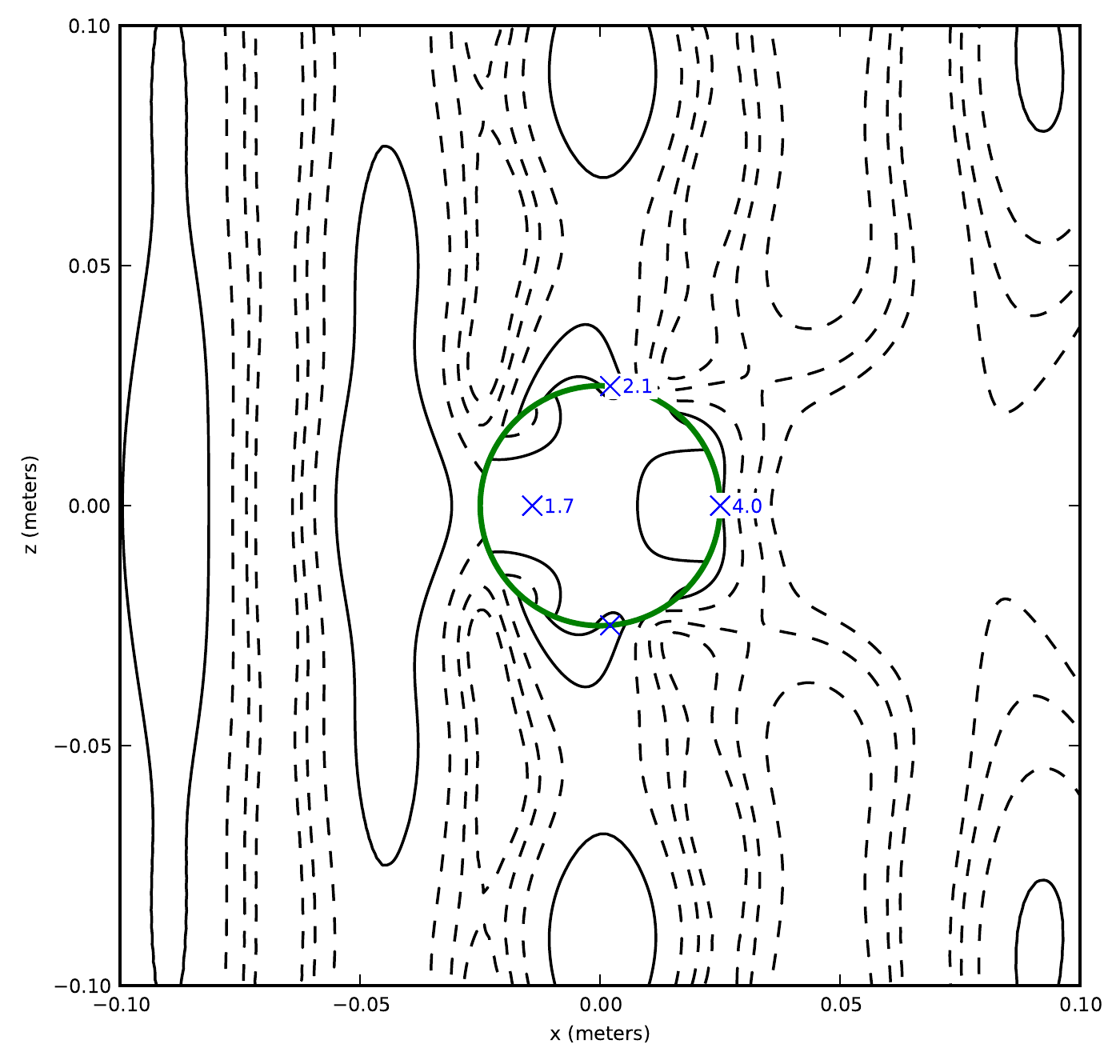}
    \caption{Inviscid, $\eta_d = 0.5$, 17.25 kHz}
  \end{center}
\end{figure}

\begin{figure}[ht]
  \begin{center}
    \includegraphics[width=\textwidth]{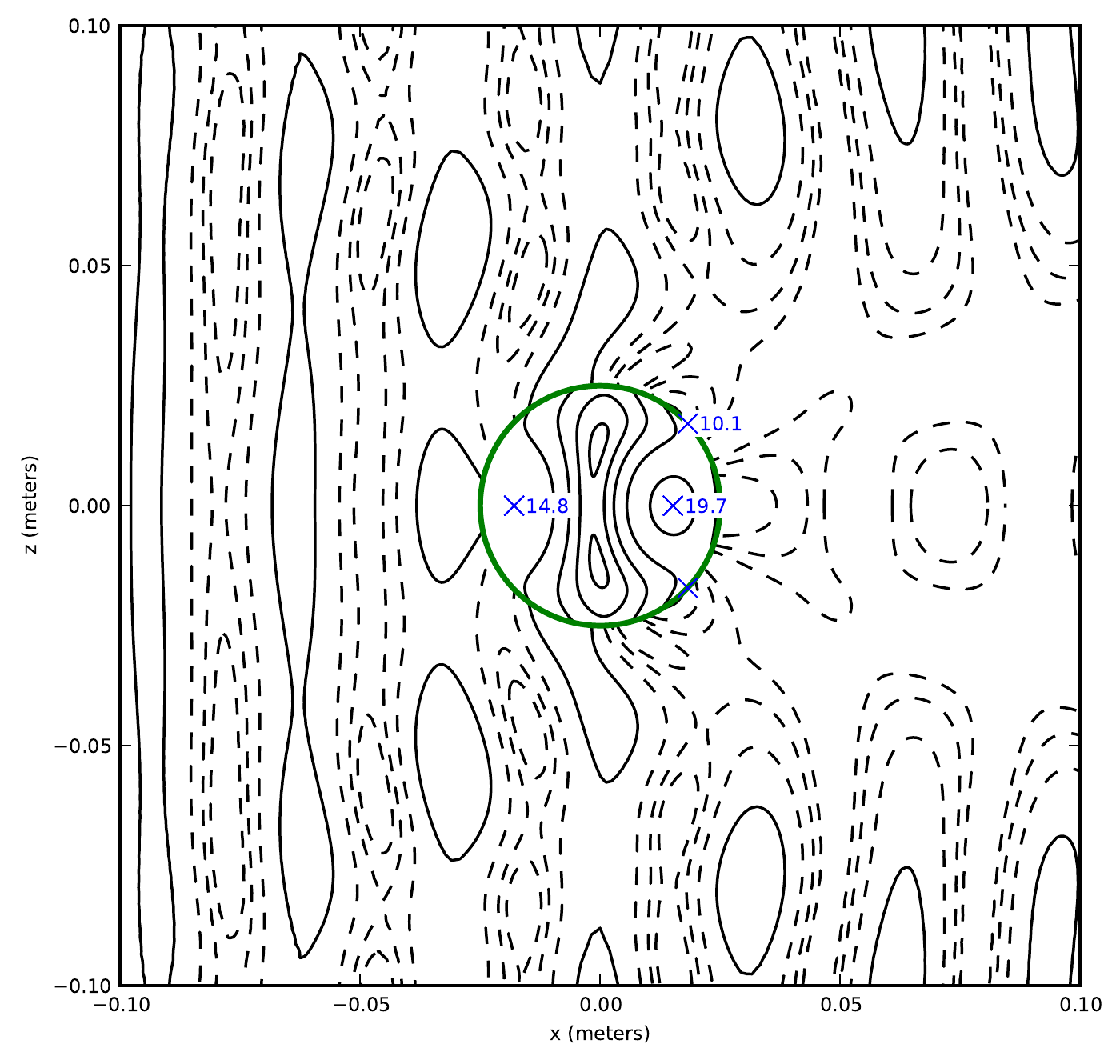}
    \caption{Inviscid, $\eta_d = 0.5$, 25.02 kHz}
  \end{center}
\end{figure}

\begin{figure}[ht]
  \begin{center}
    \includegraphics[width=\textwidth]{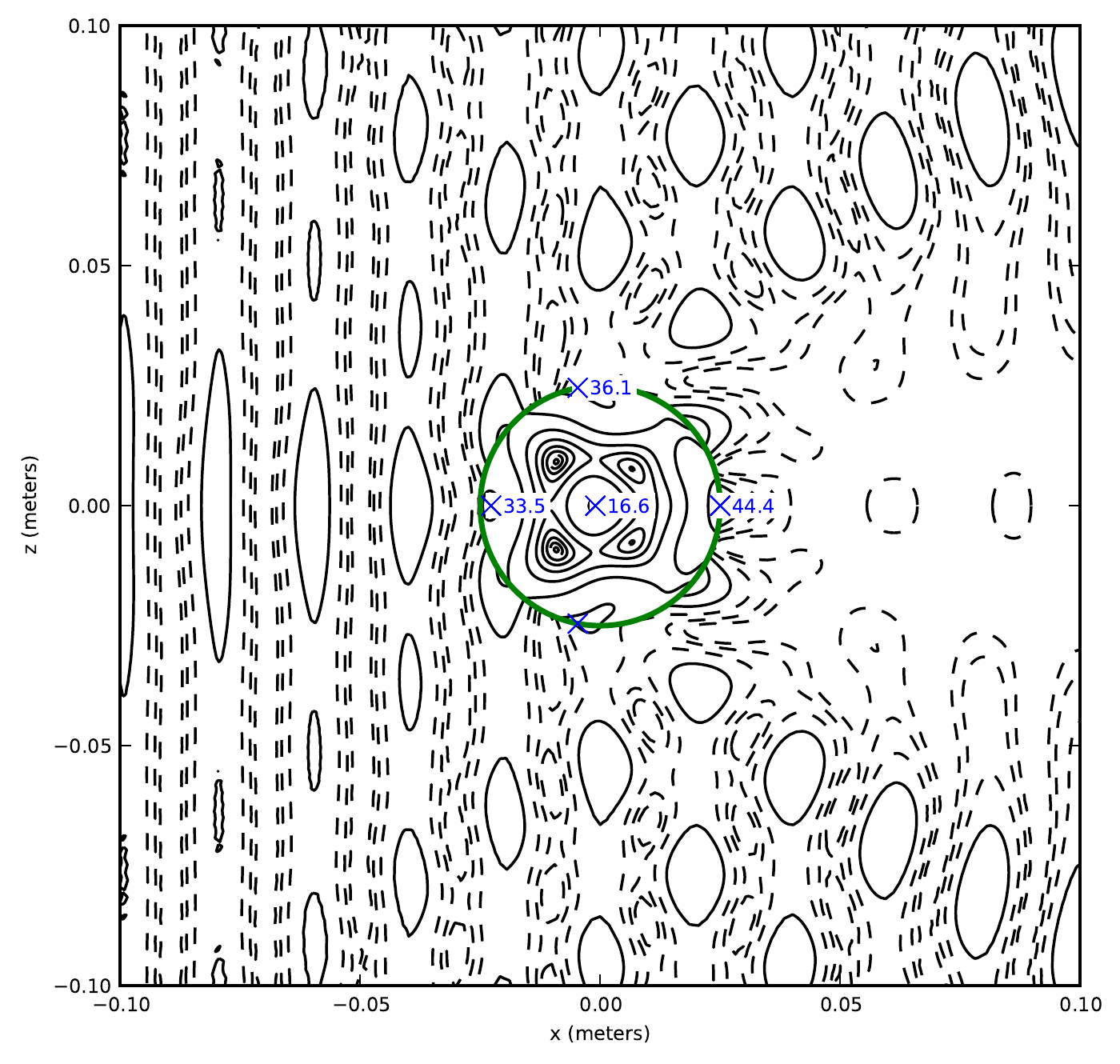}
    \caption{Inviscid, $\eta_d = 0.5$, 39.03 kHz}
  \end{center}
\end{figure}

\begin{figure}[ht]
  \begin{center}
    \includegraphics[width=\textwidth]{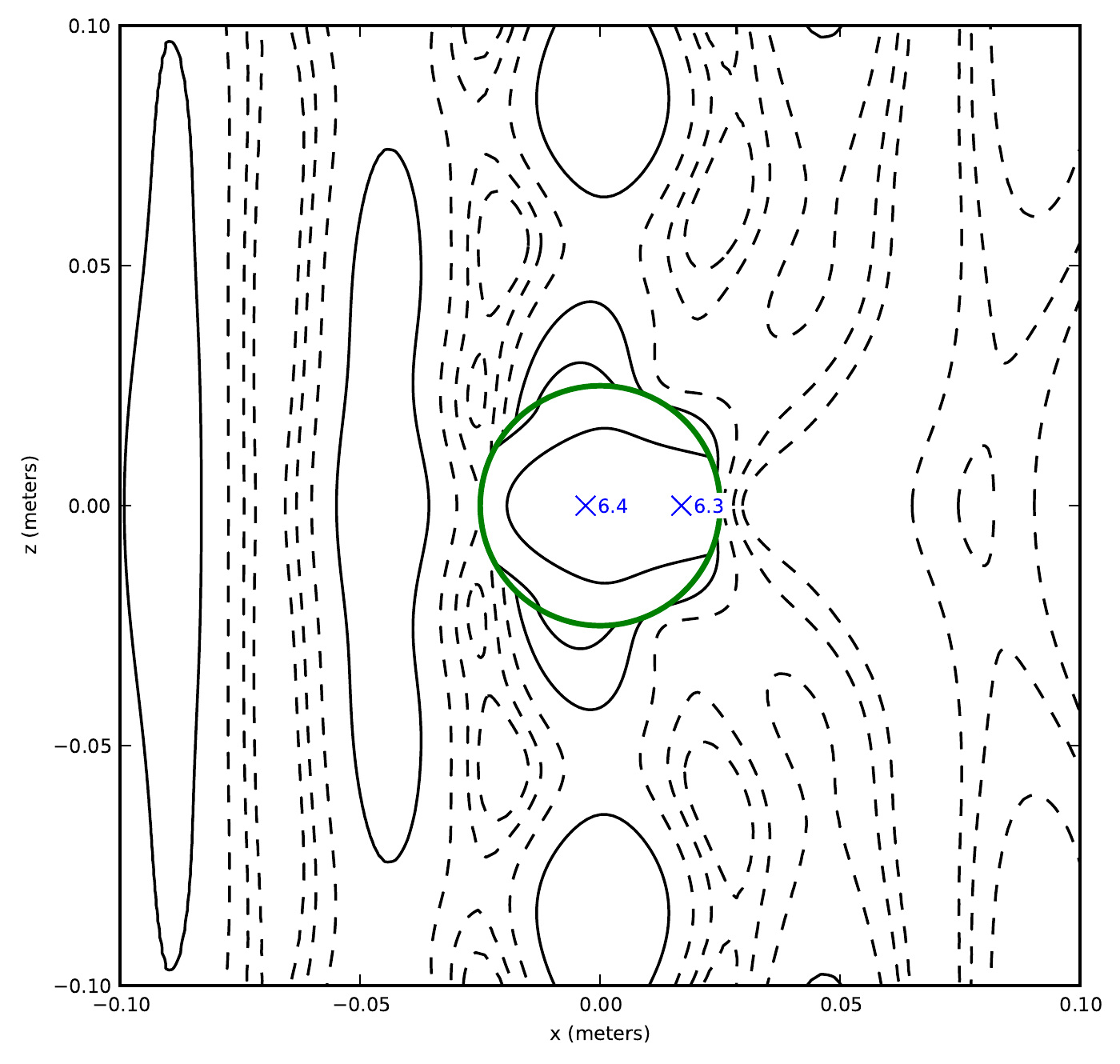}
    \caption{Inviscid, $\eta_d = 1$, 17.30 kHz}
  \end{center}
\end{figure}

\begin{figure}[ht]
  \begin{center}
    \includegraphics[width=\textwidth]{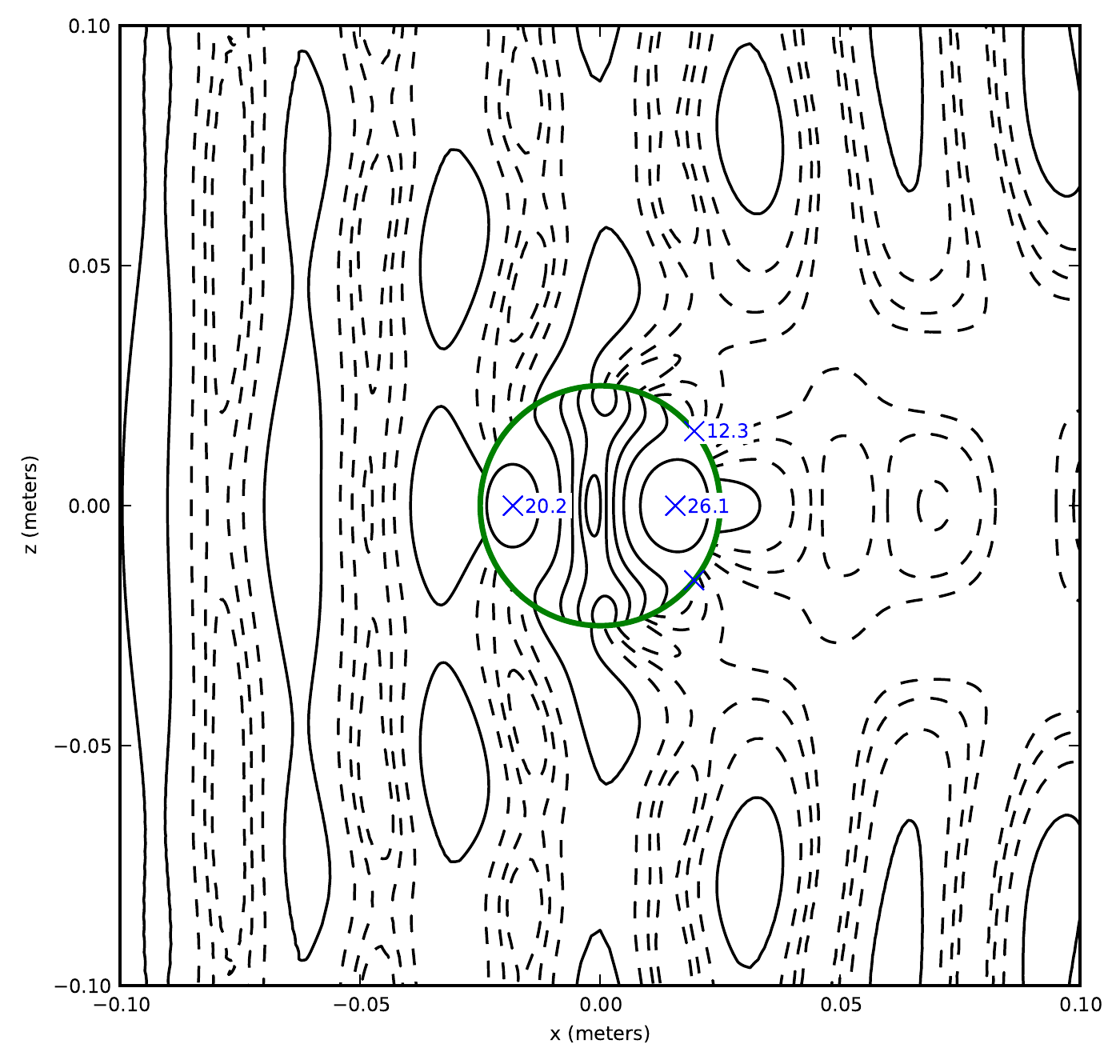}
    \caption{Inviscid, $\eta_d = 1$, 25.09 kHz}
  \end{center}
\end{figure}

\begin{figure}[ht]
  \begin{center}
    \includegraphics[width=\textwidth]{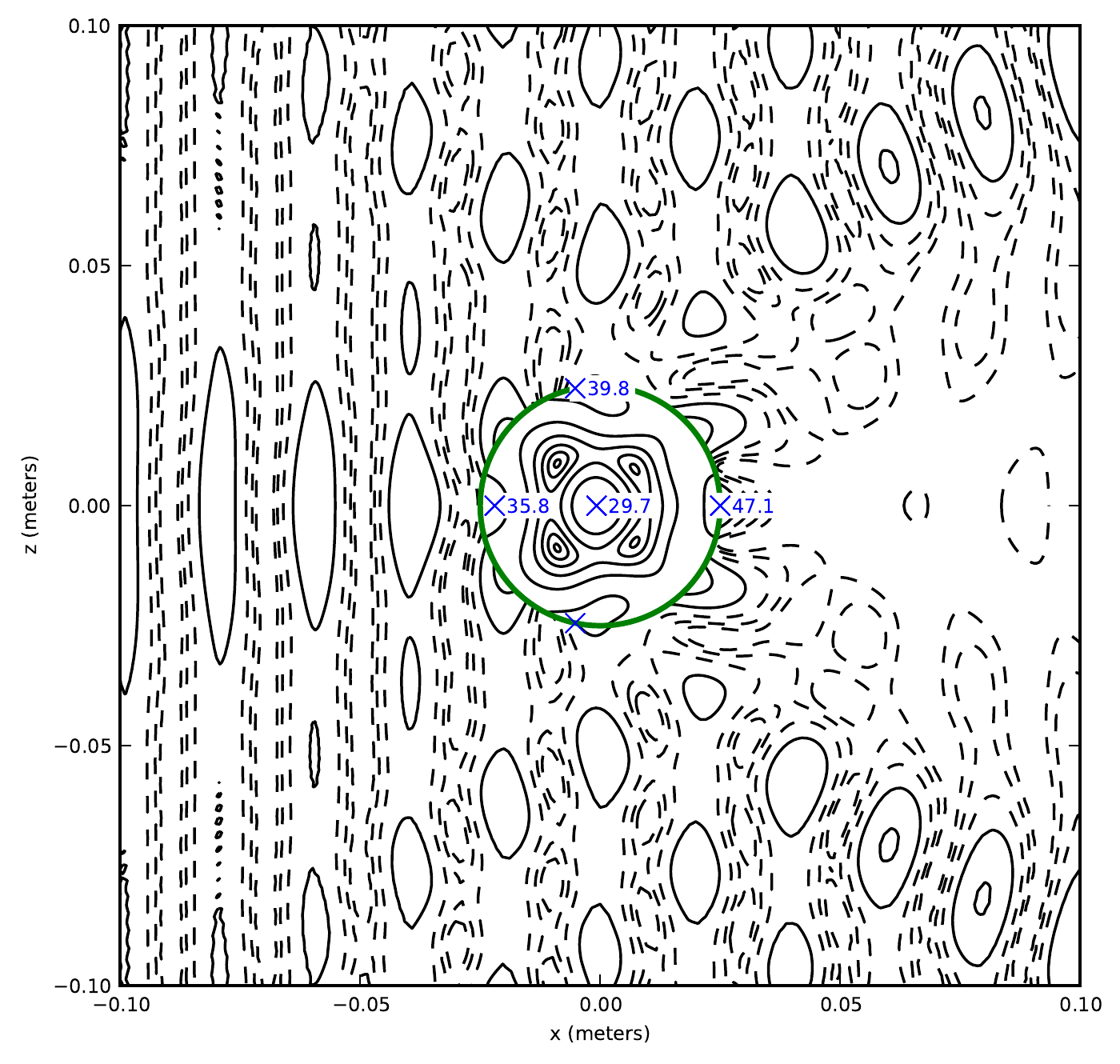}
    \caption{Inviscid, $\eta_d = 1$, 39.04 kHz}
  \end{center}
\end{figure}

\begin{figure}[ht]
  \begin{center}
    \includegraphics[width=\textwidth]{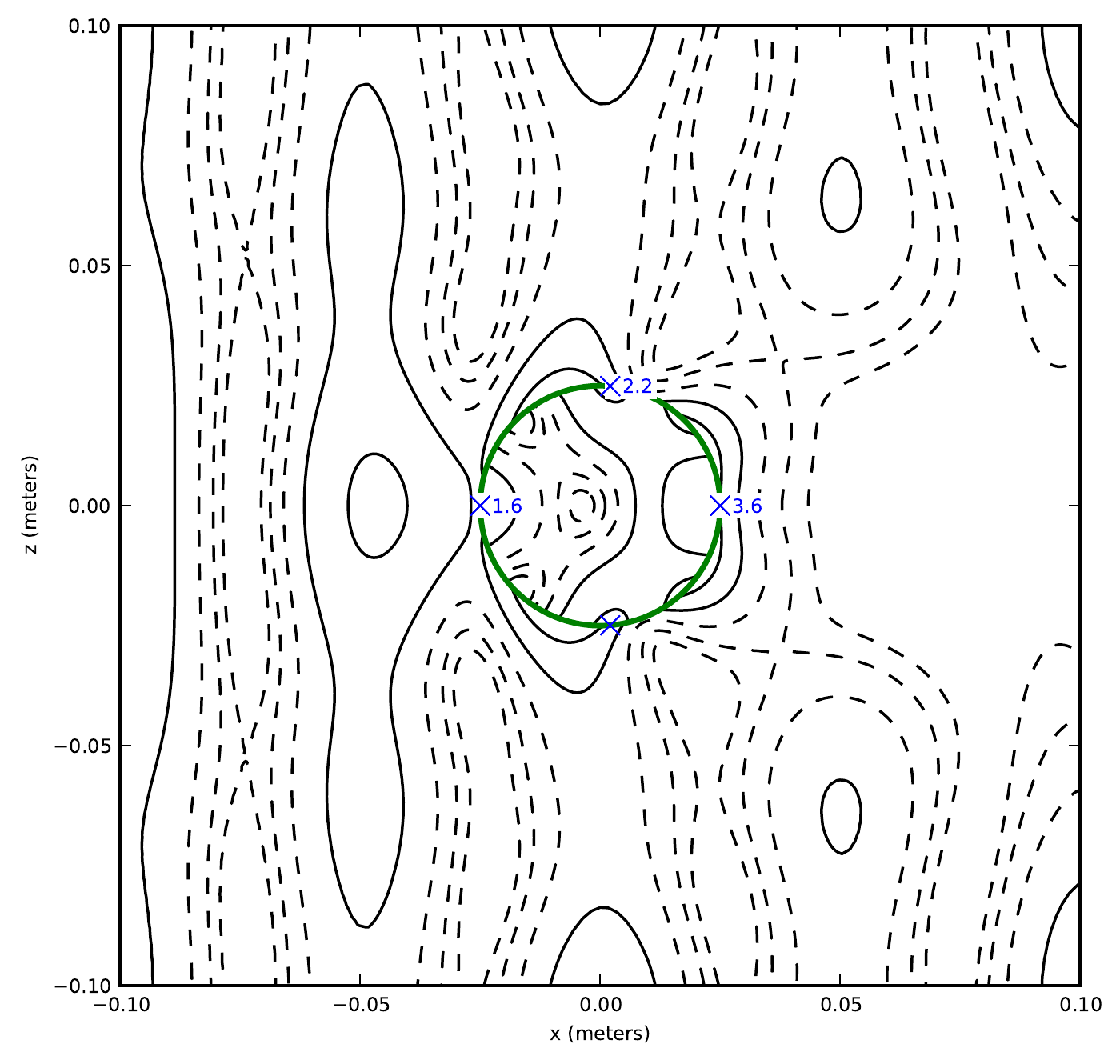}
    \caption{Viscous, $\eta_d = 0$, 15.70 kHz}
  \end{center}
\end{figure}

\begin{figure}[ht]
  \begin{center}
    \includegraphics[width=\textwidth]{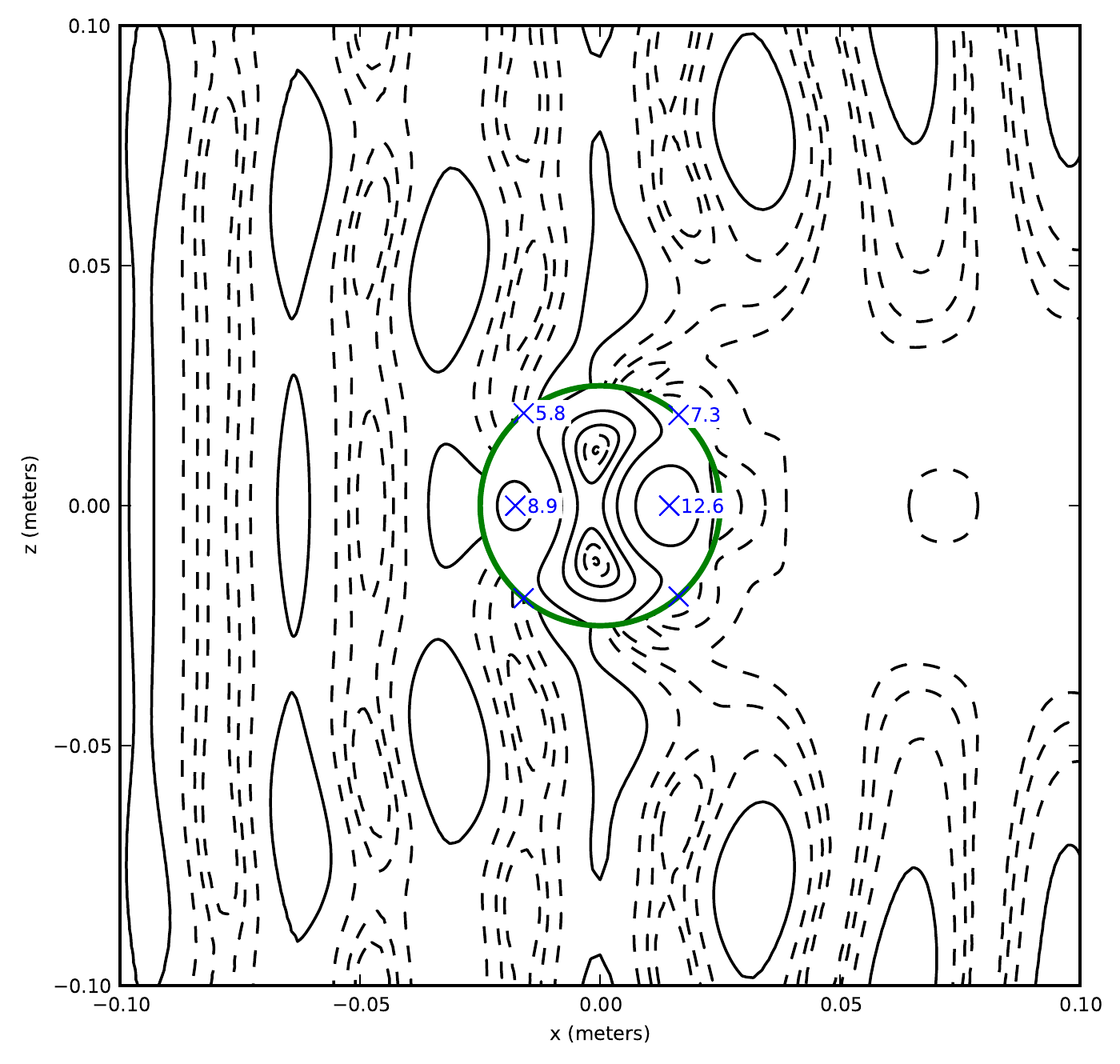}
    \caption{Viscous, $\eta_d = 0$, 24.55 kHz}
  \end{center}
\end{figure}

\begin{figure}[ht]
  \begin{center}
    \includegraphics[width=\textwidth]{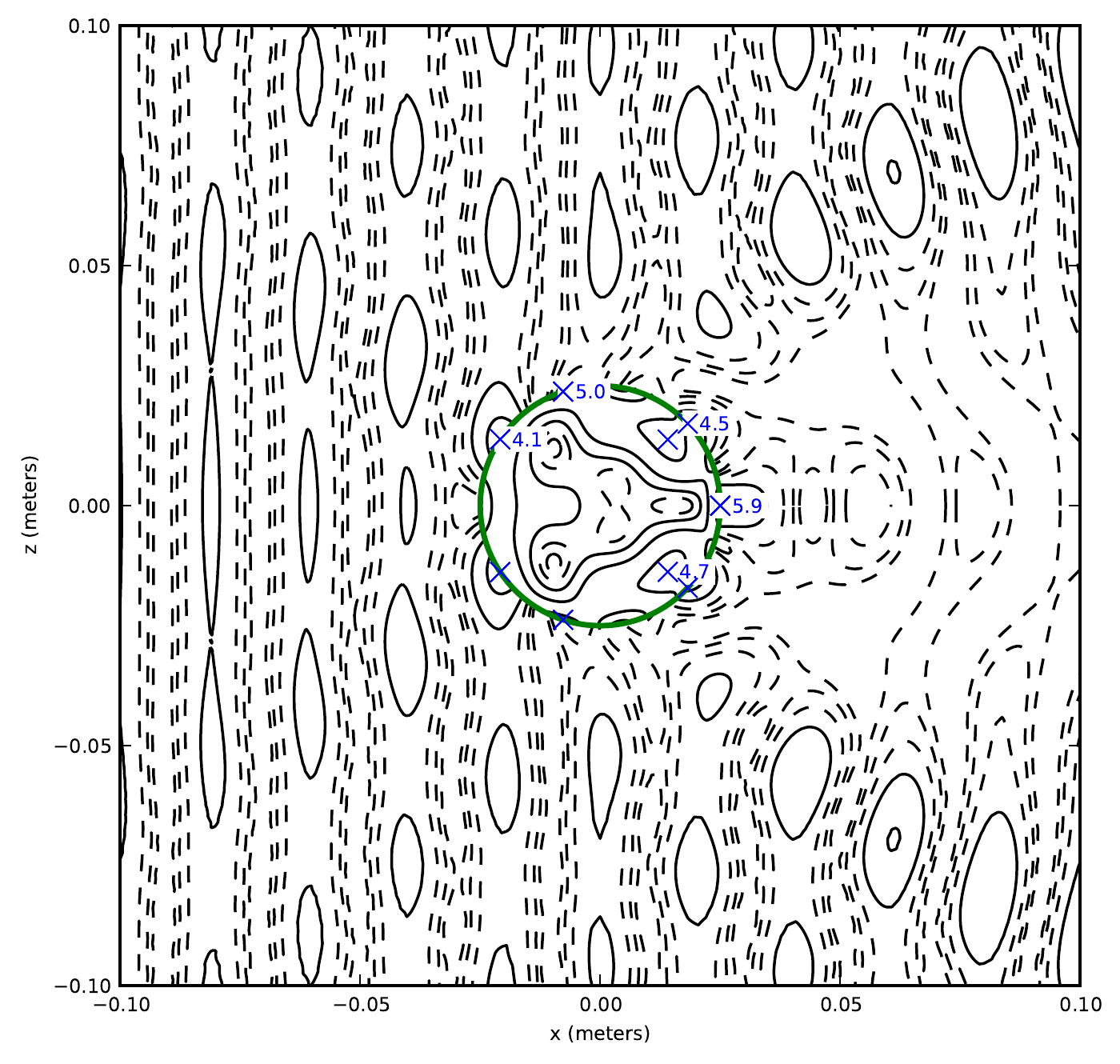}
    \caption{Viscous, $\eta_d = 0$, 38.39 kHz}
  \end{center}
\end{figure}

\begin{figure}[ht]
  \begin{center}
    \includegraphics[width=\textwidth]{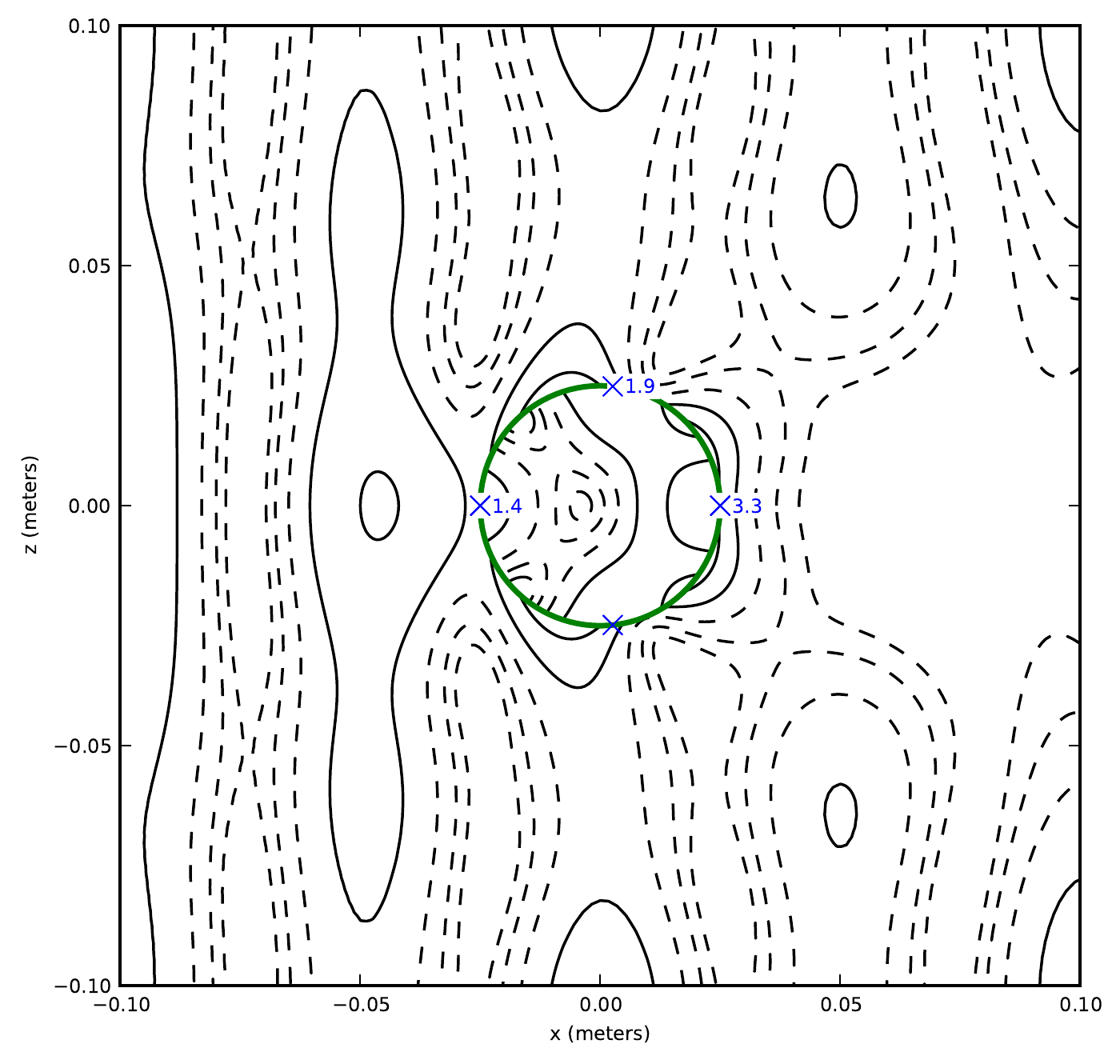}
    \caption{Viscous, $\eta_d = 0.5$, 15.80 kHz}
  \end{center}
\end{figure}

\begin{figure}[ht]
  \begin{center}
    \includegraphics[width=\textwidth]{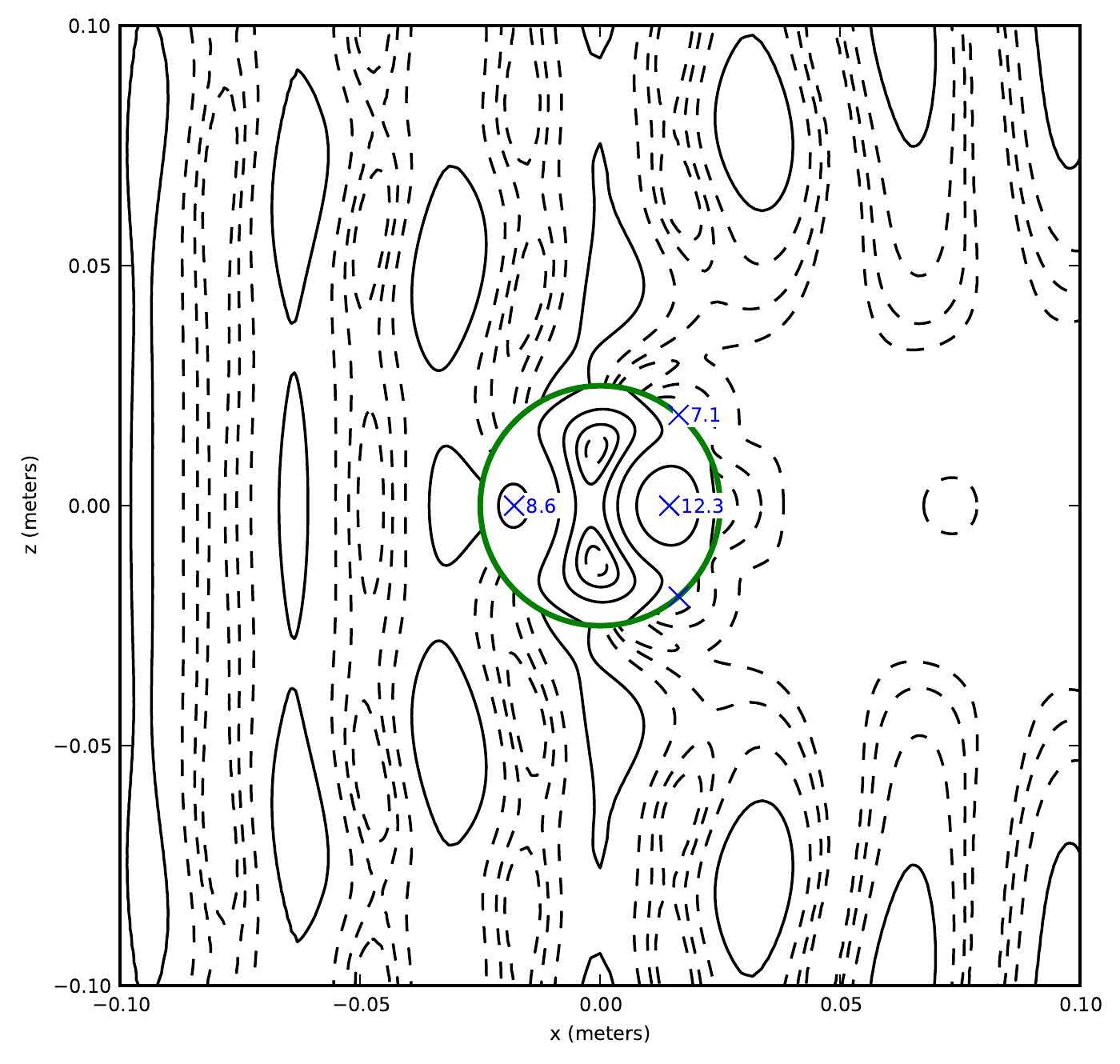}
    \caption{Viscous, $\eta_d = 0.5$, 24.54 kHz}
  \end{center}
\end{figure}

\begin{figure}[ht]
  \begin{center}
    \includegraphics[width=\textwidth]{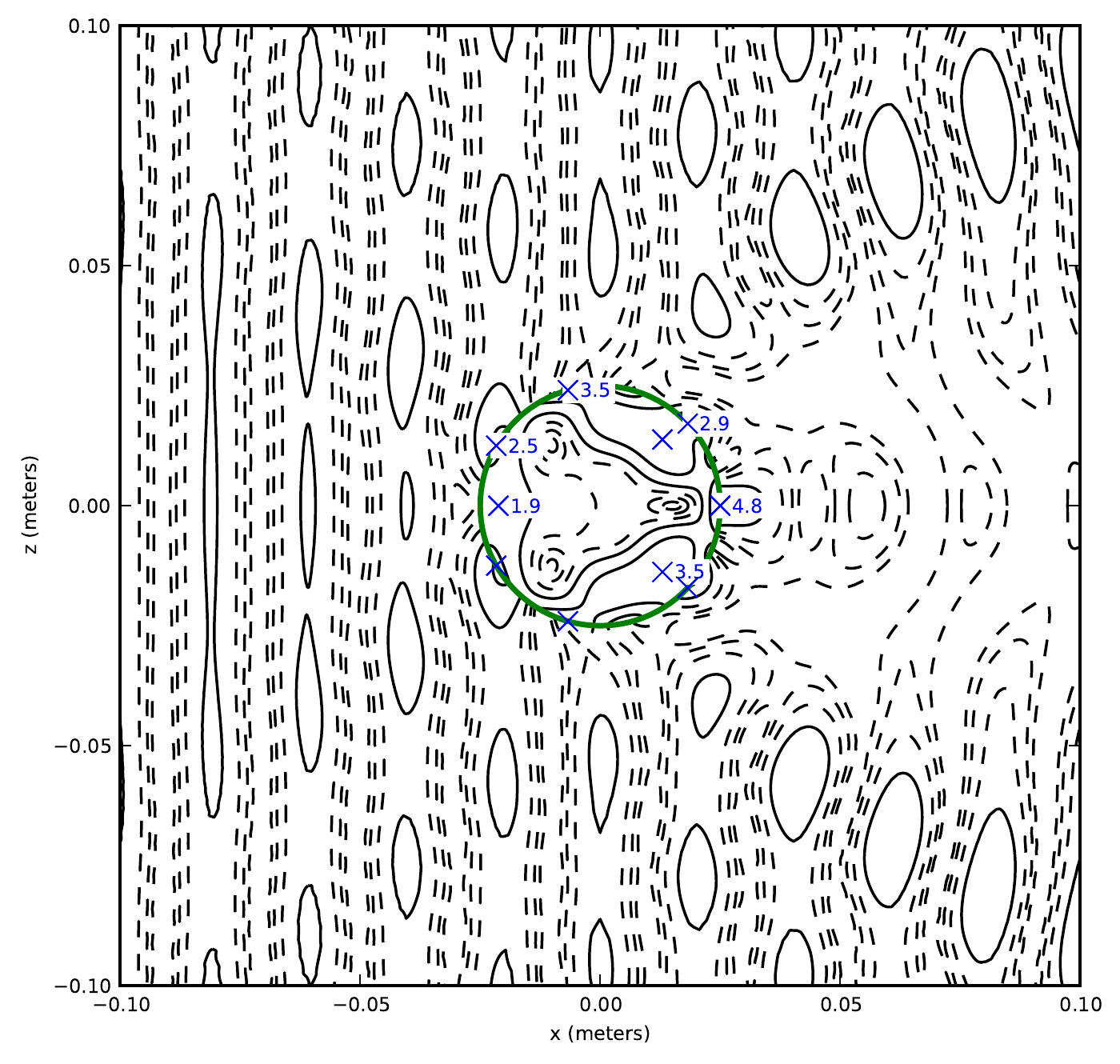}
    \caption{Viscous, $\eta_d = 0.5$, 38.35 kHz}
  \end{center}
\end{figure}

\begin{figure}[ht]
  \begin{center}
    \includegraphics[width=\textwidth]{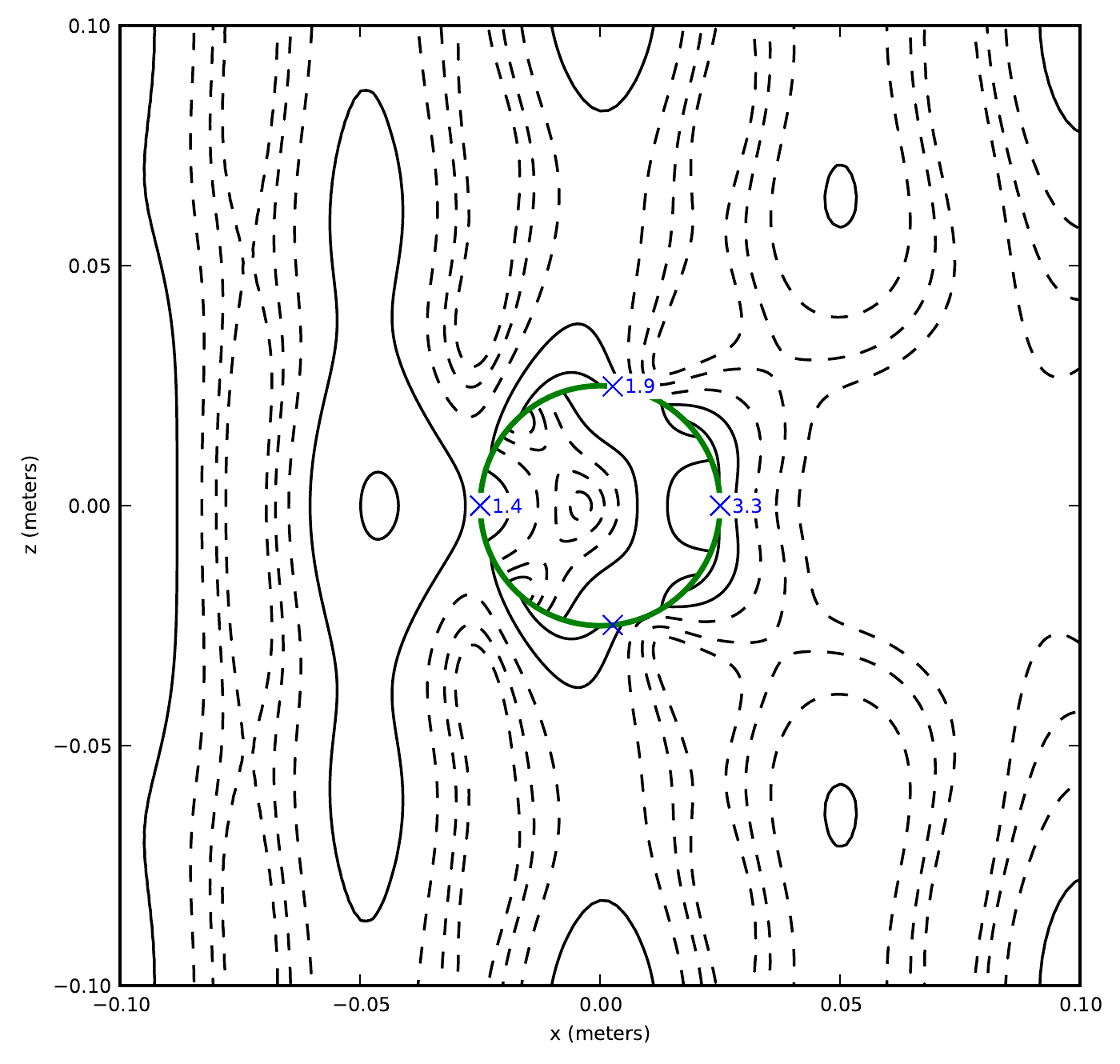}
    \caption{Viscous, $\eta_d = 1$, 15.80 kHz}
  \end{center}
\end{figure}

\begin{figure}[ht]
  \begin{center}
    \includegraphics[width=\textwidth]{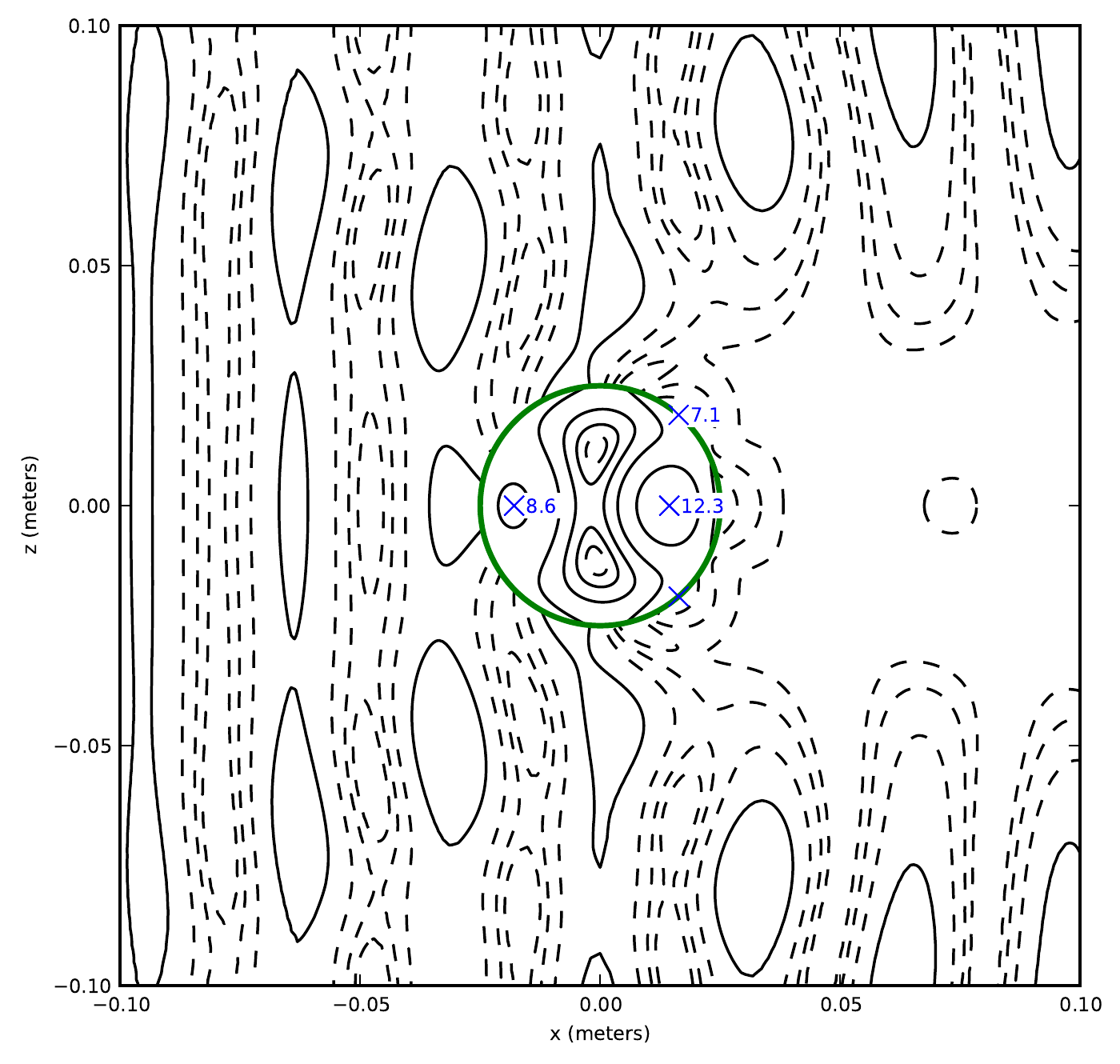}
    \caption{Viscous, $\eta_d = 1$, 24.54 kHz}
  \end{center}
\end{figure}

\begin{figure}[ht]
  \begin{center}
    \includegraphics[width=\textwidth]{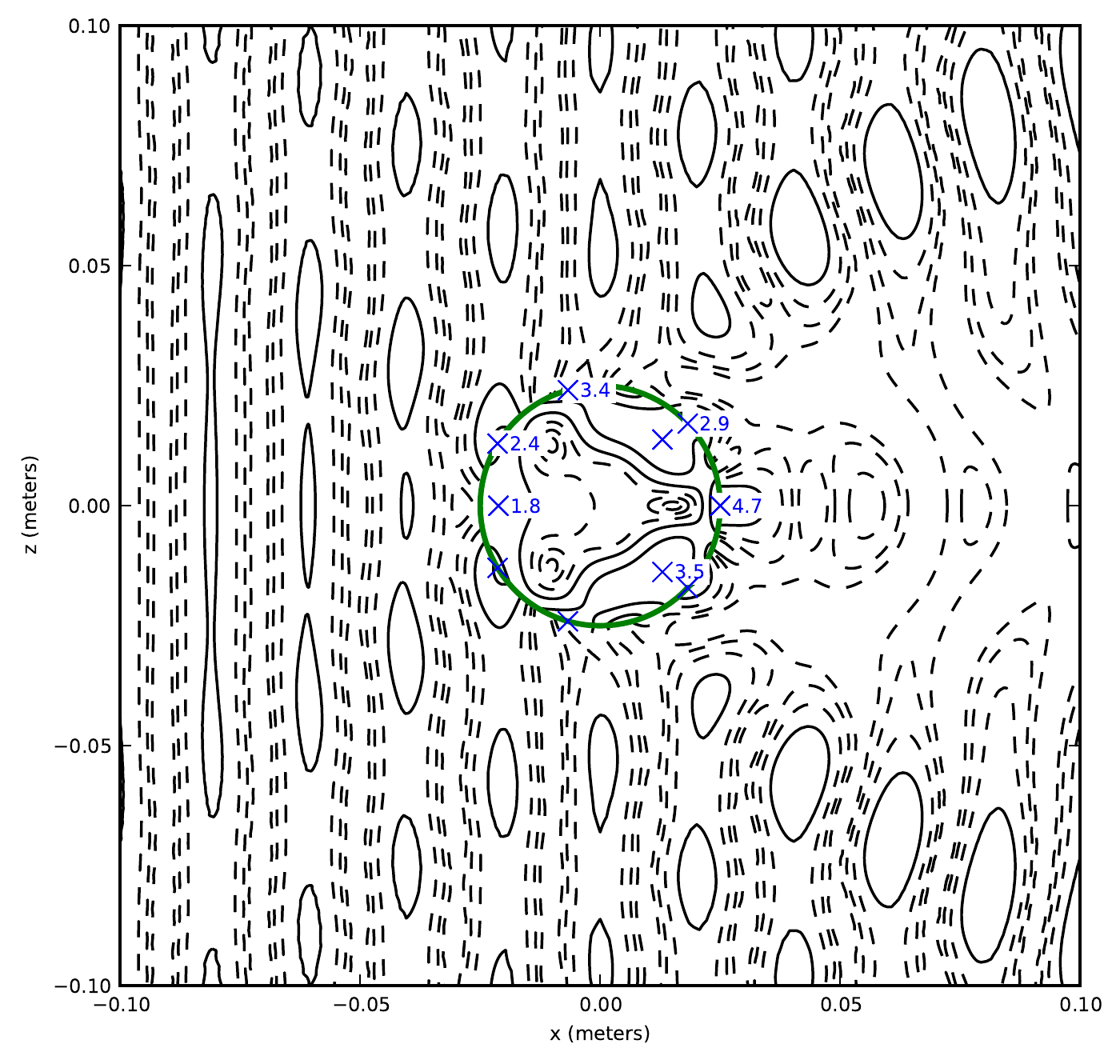}
    \caption{Viscous, $\eta_d = 1$, 38.35 kHz}
  \end{center}
\end{figure}

\end{document}